\RequirePackage[l2tabu, orthodox]{nag}

\documentclass[12pt]{amsart}
\usepackage{etoolbox}
\usepackage{blkarray}

\usepackage{fullpage,url,amssymb,enumerate,colonequals}
\usepackage[all]{xy} %
\usepackage{mathrsfs} %
\usepackage[dvipsnames]{xcolor}
\usepackage{tikz-cd}
\usepackage{breqn}
\usepackage[normalem]{ulem} %
\usepackage{textcomp}
\usepackage{microtype}
\usepackage{relsize}   %

\usepackage[
  left=25mm,
  right=25mm,
  top=20mm,
   bottom=20mm
]{geometry} %

\usepackage[
pdfauthor={Emre Can Sertöz},
pdftitle={Computing transcendence and linear relations of 1-periods},
hidelinks
]{hyperref}

\makeatletter
\renewcommand\subsection{\@startsection{subsection}{2}{\z@}
                                   {-1.5ex \@plus -.5ex \@minus -.2ex}%
                                   {1.0ex \@plus.2ex}%
                                   {\normalfont\bfseries}}%
\makeatother
\makeatletter
\g@addto@macro\bfseries{\boldmath}
\makeatother

\newcounter{nodenum}[section] 
\renewcommand{\thenodenum}{\thesubsection.\arabic{nodenum}} 
\newcommand{\node}[1]{\par\refstepcounter{nodenum}\vspace{1ex}\noindent(\thenodenum)\space \textbf{#1}}
\newcommand{\itnode}[1]{\par\vspace{1ex}\noindent\textit{#1}}
\makeatletter
 
\@addtoreset{equation}{nodenum} %
\makeatother

\usepackage[style=ieee-alphabetic,backend=biber, isbn=false, opcittracker=true,
  url=false
]{biblatex}
\AtEveryBibitem{%
  \clearfield{month}%
}

\usepackage{color}
\usepackage{marvosym} %

\usepackage[OT2,T1]{fontenc}
\DeclareSymbolFont{cyrletters}{OT2}{wncyr}{m}{n}
\DeclareMathSymbol{\Sha}{\mathalpha}{cyrletters}{"58}

\newcommand{\A}{\mathbb{A}}
\newcommand{\C}{\mathbb{C}}

\newcommand{\G}{\mathbb{G}}
\renewcommand{\H}{\operatorname{H}}

\newcommand{\Pp}{\mathbb{P}}
\newcommand{\p}{\mathbb{P}^1}

\newcommand{\Q}{\mathbb{Q}}

\newcommand{\R}{\mathbb{R}}
\newcommand{\N}{\mathbb{N}}
\newcommand{\Z}{\mathbb{Z}}

\newcommand{\Qbar}{{\overline{\Q}}}

\newcommand{\kk}{{\sf k}}

\newcommand{\xc}{\mathfrak{c}}

\newcommand{\xm}{\mathfrak{m}}

\newcommand{\xp}{\mathfrak{p}}

\newcommand{\calC}{\mathcal{C}}

\newcommand{\calF}{\mathcal{F}}

\makeatletter
\newcommand\RedeclareMathOperator{%
  \@ifstar{\def\rmo@s{m}\rmo@redeclare}{\def\rmo@s{o}\rmo@redeclare}%
}
\newcommand\rmo@redeclare[2]{%
  \begingroup \escapechar\m@ne\xdef\@gtempa{{\string#1}}\endgroup
  \expandafter\@ifundefined\@gtempa
     {\@latex@error{\noexpand#1undefined}\@ehc}%
     \relax
  \expandafter\rmo@declmathop\rmo@s{#1}{#2}}
\newcommand\rmo@declmathop[3]{%
  \DeclareRobustCommand{#2}{\qopname\newmcodes@#1{#3}}%
}
\@onlypreamble\RedeclareMathOperator
\makeatother

\DeclareMathOperator{\coker}{coker}

\DeclareMathOperator{\divv}{div}
\DeclareMathOperator{\Div}{Div}
\DeclareMathOperator{\End}{End}

\DeclareMathOperator{\Hom}{Hom}
\RedeclareMathOperator{\hom}{Hom}

\DeclareMathOperator{\im}{im}
\DeclareMathOperator{\image}{im}

\DeclareMathOperator{\Pic}{Pic}

\DeclareMathOperator{\Res}{Res}
\DeclareMathOperator{\rk}{rk}

\DeclareMathOperator{\Spec}{Spec}

\DeclareMathOperator{\supp}{supp}

\DeclareMathOperator{\tr}{tr}

\DeclareMathOperator{\Nm}{Nm}
\DeclareMathOperator{\AJ}{AJ}
\DeclareMathOperator{\poles}{poles}

\DeclareMathOperator{\Mat}{Mat}

\DeclareMathOperator{\rel}{rel}
\DeclareMathOperator{\ses}{ses}
\DeclareMathOperator{\triv}{triv}
\DeclareMathOperator{\Gy}{Gy}
\DeclareMathOperator{\Corr}{Corr}

\newcommand{\Set}{\text{(Set)}}
\newcommand{\motz}{1\text{-Mot}_\Z}
\newcommand{\mot}{1\text{-Mot}}
\newcommand{\mhsz}{(\Z,\Qbar)\text{-MHS}}
\newcommand{\mhs}{(\Q,\Qbar)\text{-MHS}}
\newcommand{\JJ}{J_{C,E}^D}
\newcommand{\ppJ}{J_{C,\psi}^\chi}

\newcommand{\AdR}{{\operatorname{AdR}}}
\newcommand{\adr}{\AdR}
\newcommand{\B}{{\operatorname{B}}}

\newcommand{\red}{{\operatorname{red}}}

\newcommand{\sat}{{\operatorname{sat}}}
\renewcommand{\ss}{\operatorname{ss}}
\newcommand{\oI}{\overline{I}}

\newcommand{\dd}{\mathrm{d}}

\usepackage{mathtools}

\newcommand{\abs}[1]{\lvert #1 \rvert}

\newcommand{\set}[1]{\{ #1 \}}

\newcommand{\tpi}{2\pi{\tt i}}
\newcommand{\tpip}{(2\pi{\tt i})}

\DeclareMathOperator{\res}{res}
\DeclareMathOperator{\id}{id}

\newcommand{\cc}{\mathcal{C}}

\newcommand{\ce}{\mathcal{E}}
\newcommand{\cf}{\mathcal{F}}

\newcommand{\co}{\mathcal{O}}
\newcommand{\cq}{\mathcal{Q}}
\newcommand{\cp}{\mathcal{P}}
\newcommand{\cR}{\mathcal{R}}

\newcommand{\cu}{\mathcal{U}}

\newcommand{\too}{\longrightarrow}
\newcommand{\tos}{\twoheadrightarrow}
\newcommand{\toi}{\hookrightarrow}
\newcommand{\isoto}{\overset{\sim}{\to}}

\newcommand{\coleq}{\colonequals}

\newcommand{\medoplus}{\mathsmaller{\bigoplus}} %
\newcommand{\wGamma}{\widetilde\Gamma}
\newcommand{\comp}{\xc_{\B,\AdR}} %

\newcounter{arrowcounter}

\newcommand{\circled}[1]{%
  \tikz[baseline=(char.base)]{%
    \node[draw, circle, inner sep=1pt] (char){\ifmmode #1\else\scriptsize #1\fi};%
  }%
}

\newcommand{\labelarr}[1]{%
  \refstepcounter{arrowcounter}%
  \label{#1}%
  \circled{\thearrowcounter}%
}

\newcommand{\circledref}[1]{%
  {\begingroup
     \hypersetup{hidelinks}%
     \hyperref[#1]{\circled{\ref*{#1}}}%
   \endgroup}%
}

 \author[E.~C.~Sertöz]{Emre Can Sertöz}
 \address{Emre Can Sertöz, Mathematical Institute, Leiden University, The Netherlands}
 \email{emre@sertoz.com}
 \urladdr{\url{https://emresertoz.com}}
 
 \author[J.~Ouaknine]{Joël Ouaknine}
 \address{Joël Ouaknine, Max Planck Institute for Software Systems, Saarland, Germany}
 \email{joel@mpi-sws.org}
 \urladdr{\url{https://people.mpi-sws.org/~joel/}}
 
 \author[J.~Worrell]{James Worrell}
 \address{James Worrell, Department of Computer Science, Oxford University, United Kingdom}
 \email{jbw@cs.ox.ac.uk}
 \urladdr{\url{https://www.cs.ox.ac.uk/people/james.worrell/home.html}}

\title{Computing transcendence and linear relations of 1-periods}

\begin{document}

\begin{abstract}
A $1$-period is a complex number given by the integral of a univariate algebraic function, where all data involved --- the integrand and the domain of integration --- are defined over algebraic numbers. We give an algorithm that, given a finite collection of $1$-periods, computes the space of $\Qbar$-linear relations among them. In particular, the algorithm decides whether a given $1$-period is transcendental, and whether two $1$-periods are equal. This resolves, in the case of $1$-periods, a problem posed by Kontsevich and Zagier, asking for an algorithm to decide equality of periods. The algorithm builds on the work of Huber and Wüstholz, who showed that all linear relations among $1$-periods arise from $1$-motives; we make this perspective effective by reducing the problem to divisor arithmetic on curves and providing the theoretical foundations for a practical and fully explicit algorithm. To illustrate the broader applicability of our methods, we also give an algorithmic classification of autonomous first-order (non-linear) differential equations.
\end{abstract}

\maketitle
\tableofcontents

\newpage 
\section{Introduction}

\node{} Periods are the values of integrals of algebraic forms over semi-algebraic domains, both defined over algebraic numbers. These numbers appear broadly in the sciences and connect arithmetic algebraic geometry to transcendental number theory. The Kontsevich--Zagier~\cite{KZ} and Grothendieck~\cite{Grothendieck1966,Andre1989,Ayoub2014} period conjectures predict that equality between periods should be decidable by algebro-geometric means. Yet these conjectures remain far out of reach in this generality.

\node{} Kontsevich and Zagier ask for an algorithm~\cite[Problem 1, p.7]{KZ} to decide equality of two periods. In this paper, we resolve this problem in the case of 1-periods: integrals of \emph{univariate} algebraic forms. We construct an algorithm that can decide equality among 1-periods and even decide transcendence. Beyond mere decidability statements, we have devised an algorithm grounded in divisor arithmetic on curves allowing us to rely entirely on well-established, practical methods.

\node{} The theoretical bedrock for our algorithm is the work of Huber and Wüstholz~\cite{HW}, who proved that all $\Qbar$-linear relations among 1-periods arise from morphisms of 1-motives. While their results are not algorithmic, the algebro-geometric nature of the statement suggests that an effective version should be possible.  Rather than attempting to translate their methods --- which would require working with semi-abelian varieties embedded in high-dimensional projective spaces --- we focus instead on generalised Jacobians and the mixed Hodge structures $\H_1(C\setminus{}D,E)$ of punctured and marked curves. In doing so, we develop a range of computational techniques that may be of independent interest for the study of algebraic curves and their periods; see for instance~\eqref{node:intro_classify}.

\node{} To give an emblematic example of a $1$-period, consider the identity
\begin{equation}
\pi = \int_{-1}^1 2\sqrt{1 - x^2} \, \dd x,
\end{equation}
which expresses the area of the unit circle as the integral of an algebraic function over an interval with algebraic endpoints. 

\node{} In general, a $1$-period may be defined as an integral $\int_a^b f(x)\, \dd x$, where $a,b$ are real algebraic numbers and the integrand $f(x)\colon [a,b] \to \C$ is continuous and algebraic over~$\Qbar$; that is, it satisfies a nonzero polynomial relation $P(x, f(x)) = 0$ for some $P \in \Qbar[x,y]$. This definition recognisably includes classical examples such as elliptic integrals and special values of the Gauss hypergeometric function.

\node{} While this formulation offers a concrete and familiar perspective, it is not the one we will adopt. For our purposes, a more flexible definition --- better suited to the geometric and computational framework developed here --- will be preferable. The two formulations are equivalent, though we will not need this fact.

\node{}\label{def:period} A \emph{$1$-period} is a complex number that can be realised as the value of an integral
$
\int_\gamma \omega,
$
where $C/\Qbar$ is a smooth projective curve, $\omega$ is a rational differential on $C$, and $\gamma$ is a smooth $1$-chain on $C(\C) \setminus \poles(\omega)$ with boundary $\partial \gamma$ supported on $C(\Qbar)$. We refer to the data $(\gamma, \omega)_C$ as a \emph{representation} of the $1$-period. 

\node{} To work with $1$-periods effectively, we must specify how the tuple $(\gamma, \omega)_C$ is represented. The curve $C/\Qbar$ is given by specifying a birational plane model for each irreducible component; equivalently, each component is represented by an irreducible polynomial $P(x,y)\in\Qbar[x,y]$. The rational differential $\omega$ is likewise represented, per component, by a polynomial $Q(x,y)\in\Qbar[x,y]$ such that the form is given by the restriction of $Q(x,y)\,\dd x$ to~$C$. We turn next to the representation of the chain $\gamma$.

\node{} For the smooth $1$-chain $\gamma$, the value of the integral depends only on its homology class. Any representation of $\gamma$ from which this class can be recovered is therefore acceptable --- for instance, one that allows the integrals of rational forms over $\gamma$ to be approximated to arbitrary precision~\eqref{node:hom_reduction}. For the basis of homology, we use $1$-chains that are rectilinear in the $x$-coordinate; see~\S\ref{sec:betti}.

\node{} We fix an embedding $\Qbar \subset \C$, and assume that algebraic numbers are represented by their minimal polynomial over~$\Q$ together with data isolating the intended complex root (e.g., a complex ball with rational center and radius).

\node{Main Theorem.} Given representations $(\gamma_1, \omega_1)_{C_1}, \dots, (\gamma_k, \omega_k)_{C_k}$ of $1$-periods $\alpha_i = \int_{\gamma_i} \omega_i$, Algorithm~\eqref{alg:main_algorithm} in Section~\ref{sec:algo_outline} computes a $\Qbar$-basis for the space of $\Qbar$-linear relations
\begin{equation}
\rel_{\Qbar}(\alpha_1, \dots, \alpha_k) \coleq \left\{ (\beta_1, \dots, \beta_k) \in \Qbar^k : \sum_{i=1}^k \beta_i \alpha_i = 0 \right\}.
\end{equation}

\node{Corollary.} Given a representation of a $1$-period $\alpha$, Algorithm~\eqref{alg:main_algorithm} decides whether $\alpha$ is transcendental, and when it is algebraic, produces a standard representation of $\alpha$ as an algebraic number.

\itnode{Proof.} Indeed, $1 = \int_{[0,1]} \dd x$ is a $1$-period, so we may compute $\rel_{\Qbar}(\alpha, 1)$, which is the $0$ vector space precisely when $\alpha \notin \Qbar$. (The basis for the $0$ vector space is the empty set.) If $\alpha \in \Qbar$, then $\rel_{\Qbar}(\alpha, 1)$ is a $1$-dimensional subspace of $\Qbar^2$, and the algorithm will return a generator $(\beta_1,\beta_2)$ such that $\alpha = -\beta_2/\beta_1$. \hfill \qed

\node{Corollary.}\label{cor:equality} Given representations of two $1$-periods $\alpha_1$ and $\alpha_2$, Algorithm~\eqref{alg:main_algorithm} decides whether $\alpha_1 = \alpha_2$.

\itnode{Proof.} Equality holds precisely when $(1, -1) \in \rel_{\Qbar}(\alpha_1, \alpha_2)$, which can be tested using a $\Qbar$-basis. \hfill \qed

\node{Corollary.} The $\Qbar$-vector space of $1$-periods is effective.

\itnode{Proof.} This means that the operations of addition, scalar multiplication, and equality are computable. The $\Qbar$-scaling of a representative is achieved by scaling the form: $\beta \cdot (\gamma,\omega)_C = (\gamma,\beta \cdot \omega)_C$. To add two representatives $(\gamma_1,\omega_1)_{C_1}$ and $(\gamma_2,\omega_2)_{C_2}$, take the disjoint union of the underlying curves, and add the $1$-chains and $1$-forms there. Equality is checked as in Corollary~\eqref{cor:equality}.

\node{Technical framework.} We construct a complete computational model for the mixed Hodge structure $\H_1(C \setminus D, E)$: the de~Rham realisation via rational differentials over the base field~\S\ref{sec:adr}; the Betti realisation via embedded graphs~\S\ref{sec:betti}; and the comparison isomorphism via certified precision integration. To circumvent arithmetic on periods, we augment this structure with the Jacobian motive of $(C \setminus D, E)$, see~\S\ref{sec:pushpull}. This is a toric extension $J_C^D$ of the Jacobian $J_C$ of~$C$, marked by divisors supported on~$E$. Divisor arithmetic on~$C$ enables exact computation with the points of $J_C^D$, which compensates for the inability to test identities among periods directly. 

To obtain the flexibility needed for ``symmetrizing'' the motive while remaining within the domain of curves, as required in~\S\ref{sec:supsat}, we generalise slightly to push-pull Jacobian motives and their mixed Hodge structures. On the theoretical side, we give an effective description of the period relations of ``sufficiently symmetric'' $1$-motives; see~\S\ref{sec:motives}. These developments require a full-stack construction: even classical components must be extended to accommodate the generality of $(C \setminus D, E)$, and functoriality and effectivity must be engineered to work together as a computable theory of $1$-periods.

\node{Classification problem.}\label{node:intro_classify} We expect the computational model for $\H_1(C\setminus{}D,E)$, together with the subroutines developed, will be of independent interest. To illustrate the point, we resolve the problem of algorithmically classifying first-order autonomous (non-linear) differential equations posed by Noordman, van der Put, and Top~\cite[p.1655]{NPT22}. See~\S\ref{sec:autonomous_ode} for the precise statement and the algorithm.

\subsection{Connections to foundational work}

\node{} Another approach to deciding equality of $1$-periods is to determine effective separation bounds. Hirata-Kohno~\cite{HirataKohno1991} and Gaudron~\cite{Gaudron2005} prove that for any given $1$-period~$\alpha$, there exists a constant $c>0$ --- effective in terms of the data defining~$\alpha$ --- such that $|\alpha| < c$ implies $\alpha = 0$. These results ultimately rely on Wüstholz's analytic subgroup theorem~\cite{Wustholz1989} and its quantitative refinements by Masser–Wüstholz~\cite{MW93a,MW93b}. In the classical setting of Baker periods (i.e., underlying curves have genus zero), the bound is completely explicit~\cite{BW93}. For periods of elliptic curves, explicit bounds have been obtained~\cite{David1995,DHK09}, and more recently, also in the case of abelian varieties~\cite{BosserGaudron19}, corresponding to the situation in~\eqref{def:period} where the $1$-form has no poles. In general, however, the constant $c$ remains only effective in principle: computing it explicitly would require navigating a long chain of intermediate constants as well as determening subtle arithmetic constants; see~\cite[p.145]{Gaudron2005} and the discussion after~\cite[Theorem 7.2]{BW08} for more on this point. %

\node{} Separation bounds permit, in principle, effective tests for equality of $1$-periods. But deciding transcendence requires exhaustive knowledge of all $\Qbar$-linear relations. In this sense, our algorithm provides the first effective method for determining transcendence of $1$-periods.

\node{} Our algorithm ultimately relies on one of the main theorems of Huber and Wüstholz --- the ``Dimension estimate''~\cite[Theorem~15.3]{HW} --- or rather on a reassembly of its proof, which gives a description of the space of period relations associated to a \emph{saturated} $1$-motive. Their arguments, which build around Wüstholz's analytic subgroup theorem~\cite{Wustholz1989}, determine the $\Qbar$-dimension for the space of $1$-periods of a saturated $1$-motive. Explicitly relating the period relations of an arbitrary motive to those of a saturated motive involves substantial additional work. An important part of our contribution is an effective method to perform the saturation step, see~\S\ref{sec:supsat}, using only arithmetic on algebraic curves --- a task that turns out to be technically intricate and requires the full machinery developed here. 

\node{} Recall that one of our key contributions is to base our computations on the divisor arithmetic on curves. To put this into context, we compare it to Chapter~14 of~\cite{HW}, where the vanishing of periods (the case $k = 1$ for our algorithm) is studied from a curve-theoretic point of view. There, it is assumed that the underlying curve has a \emph{simple Jacobian}, to eliminate the need for supersaturation --- the most challenging part of our algorithm, see~\S\ref{sec:supsat}. Nevertheless, a full translation to arithmetic on curves is not attained there: even in the first nontrivial case, one ends up with a check involving the condition~``(14.3)'' of~\cite[p.137]{HW} involving the vanishing of a differential under a motivic pullback map. We are able to replace motivic abstractions with explicit, curve-theoretic computations by methodically building up subroutines to compute with the $\H_1$ of punctured marked curves and building correspondences as the cornerstone of our morphisms. To demonstrate the effect, we unwind what our algorithm does in the aforementioned case~\eqref{node:reinterpret_143}.

\node{} A crucial computational input to our algorithm is the endomorphism algebra of a Jacobian, represented explicitly in terms of algebraic correspondences. This plays a central role in the supersaturation step, where correspondences are used to generate additional relations. For this, we rely on the practical semi-algorithm of Costa--Mascot--Sijsling--Voight~\cite{CMSV18}, which computes the endomorphism algebra as a space of correspondences and has been successfully implemented. Its correctness is unconditional, and its termination is guaranteed under the Mumford--Tate conjecture, as shown by Costa--Lombardo--Voight~\cite{CLV21}. To ensure unconditional termination in theory, we cite the effective construction of Lombardo~\cite{Lombardo2018}, though it is not intended for practical use. See~\S\ref{sec:corr} for more details. %

\node{} We rely on many other subroutines, often implicitly. However, the computation of the Riemann--Roch space
\begin{equation}
L(D) = \{f \in \kappa(C)^\times \mid D + \divv(f) \ge 0 \} \cup \{0\} \subset \kappa(C)
\end{equation}
for a divisor $D$ on an irreducible curve $C/\Qbar$ with function field $\kappa(C)$ lies at the heart of our approach. See~\cite{Hess2002} for an efficient algorithm to compute $L(D)$ and a historical overview of its development.

\subsection{Acknowledgements}

We thank Jaap Top for pointing us to his paper on autonomous differential equations, Don Zagier for valuable remarks that clarified our exposition, and Edgar Costa, Annette Huber, and Davide Lombardo for fruitful discussions. We are also grateful to our colleagues who provided feedback on earlier versions of this paper, including Steffen Müller, Riccardo Pengo, Matthias Schütt, and others.

ECS would like to thank the Institute for Algebraic Geometry at Leibniz University Hannover, where this project began. JO is also affiliated with Keble College, Oxford as Emmy Network Fellow, and supported by ERC grant DynAMiCs (101167561) and DFG grant 389792660 as part of TRR 248 (\url{https://perspicuous-computing.science}). JW is supported by UKRI Frontier Research Grant EP/X033813/1.

\newpage
\section{The algorithm}\label{sec:algo_outline}

\node{Algorithm.}\label{alg:main_algorithm} 
The following is the main algorithm of the paper, which computes all $\Qbar$-linear relations between a given finite collection of $1$-periods. Intermediate steps refer to constructions developed in the rest of the paper. 

\node{Input.} A tuple of representatives, $(\gamma_i,\omega_i)_{C_i}$, $i=1,\dots,k$, for $1$-periods $\alpha_1,\dots,\alpha_k$. 

\node{Output.} A basis for the space of $\Qbar$-linear relations between $\alpha_1,\dots,\alpha_k$, i.e., the kernel of the map $\Qbar^k \to \C : (\beta_1,\dots,\beta_k) \mapsto \sum_{i=1}^k \alpha_i \beta_i$

\node{} For each $i=1,\dots,k$, compute a finite set $D_i \subset C_i(\Qbar)$ supporting the residual divisor of $\omega_i$ and compute a finite set $E_i \subset C_i(\Qbar)$ supporting the boundary divisor of $\gamma_i$, such that $D_i$ and $E_i$ are disjoint. Let $C = \coprod_{i=1}^k C_i$, $D = \coprod_{i=1}^k D_i$, $E = \coprod_{i=1}^k E_i$ where we view $D, E \subset C(\Qbar)$.

\node{} Compute~\eqref{node:compute_rep_of_MHS} a representative of the mixed Hodge structure 
\begin{equation}
  \H_1(C\setminus{}D,E) = (\H_1^\B(C\setminus{}D,E)), \H^1_\AdR(C\setminus{}D,E)^\vee,\comp). 
\end{equation}
That means, we compute a basis of differentials for the algebraic de Rham cohomology $\H^1_\AdR(C\setminus{}D,E) \simeq \Qbar^m$, a basis of $1$-chains for the Betti homology $\H_1^\B(C\setminus{}D,E) \simeq \Q^m$, and the representatives allow for an arbitrary precision approximation of the comparison isomorphism. 

\node{} Let $\gamma = \gamma_1 + \dots + \gamma_k$ be the $1$-chain on $C$ obtained by summing the given $1$-chains. Determine~\eqref{node:hom_reduction} the coordinates of the homology class $[\gamma]$ in $\H_1^\B(C\setminus{}D,E)$ with respect to our basis of $1$-chains. View $\omega_i$ as a $1$-form on $C$ extended by $0$ to components $C_j$, $j \neq i$, and compute~\eqref{node:symbolic_reduction} the coordinates of the cohomology class $[\omega_i]$ in $\H^1_\AdR(C\setminus{}D,E)$. We now have the matrix representation of the following map:
\begin{equation}
 I \colon \Qbar^k \to \H_1^\B(C\setminus{}D,E) \otimes \H^1_\AdR(C\setminus{}D,E) \simeq \Qbar^{m \times m} : \beta \mapsto \sum_{i=1}^k \beta_i [\gamma]\otimes[\omega_i].
\end{equation}

\node{}\label{node:compute_jacobian_in_algo} Compute~\eqref{thm:endo} a representation of the endomorphism algebra $\ce = \End(J_C)$ of the Jacobian of $C$.

\node{}\label{node:supersaturate_in_algo} Let $M=J_{C,E}^D$ denote the Jacobian motive of $(C\setminus{}D,E)$. Compute~\eqref{node:supsat}, the supersaturation $M^{ss}$ of $M$ together with an explicit isogeny direct sum decomposition of $M^{ss}$ into a Baker motive $M_B$ and a saturated push-pull Jacobian $M_1 = J_{C,\psi}^\chi$. The process is period effective, meaning that the period relations of $M$ can be determined from those of $M^{ss}$, see~\eqref{node:rels_of_M_from_Mss}. The motive $M_1$ is computed together with the explicit $\ce$-action $\ce \otimes \H_1(M_1) \to \H_1(M_1)$.

\node{} Since $M_1$ is saturated, we can compute the period relations $\cR(M_1) \subset \H_1^\B(M_1)\otimes \H^1_{\AdR}(M_1)$ of $M_1$: this space is equal~\eqref{thm:expected_relations} to the expected period relations $\cR_e(M_1)$, which we can determined from the $\ce$-action on $\H_1$ above and from the trivial relations which are readily computed~\eqref{rem:trivial}. 

\node{} Compute~\eqref{lem:baker_period_relations} the period relations $\cR(M_B)$ of the Baker motive $M_B$.

\node{} We now compute $\cR(M^{ss})$, the space of period relations of the supersaturated motive. This space decomposes as follows~\eqref{cor:rels_of_Baker_and_sat}, and each summand is now computable:
\begin{equation}
  \cR(M^{ss}) = \cR_{\triv}(M^{ss}) + \cR(M_B) \oplus \cR_e(M_1) + \H_1^\B(M_B)\otimes\H^1_\AdR(M_1) + \H_1^\B(M_1)\otimes\H^1_\AdR(M_B). 
\end{equation}

\node{} Since the construction of $M^{ss}$ from $M$ was period effective~\eqref{node:rels_of_M_from_Mss}, we can determine $\cR(M) \subset \H_1^\B(M) \otimes \H^1_\AdR(M)$ by a simple procedure from that of $\cR(M^{ss})$. Compute $\cR(M)$.

\node{} From the very start we work with the identification $\H_1(J_{C,E}^D) = \H_1(C\setminus{}D,E)$. Therefore $\cR(M)$ lives in the codomain of $I$. Output a basis for the pullback $I^{-1}\cR(M) \subset \Qbar^k$.

\node{Remark.} We suspect that the most time consuming step will be~\eqref{node:compute_jacobian_in_algo}, but we can rely on the streamlined implementation in~\cite{CMSV18}. The step~\eqref{node:supersaturate_in_algo} is terribly involved and could present another computational bottleneck.

\newpage
\section{Relative algebraic de Rham cohomology of a punctured curve}\label{sec:adr}

\node{} Let $C$ be a smooth, proper curve defined over $\Qbar$, and let $D,E \subset C(\Qbar)$ be two finite, disjoint sets. The pair $(C\setminus{}D,E)$ will be called a \emph{punctured marked curve}. Topologically, this is to be viewed as $C\setminus{}D$ relative to $E$. 

\node{} This section describes how to compute and effectively work with the algebraic de~Rham cohomology $\H^1_\AdR(C\setminus{}D,E)$ of the punctured marked curve $(C\setminus{}D, E)$. Since (co)homology is additive over components of $C$, we will assume without loss of generality that $C$ is irreducible, unless stated otherwise.

\node{} The case $D=E=\emptyset$ is well-known to the experts~\cite{Katz1968a,Coleman1989}. Nevertheless, we spell-out the proof of the fact that $\H^1_\AdR(C)$ is isomorphic (as a $\Qbar$-vector space) to the space of second kind differentials on $C$ modulo exact forms (over $\Qbar$). The pairing with singular homology is via integration~\eqref{node:period_pairing}. This is classical over $\C$ and applicable in greater generality~\cite{Rosenlicht1953,AH55,Messing1975}. This allows for the practical identity
\begin{equation}\label{eq:adr_sec_intro}
  \H^1_\AdR(C) \simeq \H^0(C,\Omega_C^1((g+1)p))
\end{equation}
where $p \in C(\Qbar)$ is any non-Weierstrass point and $g$ is the genus of $C$, see~\eqref{prop:eff_C_coh}. We also give two different ``reduction algorithms'' to find the unique representative in $\H^0(C,\Omega_C^1((g+1)p))$ of any given second type differential (see~\S\ref{sec:reduction_algos}). 

\node{} For the general case when $D,E$ are non-empty, we could not find a reference to a result analogous to~\eqref{eq:adr_sec_intro}. We therefore prove Proposition~\eqref{prop:coh_of_punctured_curve} and Corollary~\eqref{cor:represent_cohom} which allows us to effectively compute a basis for $\H^1_\AdR(C\setminus{}D,E)$~\eqref{prop:effective_cohom_basis}. The aforementioned reduction algorithms in~\S\ref{sec:reduction_algos} are applicable in this generality.

\subsection{Representing irreducible curves and their differentials}\label{sec:rep_irred_curve} 

\node{} In this subsection, $C/\Qbar$ is an irreducible, smooth, proper curve.

\node{} A representation of the curve $C$ is a (reduced, irreducible) polynomial $P(x,y) \in \Qbar[x,y]$ defining an affine plane curve $Z(P) \subset \A^2$ birational to $C$. Without loss of generality, we will assume that the $y$-derivative $P_y$ of $P$ is non-zero. By default, our curves come with a finite map $x \colon C \to \p$, the projection to the $x$-axis.

\node{}\label{node:rep_curve} Giving a plane model of $C$ is equivalent to representing the function field $\kappa(C)$ of $C$ by exhibiting an element $x \in \kappa(C)$, transcendental over $\Qbar$, and an element $y \in \kappa(C)$ algebraic over $\Qbar(x)$ and generating $\kappa(C)$. Thus, $\kappa(C) = \Qbar(x)[y]/P(x,y)$. 

\node{} The transcendental element $x\in \kappa(C)$ realizes $C$ as a finite cover $x\colon C \to \p$ of degree equal to the $y$-degree of $P$. The corresponding map $x \colon C(\C) \to \p(\C)$ between the underlying Riemann surfaces realizes $C(\C)$ as a branced cover of $\p(\C)$. 

\node{} The rational differentials on $C$ are the Kähler differentials $\Omega^1_{\kappa(C)/\Qbar}$ of $\kappa(C)$ over $\Qbar$. The relation $\dd P = P_x \dd x + P_y \dd y  = 0$ and our assumption $P_y \neq 0$ allows us to eliminate $\dd y$,
\begin{equation}
  \Omega^1_{\kappa(C)/\Qbar} = \kappa(C) \dd x = \left(\Qbar(x)[y]/P(x,y)\right) \, \dd x,
\end{equation}

\node{} We will represent a rational differential $\omega \in \Omega^1_{\kappa(C)/\Qbar}$ by writing it as $\omega = f(x,y) \dd x$ where $f(x,y) \in \kappa(C) = \Qbar(x)[y]/P(x,y)$.

\subsection{Representing reducible curves and their differentials}\label{sec:rep_red_curve}

\node{} Assume $C/\Qbar$ is a smooth proper curve with irreducible components $C_1,\dots,C_k$. 

\node{} We will represent $C$ by representing the tuple $C_1,\dots,C_k$ of its components, i.e., by exhibiting a tuple of polynomials $P_i \in \Qbar[x,y]$, for $i=1,\dots,k$, each with non-vanishing $y$-derivatives. Let $x_i$ be the image of $x$ in $\kappa(C_i) = \Qbar(x)[y]/P_i(x,y)$.

\node{} We define the algebra $\kappa(C) \coleq \kappa(C_1) \oplus \dots \oplus \kappa(C_k)$ of rational functions on $C$. Let $x=(x_1,\dots,x_k) \in \kappa(C)$ be the tuple of transcendental generators. This is equivalent to identifying the codomains of the maps $x_i \colon C_i \to \p$ to form a finite map $x \colon C \to \p$.

\node{} The space of rational differentials $\Omega^1_{\kappa(C)/\Qbar}$ on $C$ is identified with $\kappa(C) \dd x$. Equivalently, a rational differential $\omega$ on $C$ is a tuple of rational differentials $(\omega_1,\dots,\omega_k)$ with $\omega_i \in \Omega^1_{\kappa(C_i)/\Qbar}$.

\subsection{Poles of forms and definition of algebraic de Rham cohomology} 

\node{} Throughout this section, $C$ is an \emph{irreducible} smooth proper curve over $\Qbar$. 

\node{} Let $\co_C$ be the sheaf of rational functions on $C$ and $\Omega^1_C$ the sheaf of Kähler differentials of $C/\Qbar$. For any divisor $\xi$ on $C$ and any open set $U \subset C$ we have $\co_C(\xi)(U) \subset \kappa(C)$ and $\Omega^1_C(\xi)(U) \subset \Omega^1_{\kappa(C)/\Qbar}$. In particular, we will view local regular sections as global rational sections. The $\Qbar$-dimension of the global sections of a sheaf on $C$ will be denoted by $h^0$.

\node{} For a rational differential $\omega$ on $C$, the \emph{polar locus} of $\omega$ is the set of poles of $\omega$. A pole of order one is called a \emph{simple pole}. The polar locus is a finite subset of $C(\Qbar)$ and is readily to computed.

\node{}\label{node:residue} For any point $p \in C(\Qbar)$ we will write $\res_p \omega \in \Qbar$ for the residue of $\omega$ at $p$. Note that the residue $\res_p \omega$ can be computed by making a Laurent series expansion of $\omega$ with respect to a local coordinate $z_p$ at $p$ and taking the coefficient of $\dd z_p / z_p$. If $\res_p \omega \neq 0$ then $p$ is a \emph{residual pole of $\omega$}.

\node{} We define the \emph{residual divisor of $\omega$} to be $\Res(\omega) \coleq \sum_{p \in C(\Qbar)} (\res_p \omega) p$. The support of $\Res(\omega)$ is the \emph{residual polar locus of $\omega$}. Since the polar locus is finite, the residual divisor is readily computed.

\node{} If $\omega$ has no poles, it is \emph{regular}; classically called a \emph{differential of the first kind}. If $\omega$ has no residual poles, $\Res(\omega)=0$, then $\omega$ is a \emph{differential of the second kind}. If all poles of $\omega$ are simple, then $\omega$ is a \emph{differential of the third kind}. Every rational differential is the sum of a differential of the second and third kind, and this decomposition is unique modulo differentials of the first kind.

\node{}\label{node:alternative_def_of_second_type} A rational differential $\omega$ is of second kind if and only if it is Zariski locally equivalent to a section of $\Omega^1_C$ modulo an exact form. The latter condition means: $\forall p \in C(\Qbar)$ there is an open neighborhood $U$ of $p$ and an exact form $\dd f$ such that $\omega + \dd f$ is regular on $U$. 

\node{} Let $D \subset C(\Qbar)$ be a finite set. If $\omega$ is a differential of the third kind with poles supported on $D$ then we will call $\omega$ a $\log(D)$-regular differential. For an open set $U \subset C$, we will say $\omega$ is $\log(D)$-regular on $U$ if all poles of $\omega|_U$ are simple and supported on $D\cap U$. Note that $\Omega^1_C(D)$ is the sheaf of $\log(D)$-regular differentials.

\node{Remark.} In the case of curves, the sheaf $\Omega^1_C(D)$ coincides with the sheaf $\Omega^1_C(\log(D))$ of logarithmic differentials at $D$. The latter is needed for algebraic de Rham cohomology in arbitrary dimensions but the former is standard when discussing linear systems on curves. We will stick to curves and use the former notation throughout.

\node{} Let $D, E \subset C(\Qbar)$ be finite disjoint sets. We recall (\cite{Grothendieck1966},~\cite[\S 5.5]{Peters2008},~\cite[\S 2.2.6]{BurgosFresan}) that the algebraic de Rham complex on $(C\setminus{}D,E)$ is %
\begin{equation}\label{eq:alg_dR_complex}
  \Omega_{(C\setminus{}D,E)}^{\bullet} : 0 \to \co_C  \to \Omega^1_C(D) \oplus \co_E  \to 0
\end{equation}
where the non-zero map is $f \mapsto (\dd f, -f|_E)$. The hyper cohomology of~\eqref{eq:alg_dR_complex} is the algebraic de Rham cohomology of $(C\setminus{}D,E)$. Our focus is on the first cohomology:
\begin{equation}
  \H^1_\AdR(C\setminus{}D,E) \coleq \H^1(C,\Omega_{(C\setminus{}D,E)}^\bullet).
\end{equation}

\subsection{Algebraic de Rham cohomology in terms of second type differentials}

\node{} Fix $D,E \subset C(\Qbar)$ finite and disjoint. 

\node{} A rational differential $\omega$ on $C$ will be called $\log(D)$-second kind if $\omega$ is Zariski locally equivalent to a section of $\Omega^1_C(D)$ modulo exact forms (compare with~\eqref{node:alternative_def_of_second_type}).

\node{Proposition.}\label{prop:coh_of_punctured_curve} We have a canonical isomorphism
\begin{equation}\label{eq:coh_of_punctured_curve}
  \H^1_\AdR(C\setminus{}D) \simeq \frac{\set{\text{$\log(D)$-second kind differentials on $C$}}}{\set{\text{exact forms on $C$}}},
\end{equation}
with the nontrivial parts of the Hodge and weight filtration given by
\begin{align}
  \label{eq:F1_puncture}  F^1\H^1_\AdR(C\setminus{}D) &= \H^0(C,\Omega^1_C(D)) = \set{\text{$\log(D)$-regular differentials on $C$}}, \\
  \label{eq:W1_puncture}  W_1\H^1_\AdR(C\setminus{}D) &= \H^1_\AdR(C) = \frac{\set{\text{second kind differentials on $C$}}}{\set{\text{exact forms on $C$}}}.
\end{align}
\itnode{Proof.} We continue to assume $C$ is irreducible without loss of generality. Let $\cu = \set{U_i}$ be an affine open cover of $C$. As usual we let $U_{ij} = U_i \cap U_j$, etc. With respect to $\cu$, the Čech-to-de-Rham double complex of~\eqref{eq:alg_dR_complex} where $E=\emptyset$ is  
\begin{equation}
  \begin{tikzcd}
    \bigoplus_i \co_C(U_i) \arrow[d] \arrow[r] & \bigoplus_i \Omega^1_C(D)(U_i) \arrow[d] \arrow[r] & \dots \\ 
    \bigoplus_{ij} \co_C(U_{ij}) \arrow[d] \arrow[r] & \bigoplus_{ij} \Omega^1_C(D)(U_{ij})  \arrow[d]  \arrow[r] & \dots \\
    \vdots & \vdots 
  \end{tikzcd}
\end{equation}
from which we form the single complex $\ce_0 \overset{\dd_0}{\to} \ce_1 \overset{\dd_1}{\to} \dots$, which is explicitly
\begin{equation}
  \medoplus_i \co_C(U_i) \overset{\dd_0}{\to} \left( \medoplus_{ij}\co_C(U_{ij}) \medoplus_k \Omega^1_C(D)(U_k) \right) \overset{\dd_1}{\to} \dots
\end{equation}
Then, we find
\begin{align}
  \ker d_1 &= \{(a_{ij}; b_k) \in \ce_1 : \text{$(a_{ij})$ is a Čech cocycle, $b_i - b_j = \dd a_{ij}$}\} \\
  \image d_0 &= \{(a_i-a_j; \dd a_k) \in \ce_1 : a_i \text{ is a regular function on $U_i$}\}.
\end{align} %
The quotient $\ker d_1 / \image d_0$ computes $\H^1_\AdR(C\setminus{}D)$. We will take the limit over $\cu$ to account for all possible covers.

If $\omega$ is a $\log(D)$-second kind differential, let $\cu = \set{U_i}$ be an open affine cover for which $\omega|_{U_i} = \eta_i + \dd f_i$, where $\eta_i \in \Omega^1_C(U_i)$, $f_i \in \kappa(C)$. Then $\eta_i-\eta_j = \dd (f_i-f_j)$ implies that $f_i-f_j \in \co_C(U_{ij})$. We may thus consider the tuple $(f_i-f_j; \eta_k) \in \ker d_1$. Another representation $\omega|_{U_i} = \mu_i + \dd g_i$ will give rise to a tuple $(g_i-g_j;\mu_k)$ where the difference $((f_i-g_i)-(f_j-g_j);\eta_k-\mu_k)$ lies in $\image d_0$ since $\eta_i - \mu_i = \dd(f_i-g_i)$. Furthermore, $\omega$ lands in $\ker d_0$ only if it is locally exact, $\omega|_{U_i} = \dd f_i$, but then $\omega$ is exact since $C$ is connected. This proves that the map from right to left of~\eqref{eq:coh_of_punctured_curve} is well defined and injective. 

To prove surjectivity, take the class of $(a_{ij};b_k) \in \ce_1$, fix an index $0$ and observe that $\omega \coleq b_0$ viewed as a global rational differential is of the second kind: $\omega|_{U_j} = b_j+\dd a_{0j}$. 

The statements about the filtrations are straightforward to prove. For~\eqref{eq:F1_puncture} we simply note that a $\log(D)$-regular differential is never exact. As for~\eqref{eq:W1_puncture}, the first identity can be deduced from the long exact sequence for the inclusion $C\setminus{}D \to C$ and the second is deduced from~\eqref{eq:coh_of_punctured_curve} for $D=\emptyset$.
\hfill \qed

\node{}\label{node:prim_E} Let $E\subset C(\Qbar)$ be a finite set. There is a restriction map
\begin{equation}
 \H^0_\AdR(C) \to \H^0_\AdR(E) 
\end{equation}
whose cokernel $\H^0_\AdR(E)_0$ is the \emph{primitive cohomology} of $E$. This notation is also used when $C$ is reducible.

\node{}\label{node:qbarE} Let $\Qbar^E \coleq \hom_{\Set}(E,\Qbar)$ be the $\Qbar$-valued functions on $E$. Thus $\H^0_\AdR(E) = \Qbar^E$. Let $\Qbar^E_0 \coleq \Qbar^E/\H^1_\AdR(C)$ where we mod out by restriction of locally constant functions on $C$. When $C$ is irreducible $\H^1_\AdR(C) = \Qbar$ and the quotient is via the diagonal action of $\Qbar$, i.e., by functions globally constant on $E$. For a tuple $a \in \Qbar^E$ write $[a] \in \Qbar^E_0$ for the class of $a$. The primitive cohomology $\H^0_\AdR(E)_0$ coincides with $\Qbar^E_0$.

\node{}\label{node:df_to_fE} We map $\kappa(C)$ on to the direct sum
\begin{equation}
  \set{\text{$\log(D)$-second kind differentials on $C$}} \oplus \Qbar^E
\end{equation}
by $f \mapsto (\dd f, -f|_E)$. We will refer to the image as \emph{exact forms on $(C,E)$}. 

\node{} Using the identification in Proposition~\eqref{prop:coh_of_punctured_curve} and the one in~\eqref{node:prim_E}, the obvious maps give an exact sequence
\begin{equation}\label{eq:ses_for_differentials}
  0 \to \H^0_\AdR(E)_0 \to  \frac{\set{\text{$\log(D)$-second kind differentials on $C$}} \oplus \Qbar^E}{\set{\text{exact forms on $(C,E)$}}} \to \H^1_\AdR(C\setminus{}D) \to 0.
\end{equation}

\node{Lemma.}\label{lem:coh_of_cde} We have a canonical isomorphism 
\begin{equation}\label{eq:coh_of_cde}
  \H^1_\AdR(C\setminus{}D,E) \simeq \frac{\set{\text{$\log(D)$-second kind differentials on $C$}} \oplus \Qbar^E}{\set{\text{exact forms on $(C,E)$}}}.
\end{equation}

\itnode{Proof.} Both sides of~\eqref{eq:coh_of_cde} fit into the center of the short exact sequence~\eqref{eq:ses_for_differentials}. To deduce the isomorphism, it will be sufficient to prove there is an isomorphism from right to left of~\eqref{eq:coh_of_cde} respecting the extension~\eqref{eq:ses_for_differentials}. 

We continue with the notation of the proof of~\eqref{prop:coh_of_punctured_curve} and make minor modifications to compute the Čech-to-de-Rham double complex of~\eqref{eq:alg_dR_complex}. If $\cu = \set{U_i}$ is an affine open cover of $C$, let $E_i = E \cap U_i$. Then the associated single complex is of the form
\begin{equation}
  \medoplus_i \co_C(U_i) \overset{\dd_0}{\to} \left( \medoplus_{ij}\co_C(U_{ij}) \medoplus_k \Omega^1_C(D)(U_k) \oplus \Qbar^{E_k} \right) \overset{\dd_1}{\to} \dots
\end{equation}
Once again we write 
\begin{align*}
  \ker d_1 &= \{(a_{ij}; b_k,c_k) : \text{$(a_{ij})$ is a Čech cocycle, $b_i - b_j = \dd a_{ij}$, $(c_i-c_j)|_{E_{ij}} = -a_{ij}|_{E_{ij}}$}\}, \\
  \image d_0 &= \{(a_i-a_j; \dd a_k) : a_i \text{ is a regular function on $U_i$}\}.
\end{align*} %
The quotient $\ker \dd_1 / \image \dd_0$ computes $\H^1_\AdR(C\setminus{}D,E)$. We will take the limit over $\cu$ to account for all possible covers. 

Let $(\omega,c) \mod \dd \kappa(C)$ be an element in the right hand side of~\eqref{eq:coh_of_cde}. As in the proof of Proposition~\eqref{prop:coh_of_punctured_curve}, $\omega$ defines an affine open cover $\cu = \set{U_i}$, sections $\eta_i \in \Omega^1_C(D)(U_i)$, and $f_i \in \kappa(C)$, with $\omega|_{U_i} = \eta_i + \dd f_i$. The element $c \in \Qbar^E$ can be restricted to $c_i \coleq c|_{E_i}$. As in Proposition~\eqref{prop:coh_of_punctured_curve}, this gives an element $(f_i-f_j; \eta_k, c_k) \in \ker \dd_1$ which is well defined modulo $\image \dd_0$. Compatibility of this map with both sides of~\eqref{eq:ses_for_differentials} is clear.
\hfill \qed

\node{}\label{node:period_pairing} Let $\omega$ be a $\log(D)$-second type differential on $C$ and $a \in \Qbar^E$. Let $\gamma$ be a smooth $1$-chain on $C(\C)$, avoiding the poles of $\omega$, and with boundary $\partial \gamma$ supported on $E$. Extend $a$ linearly to divisors supported on $E$, $a(\sum n_i p_i) \coleq \sum n_i a(p_i)$. It is easy to check that the integration pairing 
\begin{equation}
  \int_\gamma (\omega,a) \coleq \int_{\gamma} \omega + a(\partial \gamma),
\end{equation}
is well-defined on the singular homology class $[\gamma] \in \H_1^\B(C\setminus{}D,E)$ and on the class of $(\omega,a)$ modulo $(C,E)$-exact forms. Via the isomorphism in Lemma~\eqref{lem:coh_of_cde}, we can describe the ``period pairing''
\begin{align}%
  \H_1^\B(C\setminus{}D,E) \otimes \H^1_\AdR(C\setminus{}D,E) &\to \C \\
  [\gamma] \otimes [\omega,a] \mapsto \int_{\gamma} (\omega,a).
\end{align}
This pairing agrees with the usual comparison isomorphism between $\H^1_\AdR \otimes \C$ and $\H^1_\B \otimes \C$, as can be checked by expanding the proof of the comparison isomorphism~\cite[Theorem~2.150]{BurgosFresan}. %

\subsection{An effective basis for cohomology}

\node{Definition.}\label{node:differential_basis} A tuple $(\omega_1,a_1), \dots, (\omega_m,a_m)$ is called a \emph{differential basis} for the cohomology $\H^1_\AdR(C\setminus{}D,E)$ if $\omega_i$'s are $\log(D)$-second type differentials, $a_i \in \Qbar^E$ are functions on $E$, and the $(\omega_i,a_i)$'s map to a basis via the identification~\eqref{lem:coh_of_cde}. We will refer to a differential basis as a representation of $\H^1_\AdR(C\setminus{}D,E)$. 

\node{} In this subsection and the next we will show how to compute a differential basis together with the Hodge and weight filtrations on cohomology. See~\eqref{prop:effective_cohom_basis} for a summary.

\node{}\label{node:weierstrass} A point $p \in C(\Qbar)$ is a Weierstrass point if $h^0(g p) > 1$ where $g$ is the genus of $C$. There are at most $g(g^2-1)$ Weierstrass points on a curve~\cite[p.274]{Griffiths1994}, hence a non-Weierstrass point can be found in finite time by trial-and-error. 

\node{Proposition.}\label{prop:eff_C_coh} Let $C$ be an irreducible curve of genus $g$ and $p \in C(\Qbar)$ a non-Weierstrass point. Then, there is a natural isomorphism
\begin{equation}\label{eq:simple_cohomology}
  \H^1_\AdR(C) \simeq \H^0(C,\Omega^1_C((g+1) p)).
\end{equation}
If $p \in C(\Qbar)$ is arbitrary, then there exists an integer $d_0$, $g < d_0 \le 2g $, such that for every $d \ge d_0$ we have $h^0((d-1)p) = d-g$. Then, for all $d \ge d_0$, we have a natural isomorphism
\begin{equation}\label{eq:simple_cohomology_v2}
  \H^1_\AdR(C) \simeq \frac{\H^0(C,\Omega^1_C(d p))}{\dd \H^0(C,\co_C((d-1)p))}.
\end{equation}
\itnode{Proof.} All cohomology classes can be represented by a differential of the second kind with poles only on $p$. To see this repeat the proof of Proposition~\eqref{prop:coh_of_punctured_curve} with the cover $\cu=\{U_0,U_1\}$ where $U_0 = C\setminus{}p$ and $U_1$ an affine neighbourhood of $p$. 

All global sections of $\Omega^1_C(dp)$ for $d \ge 0$ are of the second kind since, by the residue theorem, a rational form with poles supported at a single point can not have any residue. 

Thus, we have a natural isomorphism
  \begin{equation}\label{eq:point_pole}
    \H^1_\AdR(C) \simeq \frac{\H^0(C,\Omega^1_C(d p))}{\dd \H^0(C,\co_C((d-1)p))},\quad d \gg 0.
  \end{equation}
Since $\dim \H^1_\AdR(C) = 2g$, we can pick any $d$ for which the right hand side is of dimension $2g$.

By Riemann--Roch we have
\begin{equation}
  h^0(\Omega^1_C((g+1)p)) = 2g.
\end{equation}
If $p$ is non-Weierstrass, then $h^0(gp)=1$ and the space $\H^0(C,\co_C(gp))$ consists only of constants, which are killed by derivation $\dd$. Thus we may take $d=g+1$ and this proves~\eqref{eq:simple_cohomology}.

We now consider the case when $p$ is arbitrary. The space in the denominator of~\eqref{eq:point_pole} is of dimension $h^0((d-1)p)-1$ since derivation kills constants. The space in the numerator is of dimension $h^0(\Omega^1_C(dp)) = g+d-1$ for $d > 0$ by Riemann--Roch. We want to solve for $d>0$ such that
  \begin{equation}\label{eq:desired_equality}
    2g = h^0(\Omega^1_C(dp)) - h^0((d-1)p) + 1,
  \end{equation}
or equivalently $h^0((d-1)p) = d-g$. By Riemann--Roch this holds precisely when 
\begin{equation}
h^0(\Omega^1_C((-d+1)p) = 0.
\end{equation}
This is true whenever $d \ge 2g$ but it can never hold unless $g < d$ since $h^0(\Omega^1_C) = g$. Hence the minimal $d$ will satisfy $g<d\le 2g$.
\hfill\qed

\node{Remark.} A restatement of Proposition~\eqref{prop:eff_C_coh}, in light of Proposition~\eqref{prop:coh_of_punctured_curve}, is that a differential of the second kind on $C$ is equivalent, modulo exact forms, to a \emph{unique} differential with poles only at $p$, provided $p$ is non-Weierstrass. This is analogous to the situation with forms of the first kind where regular differentials uniquely represent their cohomology classes.

\node{} For a non-Weierstrass point $p \in C(\Qbar)\setminus{}D$, define
\begin{equation}\label{def:WDp}
  W_{D,p} \coleq \ker \left( \res_p \colon \H^0(C,\Omega^1_C(D+(g+1)p)) \to \Qbar \right).
\end{equation}
By construction, any element of $W_{D,p}$ is regular on $C\setminus{}\left(D \cup \set{p}\right)$, has at most a simple poles on $D$ and any pole on $p$ is non-residual. In particular, $W_{D,p}$ consists of differentials of the $\log(D)$-second kind.

\node{Proposition.}\label{prop:WDp_coh} For $p \in C(\Qbar)\setminus{}D$ non-Weierstrass, the identification~\eqref{prop:coh_of_punctured_curve} gives a natural isomorphism, 
\begin{equation}\label{eq:W_map}
  W_{D,p} \to \H^1_\AdR(C\setminus{}D) : \eta \mapsto \eta \text{ modulo exact forms}.
\end{equation}
\itnode{Proof.}
Let $s \coleq \# D - 1$. Using Riemann--Roch and the residue theorem we find 
\begin{equation}
  \dim W_{D,p} = \begin{cases}
    2g & D = \emptyset \\ 2g+s & s \ge 0
  \end{cases}.
\end{equation}
This agrees with the dimension of $\H^1_\AdR(C\setminus{}D)$ since it fits into the exact sequence
\begin{equation}
 0 \to \H^1_\AdR(C) \to \H^1_\AdR(C\setminus{}D) \overset{\Res}{\to} {}_0\Qbar^D \to 0,
\end{equation}
where ${}_0\Qbar^D \subset \Qbar^D$ is the space of functions on $D$ whose value on $D$ sum to $0$ (the dual of~$\Qbar^D_0$).

Therefore, to show that~\eqref{eq:W_map} is an isomorphism, it will be sufficient to show that it is surjective. 

We proved in Proposition~\eqref{prop:eff_C_coh} that the subspace
\begin{equation}
  \H^0(C,\Omega^1_C((g+1)p)) \subset W_{D,p}
\end{equation}
surjects onto $\H^1(C)$. And when $s>0$, $h^1(\Omega^1_C(D))=0$ implies $\H^0(C,\Omega^1_C(D)) \subset W_{D,p}$ surjects onto ${}_0\Qbar^D$.
\hfill\qed

\node{Corollary.}\label{cor:represent_cohom} For $p \in C(\Qbar)\setminus{}D$ non-Weierstrass, we have a natural isomorphism  
\begin{equation}\label{eq:split_cohom}
  \H^1_\AdR(C\setminus{}D,E) \simeq W_{D,p} \oplus \Qbar^E_0 ,
\end{equation}
which gives a splitting of the short exact sequence
\begin{equation}
  0 \to \H_\AdR^0(E)_0 \to \H^1_\AdR(C\setminus{}D,E) \to \H^1_\AdR(C\setminus{}D) \to 0
\end{equation}
that depends only on $p$ (but it \emph{does} depend on $p$). 
\itnode{Proof.} Since~\eqref{eq:W_map} is an isomorphism, there are no exact forms in $W_{D,p}$. Thus the map $W_{D,p} \to \H^1_\AdR(C\setminus{}D,E) : \eta \mapsto [\eta,0]$, using the identification~\eqref{eq:coh_of_cde}, is injective. This gives the splitting in question, in light of Proposition~\eqref{prop:WDp_coh}. The identification $\H_\AdR^0(E)_0 \simeq \Qbar^E_0$ is definitional~\eqref{node:prim_E}. 
\hfill \qed 

\node{}\label{node:explicit_cohom_filtrations} Using the splitting in~\eqref{eq:split_cohom}, we can make the weight and Hodge filtrations on $\H^1_\AdR(C\setminus{}D,E)$ explicit. The non-trivial pieces are thus:
\begin{equation}
  W_0 = 0 \oplus \Qbar^E_0, \quad W_1 =  W_{\emptyset,p} \oplus \Qbar^E_0, \quad  F^1 = \H^0(C,\Omega^1_C(D)) \oplus 0.
\end{equation}

\subsection{Computing a basis for algebraic de Rham cohomology}

\node{} In this subsection, we will give references to algorithms that allow us to work effectively with linear systems on curves. We wish to emphasize that a basis for $W_{D,p}$~\eqref{def:WDp} can be computed effectively as a result.

\node{} For $f\in \kappa(C)$ let $\divv(f)$ be the associated divisor of zeros and poles of $f$. For any divisor $D'$ on $C$, we write
\begin{equation}
L(D') \coleq \{f \in \kappa(C) \mid \divv(f) + D' \ge 0 \} \subset \kappa(C),
\end{equation}
and note that $L(D') = \H^0(C,\co_C(D'))$.

\node{}\label{node:effective_linear_system} We can algorithmically~\cite{Hess2002} evaluate the map
\begin{equation}
  D' \mapsto L(D'),
\end{equation}
which takes a divisor supported on $\Qbar$ points of $C$ and outputs (a basis for) the linear system $L(D')$.

\node{}\label{node:basis_of_holo_differentials} Given $C$ (or, rather, a birational plane model $P(x,y)=0$) we can determine a basis for the global sections of the canonical bundle as lying in the space of rational differentials:
\begin{equation}
  \H^0(C,\Omega^1_C) \subset \Omega^1_{\kappa(C)/\Qbar} = \kappa(C) \dd x,
\end{equation}
this is explained in~\cite[Appendix A, \S 2]{GACI}.

\node{}\label{node:basis_of_differentials} For any divisor $D_1$ we can compute a basis for the global sections of $\Omega^1_C(D_1)$ as follows. Find~\eqref{node:basis_of_holo_differentials} a non-zero section $\omega_0 \in \H^0(C,\Omega^1_C)$ and determine its divisor of zeros $D_0$. The space $\H^0(C,\Omega^1_C(D_1))$ is identified with $L(D_0+D_1) \subset \kappa(C)$ through multiplication with $\omega_0$, i.e.,
\begin{equation}
  \H^0(C,\Omega^1_C(D_1)) = L(D_0+D_1)\cdot \omega_0 \subset \Omega^1_{\kappa(C)/\Qbar}. 
\end{equation}
Now use~\eqref{node:effective_linear_system}.

\node{Proposition.}\label{prop:effective_cohom_basis} Given disjoint finite sets $D,E \subset C(\Qbar)$, we can effectively compute a differential basis for~$\H^1_\AdR(C\setminus{}D,E)$ together with the Hodge and weight filtrations. 
\itnode{Proof.} Find a non-Weierstrass point $p \in C(\Qbar)\setminus{}D$~\eqref{node:weierstrass}.
Compute a basis for $\H^0(C,\Omega^1_C(D + (g+1)p))$~\eqref{node:basis_of_differentials}. Compute the residue at $p$ for each basis element~\eqref{node:residue} and compute the kernel $W_{D,p}$. Now use the identity~\eqref{eq:split_cohom} to get a differential basis. The Hodge and weight filtrations are~\eqref{node:explicit_cohom_filtrations}. \hfill\qed

\subsection{Reduction algorithms}\label{sec:reduction_algos}

\node{} A crucial point for our algorithms is to represent a rational differential in terms of a basis for an algebraic de Rham cohomology. We will call this a reduction algorithm.

\node{} The curve $C$ is irreducible smooth projective and defined by a plane model $P(x,y) =0$. Fix $D, E \subset C(\Qbar)$, disjoint finite sets. Let $p \in C(\Qbar)\setminus{}D$ be non-Weierstrass.

\node{} Any rational differential $\mu \in \Omega^1_{\kappa(C)/\Qbar}$ with $\Res(\mu)$ supported on $D$ defines a cohomology class $[\mu,0] \in \H^1_\AdR(C\setminus{}D,E)$ via~\eqref{lem:coh_of_cde}. By Corollary~\eqref{cor:represent_cohom}, there exists a unique element $\eta \in W_{D,p}$ and a rational function $f \in \kappa(C)$ such that $\eta \equiv \mu + \dd f$. Thus $(\mu,0) \equiv (\eta,-f|_E)$ in $\H^1_\AdR(C\setminus{}D,E)$. A \emph{reduction algorithm} takes $\mu$ as input and returns the pair $\eta \in W_{D,p}$ and $[-f|_E] \in \Qbar^E_0$.

\node{A symbolic reduction algorithm.} \label{node:symbolic_reduction}
Express $\mu = h \dd x$ and write $\dd f = f' \dd x$ (for unknown $f$ and $f'$, where $f' = f_x - f_yP_x/P_y$). We recall that $W_{D,p} \subset \H^0(C,\Omega^1_C(D + (g+1)p))$.

We have $\eta = \mu + \dd f = (h+f') \dd x$. Here $f'$ is uniquely determined by the property that $\divv(\eta) + D+(g+1)p \ge 0$, that is
\begin{equation}\label{eq:hf}
  \divv(h+f') + D' \ge 0, \quad \text{where } D' \coleq \divv(\dd x) + D + (g+1)p.
\end{equation}
But $\divv(f') \ge \min(\divv(h),\divv(h+f'))$, where minimum is applied pointwise. Let 
\begin{equation}
  D'' \coleq \min(\divv(h),-D').
\end{equation}
Therefore, $f' \in L(D'')$ is the unique element (up to scaling) such that $\divv(h+f') + D' \ge 0$. Note that the residue at $p$ is automatically $0$. 

We can compute the Laurent tails of $h$ and of a basis for $L(D'')$ along the negative part of $D''$. Since $f'$ is the element whose Laurent tail kills the offending poles of $h$, we can now determine the correct multiple of $f'$.  

Given the $\log(D)$-second kind differential $\mu=h \dd x$, the output of this computation will be $\dd f= f' \dd x$ such that $\eta \coleq \mu + \dd f \in W_{D,p}$. 

From the poles of $f'$, we can determine a finite dimensional linear system containing $f$, and thus determine $f$ up to additive constants. Restricting $f$ gives the well-defined class $[f|_E]\in\Qbar^E_0$.

\node{A numerical reduction algorithm.} %
This time, we will reconstruct $\eta$ from numerical approximations. Despite the heuristic intermediate steps, the output certifies itself.

Find a basis $\{\omega_i\}_{i=1}^k $ for $W_{D,p} \subset \H^0(C,\Omega^1_C(D + (g+1)p))$ and a basis of rectilinear chains for $\{\gamma_i\}_{i=1}^k$ for $\H_1^\B(C\setminus{}D,E)$ avoiding $p$~\eqref{prop:rectilinear_basis_computable}.
Compute (a numerical approximation of) the period matrix $\left( \int_{\gamma_j} \omega_i \right)_{i,j = 1,\dots,k}$ and (a numerical approximation of) the row vector $\left( \int_{\gamma_j} \mu \right)_{j=1,\dots,k}$. Let $a=(a_i) \in \Qbar^k$ be the row vector satisfying
\begin{equation}
  a = \left( \int_{\gamma_j} \mu \right)_{j=1,\dots,k} \left( \int_{\gamma_j} \omega_i \right)_{i,j = 1,\dots,k}^{-1}.
\end{equation}
Then $\eta = \sum_i a_i \omega_i $.

In practice, we will reconstruct, i.e., guess, $a$ from a numerical approximation. To prove our guess, we make the following observations. If $\eta = \mu + \dd f$ then
\begin{equation}
  \divv(f') = \divv(\eta-\mu) - \divv(\dd x) =: D'''  
\end{equation}
Therefore, $D'''$ is a principal divisor and $f' \in L(D''') \simeq \Qbar$. We can recover the correct multiple of $f'$ in $L(D''')$ by making sure the identity $f' \dd x= \eta - \mu$ holds at a general point $q \in C(\Qbar)$. 

For $\eta$ that is merely formulated as a guess, we can construct the linear system $L(D''')$ and try to find $f'$ such that $\eta = \omega + f' \dd x$. Failure to find $f'$ means the guess is incorrect. We must go back and try to reconstruct $a$ using higher precision and try to find $f'$ again. Since $\eta$ with desired properties exists, this process must terminate.

\subsection{Cohomological pullback and the transfer map}\label{sec:coh_pullback}

\node{}\label{node:rep_curve_map} Let $C_1$ and $C_2$ be two irreducible smooth proper curves over $\Qbar$. Let $f \colon C_1 \to C_2$ be a non-constant map between the curves. The curves are represented by their function fields $\kappa(C_i)$~\eqref{node:rep_curve} and so we define a \emph{representation} of $f$ to be the inclusion of function fields $f^* \colon \kappa(C_2) \to \kappa(C_1)$ in terms of the generators. The effective pullback map on functions gives an effective pullback map on differential forms $f^* \colon \Omega^1_{\kappa(C_1)/\Qbar} \to \Omega^1_{\kappa(C_2)/\Qbar}$ via $f^*(r(x,y)\dd x) = f^*(r(x,y))\dd f^*(x)$.

\node{}\label{node:coh_pullback} For each $i=1,2$ let $D_i, E_i \subset C_i(\Qbar)$ be disjoint sets. Suppose $f^{-1}(D_2) \subset D_1$ and $f(E_1) \subset E_2$. Then, $f$ can be viewed as a map $f \colon (C_1\setminus{}D_1,E_1) \to (C_2\setminus{}D_2,E_2)$. Thus, we have the pullback map on cohomology $f^*_\AdR \colon \H^1_\AdR(C_2\setminus{}D_2,E_2) \to \H^1_\AdR(C_1\setminus{}D_1,E_1)$. 

\node{Proposition.}\label{prop:coh_pullback_matrix} Given $f$ as in~\eqref{node:coh_pullback} and bases for the cohomology groups $\H^i_\AdR(C_i\setminus{}D_i,E_i)$ as in~\eqref{cor:represent_cohom}, we can determine the matrix representing the cohomological pullback $f^*_\AdR$ with respect to the given bases. 
\itnode{Proof.}
Using the representation of $f$~\eqref{node:rep_curve_map}, compute the pullback $f^* \colon W_{D_2,p_2} \to \Omega^1_{\kappa(C_1)/\Qbar}$. Now apply a reduction algorithm~\eqref{node:symbolic_reduction} to express a basis for $W_{D_2,p_2}$ in terms of the basis for $W_{D_1,p_1}\oplus\Qbar_0^{E_1}$. The map $\Qbar_0^{E_2} \to \Qbar_0^{E_1}$ descends from the pullback map of $f|_{E_1} \colon E_1 \to E_2$.
\hfill\qed

\node{} Given non-constant $f \colon C_1 \to C_2$, the field extension $f^* \colon \kappa(C_2) \to \kappa(C_1)$ is finite and hence we have the \emph{trace map} $\tr_f \colon \kappa(C_1) \to \kappa(C_2)$. The trace is effective on individual elements: given a representation of $f$ and an element $a \in \kappa(C_1)$ we can compute the trace $\tr_f(a)$. 

\node{} The trace map extends to Kähler differentials $\tr_f \colon \Omega^1_{\kappa(C_1)/\Qbar} \to \Omega^1_{\kappa(C_2)/\Qbar}$ where we have $\tr_f(r\dd x) = \tr_f(r) \dd\tr_f(x)$. %

\node{}\label{node:coh_transfer} For each $i=1,2$ let $D_i, E_i \subset C_i(\Qbar)$ be disjoint sets. This time, suppose $f^{-1}(E_2) \subset E_1$ and $f(D_1) \subset D_2$. The trace map on differentials induces the \emph{transfer map}
\begin{equation}
  f_!^\AdR \colon \H^1_\AdR(C_1\setminus{}D_1,E_1) \to \H^1_\AdR(C_2\setminus{}D_2,E_2).
\end{equation}

\node{Remark.} The transfer map $f_!^\AdR$ is also known as the ``integration over fibers'' or the ``Gysin'' map. 

\node{Remark.}\label{rem:poincare_lefschetz} Observe that our hypotheses give a map $f \colon (C_1\setminus{}E_1,D_1) \to (C_2\setminus{}E_2,D_2)$. The Poincaré--Lefschetz duality takes the form $\H^1_\AdR(C_i\setminus{}E_i,D_i) \simeq \H^1_\AdR(C_i\setminus{}D_i,E_i)^\vee$. Taking the dual of cohomological pullback $f^*_\AdR$ and composing with the Poincaré--Lefschetz dualities gives the transfer map. 

\node{Proposition.}\label{prop:coh_transfer_matrix} As in~\eqref{prop:coh_pullback_matrix}, we can compute the matrix representing the transfer map $f_!^\AdR$ in given bases.
\itnode{Proof.}
The precise statement and the proof is similar to that of~\eqref{prop:coh_pullback_matrix}. On $W_{C_1,p_1} \to \Omega^1_{\kappa(C_2)/\Qbar}$ we use the trace map and then apply the reduction algorithm. Write $f^* \colon \Div(C_2) \to \Div(C_1)$ for the pullback on divisors and extend the action of $\Qbar^{E_1}$ additively to divisors supported on $E_1$. Then the restriction $\Qbar^{E_1}_0 \to \Qbar^{E_2}_0$ of the trace map descends from 
\begin{equation}
  a \mapsto \left(e \mapsto a(f^*(e)\right), 
\end{equation}
which is clearly effective.
\hfill\qed

\newpage
\section{Betti homology and the mixed Hodge structure of a punctured marked curve}\label{sec:betti}

\node{} Let $(C \setminus D, E)$ be a punctured marked curve over $\Qbar$. In this section, we complete the effective construction of the mixed Hodge structure 
\begin{equation}
\H_1(C \setminus D, E) = (\H_1^\B(C\setminus{}D,E),\H^1_\AdR(C\setminus{}D,E)^\vee,\comp)
\end{equation}
by giving an algorithm to compute a basis for the Betti (singular) homology~\eqref{prop:rectilinear_basis_computable} so that the comparison isomorphism can be approximated to arbitrary precision~\eqref{node:approx_integration}. We also give an algorithm to compute the coordinates of the homology class of a smooth $1$-chain in the given basis of homology~\eqref{node:hom_reduction}. Furthermore, for a map $p$ between punctured marked curves, we will describe how to compute the pushforward $p_*$ and transfer $p^!$ maps in homology, see~\eqref{prop:effective_hom_pushforward} and~\eqref{prop:effective_hom_transfer}. Combined with the corresponding results for $\H^1_\AdR$, we get effective pushforward and transfer maps on the mixed Hodge structure $\H_1$ of a punctured marked curve. In \S\ref{sec:corr}, we will use these to determine the action of correspondances on $\H_1$.

\node{} The idea to compute the Betti homology is to embed a graph $\Gamma$ in $C\setminus{}D$ with ``rectilinear edges''~\eqref{node:rectilinear_graph} and containing $E$ such that cycles on $(\Gamma,E)$ generate the integral Betti homology $\H_1^\B(C\setminus{}D,E;\Z)$. We construct the graph so that the bounding cycles on $C\setminus{}D$ are evident~\eqref{node:construct_rep_C}. The graph can be constructed to  avoid a given finite set $S$ of points in $C(\Qbar)\setminus{}E$. We pick $S$ to support the poles of a differential basis~\eqref{node:differential_basis} for $\H^1_\AdR(C\setminus{}D,E)$. Thus, we may pair the cycles on $(\Gamma,E)$ against the differential basis via integration to determine the period pairing~\eqref{node:period_pairing}.

\node{} The problem of computing the singular homology of an algebraic curve from its defining equation, and pairing it with algebraic de Rham cohomology, has a long history, even algorithmically. Foundational work by Tretkoff--Tretkoff~\cite{Tretkoff1984} introduced combinatorial methods via monodromy representations; later, Deconinck--van Hoeij~\cite{Deconinck2001} developed a general symbolic-numeric algorithm for computing homology bases and period matrices. More recent work use embedded graphs and certified analytic continuation to construct homology bases, particularly for superelliptic or general plane curves, enabling high-precision period computations~\cite{Bruin2019,Molin2018}. 

\node{} Our construction of $\H_1^\B(C \setminus D, E)$ generalizes and streamlines several of the ideas outlined above. The strategy is inspired by the work of Bruin, Sijsling, and Zotine~\cite{Bruin2019}, who embed a graph on the curve to facilitate homology computations. We depart from their method by explicitly constructing bounding cycles supported on the graph, rather than computing the radical of an intersection pairing. Our formulation is also slightly more abstract, and we refer to their work for implementation-level details.

\subsection{Representing the mixed Hodge structure}\label{sec:mhs}

\node{} The $k$-th Betti homology group $\H_k^\B(C\setminus D, E;\Z)$ denotes the singular homology of $C(\C) \setminus D$ relative to $E$, with integral coefficients. We will describe how to compute a representation of this homology group later in the section. We write $\H^1_\B$ for the dual $(\H_1^\B)^\vee$. The tuple $\H_1 \coleq ((\H^1_\B,W),(\H^1_\AdR,W,F),\varphi_{\AdR,\B})$ is a $(\Q,\Qbar)$-mixed Hodge structure. Here, we explain how to represent the mixed Hodge structure.

\node{}\label{node:representation_of_MHS} We \emph{represent} the mixed Hodge structure
\begin{equation}
\H_1(C \setminus D, E) \coleq \left( \H_1^\B(C\setminus D, E), \, \H^1_\AdR(C\setminus D, E)^\vee, \, \comp \right)
\end{equation}
by representing $\H^1_\AdR(C\setminus{}D,E)$~\eqref{node:differential_basis} and $\H_1^\B(C\setminus{}D,E)$~\eqref{prop:rep_hom}. 

\node{}\label{node:compute_rep_of_MHS} We can compute a representation of $\H_1(C \setminus D, E)$. Compute~\eqref{prop:effective_cohom_basis} a differential basis $(\omega_1,a_1),\dots,(\omega_m,a_m)$ for $\H^1_\AdR(C\setminus{}D,E)$. Let $S \subset C(\Qbar) \setminus (D\cup E)$ be a finite set supporting the poles of $\omega_i$'s, e.g., $S$ consists of one non-Weierstrass point per component. Compute~\eqref{prop:rectilinear_basis_computable} a rectilinear basis $\gamma_1,\dots,\gamma_m$ for $\H_1^\B(C\setminus{}D,E;\Z)$ avoiding $S$.

\node{} A representation of $\H_1(C\setminus{}D,E)$ allow us to determine the comparison isomorphism $\comp$ up to desired precision. If the given rectilinear basis $\gamma_i$'s do not avoid the poles $S$ of the differential basis, we can compute a rectilinear basis that does, and then compute~\eqref{prop:effective_hom_pushforward} the change of coordinates between the two rectilinear bases. Therefore, we will now assume our rectilinear avoids the poles of the differential basis.

\node{}\label{node:approx_integration} The period pairing~\eqref{node:period_pairing} is represented by the matrix
\begin{equation}
  \begin{pmatrix}
    \int_{\gamma_j} (\omega_i,a_i)
  \end{pmatrix}_{i,j = 1,\dots,N}
\end{equation}
The integrals $\int_{\gamma_j} \omega_i$ are readily approximated to desired degree of precision from the representations of $\gamma_j$ and $\omega_i$, see for instance~\cite{Bruin2019}. In fact, we can compute the integrals with rigorous error bounds, i.e., as complex balls containing the true value~\cite{Mezzarobba2016,Lairez2019a}. 

\node{} The period matrix represents the comparison isomorphism
\begin{equation}
 \comp \colon \H_1^\B(C\setminus{}D,E;\Z)\otimes_\Z \C \isoto \H^1_\AdR(C\setminus{}D,E)^\vee\otimes_\Qbar \C.
\end{equation}

\node{} In theory, we typically work with the rational vector space $\H_1^\B(C\setminus D,E) \coleq \H_1^\B(C\setminus D,E;\Z) \otimes_\Z \Q$, but the integral structure is crucial for algorithms involving lattices.

\subsection{Rectilinear chains}\label{sec:rectilinear}

\node{Definition.}\label{node:line_segment_p1} Consider the interval $[0,1] \subset \R$ embedded into $\p(\C)$ via $t \mapsto [t:1]$. Composing by a linear map $A \colon \p_\Qbar \isoto \p_\Qbar$ defines a smooth path $\ell_A \colon [0,1] \to \p(\C)$, which we will call a \emph{line segment} on $\p$.

\node{}\label{node:alg_end_points} The end points $\ell_A(0) = A([0:1]),\, \ell_A(1) = A([1:1]) \in \p(\Qbar)$ of a line segment are necessarily algebraic.

\node{} We will \emph{represent} a line segment $\ell_A$ on $\p$ by a $2\times 2$ matrix of algebraic numbers inducing the linear map $A$. 

\node{} Let $C/\Qbar$ be a smooth proper curve and $x \colon C \to \p$ a finite map. Let $d = \deg x$ and $B \subset \p(\Qbar)$ be the branch locus of $x$.

\node{} Let $\ell_A$ be a line segment on $\p(\C)$ with the interior $\ell_A(0,1)$ disjoint from $B$. Then $x^{-1}(\ell_A(0,1))$ has $d$ disjoint components. The choice of a point $r \in x^{-1}(\ell_A(1/2))$ identifies the unique component of $x^{-1}(\ell_A(0,1))$ passing through $r$.  Let $\ell_{A,r} \colon [0,1] \to C(\C)$ be the analytic continuation of $x^{-1}\circ \ell_A$ passing through $r$. The map $\ell_{A,r}$ is smooth on $(0,1)$ and continuous at the boundary. 

\node{Definition.}\label{node:line_segment_C} Any map of the form $\ell_{A,r}$ is a \emph{line segment} on $C$ with respect to $x$. We will not mention $x$ when it is clear from context. The tuple $(A,r)$ \emph{represents} $\ell_{A,r}$.

\node{} Observe that $\ell_A(1/2) \in \p(\Qbar)$ and thus $r \in C(\Qbar)$. Therefore, the tuple $(A,r)$ is finitely presented. %

\node{} Let $\omega$ be a rational differential on $C$ and suppose that the line segment $\ell_{A,r}$ avoids the poles of $\omega$. Then $\ell_{A,r}^*\omega$ is continuous on $[0,1]$ because the $x$-coordinate of $\ell_{A,r}$ is smooth at the boundary and $\omega = f(x,y)\dd x$. In particular, the integral of $\omega$ over $\ell_{A,r}$ is well-defined.

\node{Definition.}\label{node:rectilinear_def} A \emph{rectilinear chain} $\gamma = \sum_{i} a_i \ell_{A_i,r_i}$ on $C$ with respect to $x$ is a formal $\Z$-linear combination of line segments on $C$ with respect to $x$. We will not mention $x$ when it is clear from context. The \emph{representation of $\gamma$} is the formal $\Z$-linear combination $\sum_i a_i (A_i,r_i)$ of the representations of the line segments.

\node{}\label{node:corners} All ``corners,'' and in particular the end points, of a rectilinear chain are algebraic~\eqref{node:alg_end_points}.

\subsection{Rectilinear embedded graphs in curves}

\node{Definition.}\label{def:rectilinear_basis} A \emph{rectilinear basis} for $\H_1^\B(C\setminus{}D,E)$ is a tuple $\gamma_1,\dots,\gamma_m$ of rectilinear chains~\eqref{node:rectilinear_def} whose homology classes form a basis. We will say that a rectilinear basis \emph{avoids a set $S \subset C(\C)$} if the support of the chains are disjoint from $S$.

\node{} In this subsection, we will reformulate the problem of finding a rectilinear basis~\eqref{node:hom_basis}. In the next subsection, we will show how to compute a rectilinear basis~\eqref{prop:rectilinear_basis_computable}.

\node{} A finite oriented graph $\Gamma = (s,t \colon \Gamma_E \to \Gamma_V )$ is an ordered pair of maps between two finite sets, from the ``edge set'' $\Gamma_E$ to the ``vertex set'' $\Gamma_V$, with the maps $s,t$ assigning the two ``ends'' of an edge. 

\node{} Let $\Z[\Gamma_E],\Z[\Gamma_V]$ be the formal $\Z$-linear combinations of the edges and vertices of $\Gamma$.
The simplicial homology of the associated topological space $|\Gamma|$ is computed via the homology of the complex
\begin{equation}\label{eq:chain_complex_of_gamma}
 0 \to \Z[\Gamma_E]  \overset{t-s}{\too} \Z[\Gamma_V] \to 0.
\end{equation}
Thus $\H_1^\B(|\Gamma|;\Z) \simeq \ker(t-s)$ and $\H_0^\B(|\Gamma|;\Z) \simeq \coker(t-s)$. Evidently, given $\Gamma$, the homology groups are computable.

\node{} Let $A \subset \Gamma_V$ be a subset. The relative homology of the pair $(|\Gamma|,A)$ is computed by the complex attained by replacing $\Z[\Gamma_V]$ with $\Z[\Gamma_V]/\Z[A]$ in~\eqref{eq:chain_complex_of_gamma}. %

\node{} An \emph{embedded graph} of the punctured marked curve $(C\setminus{}D,E)$ is a finite oriented graph $\Gamma$, a subset $A \subset \Gamma_V$, and a topological embedding $\iota \colon (|\Gamma|,A) \to (C\setminus{}D,E)$. 

\node{} For each edge $e \in \Gamma_E$ of $\Gamma$, denote the corresponding $1$-simplex by $\phi_e \colon [0,1] \to |\Gamma|$.

\node{}\label{node:rectilinear_graph} Let $x \colon C \to \p$ be a finite map. An embedded graph $(\Gamma,A,\iota)$ is rectilinear with respect to $x$ if each parametrized edge $\iota\circ \phi_e \colon [0,1] \to C(\C)$ is a line segment with respect to $x$~\eqref{node:line_segment_C}. 

\node{Definition.}\label{node:rep_of_hom} A \emph{representation} of $\H_1^\B(C\setminus{}D,E;\Z)$ is a rectilinear embedded graph $\iota\colon(\Gamma,A) \toi (C\setminus{}D,E)$ and a subgroup $K \subset \H_1^\B(|\Gamma|,A;\Z)$ such that the pushforward map
\begin{equation}
  \iota_* \colon \H_1^\B(|\Gamma|,A;\Z) \to \H_1^\B(C\setminus{}D,E;\Z)
\end{equation}
induces an isomorphism
\begin{equation}
 \H_1^\B(|\Gamma|,A;\Z)/K \simeq \H_1^\B(C\setminus{}D,E;\Z).
\end{equation}

\node{}\label{node:hom_basis} Given a representation of homology, we can immediately determine a rectilinear basis for homology. Indeed, using the notation from~\eqref{node:rep_of_hom}, find a basis for $\H_1^\B(|\Gamma|,A;\Z)/K$, split the surjection $\H_1^\B(|\Gamma|,A;\Z) \to \H_1^\B(|\Gamma|,A;\Z)/K$, and embed the generators to $C(\C)$ via $\iota$.

\subsection{Constructing a representation of the homology of a curve}

\node{} Let $C$ be a smooth proper curve and $D,E \subset C(\Qbar)$ disjoint sets. Fix a finite map $x \colon C \to \p$ and a finite set $S \subset C(\Qbar) \setminus E$. Let $B \subset \p(\Qbar)$ be the branch locus of $x$. This subsection is dedicated to proving the following statement.

\node{Proposition.}\label{prop:rectilinear_basis_computable} Given $(C\setminus{}D,E)$ we can effectively compute a rectilinear basis~\eqref{def:rectilinear_basis} for $\H_1^\B(C\setminus{}D,E;\Z)$ avoiding $S$. The weight filtration~\eqref{node:betti_weight_filtration} on homology is effective.

\node{} In light of~\eqref{node:hom_basis}, it suffices to construct a representation of homology. We will first perform this construction when $D=E=\emptyset$.

\node{} If $B = \emptyset$ then $C$ is a disjoint union of $\p$'s and $x$ is an isomorphism on each component. In this case, $\H_1^\B(C) = 0$ and we may take the empty graph to represent the homology. We will exclude this case below by assuming $B \neq \emptyset$.

\node{} Let $\Gamma \subset \p\setminus{}B$ be a rectilinear embedded graph and a deformation retract of $\p\setminus{}B$. This is effective, e.g., take the boundary of a Voronoi diagram around $B \neq \emptyset$. For each $b\in B$, there is a unique component $R_b$ of $\p(\C)\setminus{}\Gamma$ containing $b$. Each $R_b$ is homeomorphic to a disk with boundary $\partial R_b$ on $\Gamma$. We thus realize $\p(\C)$ as a $2$-dimensional CW-complex.

\node{}\label{node:construct_rep_C} The preimage $\wGamma \coleq x^{-1}\Gamma$ in $C(\C)$ is a rectilinear embedded graph of $C$. Define $K \subset \Z[\wGamma_E]$ to be the subgroup generated by the boundary components of $x^{-1} R_b$ for each $b\in B$.

\node{Remark.} Lifting the edges constituting $\partial R_b$ and describing how they connect is part of the computation required in constructing $\wGamma$, thus obtaining $K$ is essentially free. 

\node{Proposition.}\label{prop:rep_hom_C}  Then, $\wGamma \toi C$ together with $K$ as constructed in~\eqref{node:construct_rep_C} is a representation~\eqref{node:rep_of_hom} of $\H_1^\B(C;\Z)$.
\itnode{Proof.}
For each $b\in B$, the map $x^{-1}R_b \to R_b$ is a finite cover of a disk branched only over $b$. The Galois theory of a punctured disk implies that $x^{-1}R_b$ is a union of disks (with boundary on $\wGamma$). The graph $\wGamma$ together with the components of $x^{-1}R_b$ for each $b$ give $C$ the structure of a CW-complex. The $1$-cells are $\Z[\wGamma_E]$ and the boundary $1$-cells $K$ are precisely those coming from the boundaries of $x^{-1}R_b$'s.
\hfill\qed

\node{Remark.}\label{node:compare_to_bruin} In~\cite{Bruin2019}, the kernel $K$ is computed by finding a lift of the intersection product on $\H_1^\B(C;\Z)$ to $\H_1^\B(|\Gamma|;\Z)$ and then taking $K$ to be the radical of this intersection pairing.

\node{Remark.} There is another strategy for constructing a graph $\Gamma$ in $C$ to determine its homology, which has proven exceptionally efficient for ``superelliptic'' covers of $\p$. Instead of the pullback of a Voronoi diagram in $\p$ around branch points, this construction takes a tree with vertices on the branch locus. The complement of the tree is homeomorphic to a disk and its preimage in $C$ will be (unbranched) copies of that disk. The boundaries of those disks will give the kernel $K$. %

\node{} Now we fix two finite disjoint sets $D,E \subset C(\Qbar)$. %

\node{}\label{node:relative_graph} Start with a rectilinear graph $\Gamma$ in $\p$ which is a deformation retract of $\p(\C)\setminus{}\left(B \cup x(D)\right)$. Enlarge $\Gamma$ by adding, for each point $q$ in $x(E)$, a new vertex at $q$ and a line segment connecting $q$ to the nearest vertex of $\Gamma$ (if $q$ is not already on an edge, otherwise split the edge at $q$). Let $\Gamma'$ be this new graph. Once again, for each $b \in B \cup x(D)$ there is a unique component $R_b$ of $\p(\C)\setminus{}\Gamma'$ containing $b$ and homeomorphic to a disk. 

\node{} Let $\wGamma' \coleq x^{-1} \Gamma'$ be the preimage of $\Gamma'$. 
The pair $(\wGamma',E)$ is a rectilinear embedded graph of $(C\setminus{}D,E)$. Let $K \subset \Z[\wGamma'_E]$ be spanned by the boundaries of the components of $x^{-1} R_b$ for $b \in (B\cup x(D))$ \emph{except} when the component contains a point of $D$, in which case we do not add its boundary to $K$. 

\node{Proposition.}\label{prop:rep_hom} The embedded graph $(\wGamma',E)$ of $(C\setminus{}D,E)$ together with the kernel $K$ constructed above is a \emph{representation} of the homology $\H_1^\B(C\setminus{}D,E;\Z)$.
\itnode{Proof.} 
The proof is analogous to that of Proposition~\eqref{prop:rep_hom_C}. This time, we realized a deformation retract of $C\setminus{}D$ as a CW complex (by omitting the disks containing points of $D$). This completes the proof when $E=\emptyset$.

In general, on $(C\setminus{}D,E)$ we will have a relative cycle connecting any two points of $E$ lying on the same component of $C$. 

We only need to show that the same is true on $(\wGamma',E)$. Which follows if $\wGamma'$ is connected whenever $C$ is connected. But the $\H_0^\B$ of both $\wGamma'$ and $C$ agree as the $2$-cells of $C$ do not change $\H_0$.
\hfill\qed

\node{Remark.}\label{rem:avoid_S} Clearly, we can choose $\Gamma'$ so that it avoids the images of a given finite set $S \subset C(\Qbar)\setminus{}E$ of points. For instance, start with a Voronoi diagram $\Gamma$ around $B \cup x(D) \cup x(S)$, form $\Gamma'$ by adding edges ending in $E$ as before~\eqref{node:relative_graph}. Then $\wGamma'$ is obtained by deleting from $x^{-1}\Gamma'$ any leaves ending in $S$. These leaves appear only when $x(S)$ and $x(E)$ intersect. The preimages of the $2$-cells in $\p$ provide us with the generators for $K$ as before. 

\node{}\label{node:betti_weight_filtration} The weight filtration on $\H_1^\B(C\setminus{}D,E;\Z)$ is as follows
\begin{equation}
\renewcommand{\arraystretch}{1.3}
  \begin{array}{ccc}
    W_{-2} & W_{-1} & W_0 \\ 
    \H_0^\B(D)_0 & \H_1^\B(C\setminus{}D) & \H_1^\B(C\setminus{}D,E)
  \end{array}
\renewcommand{\arraystretch}{1} 
\end{equation}
where $\H_0^\B(D)_0 \coleq \H_0^\B(D)/\H_0^\B(C)$ denotes the primitive homology of $D$ in $C$, and its image in the homology of $C\setminus{}D$ is generated by the loops around $D$. 

\node{} Our graph $(\wGamma',E)$ can be used to determine the weight filtration on $\H_1^\B(\wGamma',E)/K$, simply take $W_{-2}$ to be generated by the loops around $D$ and $W_{-1}$ to be generated by cycles without boundary. 

\node{} This concludes the proof of Proposition~\eqref{prop:rectilinear_basis_computable}.

\subsection{Reduction, pushforward, and transfer maps in homology}

\node{} Let $f \colon C_1 \to C_2$ be a morphism of smooth proper curves over $\Qbar$. For each $i=1,2$, let $D_i, E_i \subset C_i(\Qbar)$ be finite disjoint sets. Recall that $f$ is represented by giving the inclusion of the function fields of respective components.

\node{} If $f^{-1}(D_2) \subset D_1$ and $f(E_1) \subset E_2$ then we have a map $f \colon (C_1\setminus{}D_1, E_1) \to (C_2\setminus{}D_2, E_2)$ which induces the homological pushforward map
\begin{equation}
  f_*^\B \colon \H_1^\B(C_1\setminus{}D_1, E_1) \to \H_1^\B(C_2\setminus{}D_2, E_2).
\end{equation}

\node{Proposition.}\label{prop:effective_hom_pushforward} Suppose that for each $i=1,2$, we are given a rectilinear basis~\eqref{def:rectilinear_basis} for the homology of $(C_i\setminus{}D_i,E_i)$ with respect to a map $x_i \colon C_i \to \p$.  Given a representation of $f \colon (C_1\setminus{}D_1, E_1) \to (C_2\setminus{}D_2, E_2)$, we can effectively compute the integer matrix representing the homological pushforward $f_*^\B$ in the given rectilinear bases. In particular, when $f$ is the identity map, we can compute the change of coordinates between two bases for the homology.
\itnode{Proof.}
Compute~\eqref{prop:effective_cohom_basis} a basis for the cohomology groups $\H^1_\AdR(C_i\setminus{}D_i, E_i)$ and compute~\eqref{prop:coh_pullback_matrix} the matrix representing the cohomological pullback $f^*_\AdR$. 

Integrating the cohomology basis against the rectilinear bases, compute~\eqref{node:approx_integration} a complex ball matrix representing the comparison isomorphisms 
\begin{equation}
  \xc_i \colon \H^1_\AdR(C_i\setminus{}D_i,E_i)\otimes \C \isoto \H^1_\B(C_i\setminus{}D_i,E_i) \otimes \C.
\end{equation}

Compatibility between the pullback maps in cohomology gives the identity 
\begin{equation}\label{eq:pullback_compatibility}
  f^*_\B = \xc_1 \circ f^*_\AdR \circ \xc_2^{-1}.
\end{equation}
Our approximation of the $\xc_i$'s will give a complex ball matrix representing the right hand side of~\eqref{eq:pullback_compatibility} in the basis dual to the rectilinear bases. Increasing the precision of integration if necessary, we may assume the radius of the complex balls are less than $1/2$. The unique integers contained in each of these complex balls recovers the integer matrix representing $f^*_\B$, with transpose~$f_*^\B$.
\hfill\qed

\node{Remark.} Even when $C \coleq C_1 = C_2$ and $f = \id$, if $x_1, x_2 \colon C \to \p$ are different maps, it would be extremely cumbersome to try and find a homological equivalence between two rectilinear bases with respect to $x_1$ and $x_2$. This would involve constructing $2$-chains bounded by ``semi-algebraic'' $1$-chains. The path through cohomology is cleaner.

\node{} Suppose $f \colon C_1 \to C_2$ is such that $f^{-1}(E_2) \subset E_1$ and $f(D_1) \subset D_2$. Then we get the pushforward map
\begin{equation}
  f^*_\B \colon \H_1^\B(C_1\setminus{}E_1,D_1) \to \H_1^\B(C_2\setminus{}E_2,D_2).
\end{equation}
Take duals and use the Poincar\'e--Lefschetz duality $\H_1^\B(C_i\setminus{}E_i,D_i)^\vee \simeq \H_1^\B(C_i\setminus{}D_i,E_i)$ to get the \emph{homological transfer map}
\begin{equation}
  f^!_\B \colon \H_1^\B(C_2\setminus{}D_2,E_2) \to \H_1^\B(C_1\setminus{}D_1,E_1).
\end{equation}
Unwinding the definition, this map can be described by taking the total preimage of chains from $C_2$ to $C_1$. 
The dual of this map is called the \emph{cohomological transfer map} $f_!^\B$. 

\node{Proposition.}\label{prop:effective_hom_transfer} We can effectively compute the integer matrix representing the homological transfer map with respect to given rectilinear homology bases.
\itnode{Proof.}
Use the same outline of proof as~\eqref{prop:effective_hom_pushforward}. Except, we evoke the compatibility of the transfer maps $f_!^\AdR$ and $f_!^\B$.
\hfill\qed

\node{Reduction algorithm for $1$-chains.}\label{node:hom_reduction} Let $\gamma$ be a smooth $1$-chain on $C(\C)\setminus D$ with boundary on $E$. Given a rectilinear basis for $\H_1^\B(C\setminus{}D,E)$, we can effectively compute the integer coordinates of the homology class $[\gamma]$ in the given basis. The only assumption on the representation of $\gamma$ is that we should be able to approximate integrals of rational differentials along $\gamma$ to arbitrary precision.
\itnode{Proof.} 
We may assume $C$ is irreducible. Find a non-Weierstrass $p \in C(\Qbar)$ not lying on $\gamma$ and not on the given rectilinear basis~\eqref{node:weierstrass}. Compute a basis for $\H^1_\AdR(C\setminus{}D,E)$ with poles only on $p$~\eqref{eq:split_cohom},~\eqref{node:basis_of_differentials}. By integration~\eqref{node:approx_integration}, approximate the comparison isomorphism $\comp$ with respect to the bases at hand. Integrate the cohomology basis against $\gamma$ to get a complex ball vector $\int_\gamma \in (\H^1_\AdR \otimes \C)^\vee$. Now apply the inverse of $\comp$ to $\int_\gamma$. The unique integer vector in the resulting complex ball vector containing $\comp^{-1}(\int_\gamma)$ gives the desired coordinates. If the integer vector in the complex ball vector is not unique, increase the precision with which $\comp$ is approximated and repeat.
\hfill\qed 

\newpage
\section{Relations between the periods of $1$-motives}\label{sec:motives}

\node{} This section first introduces the notation we will use when working with $1$-motives. It then proves Theorem~\eqref{thm:expected_relations}, which states that, roughly speaking, if a $1$-motive $M$ is at least as symmetric as the underlying abelian variety, then all relations between the periods of $M$ arise from endomorphisms of $M$ and from trivial relations~\eqref{def:triv_rels}.

\node{} The focus of this section is purely theoretical: we do not concern ourselves with effectivity problems here. Our main result, Theorem~\eqref{thm:expected_relations}, is a reinterpretation of parts of the arguments used in the proof of the ``Dimension estimate''~\cite[Theorem~15.3]{HW}, one of the main theorems of \emph{loc.cit.} Instead of estimating the dimension of the space of periods, we give a direct and explicit description of the space of period relations. We will later use this theoretical description to compute period relations explicitly, see Remark~\ref{rem:how_to_use}.

\subsection{Definition of a $1$-motive}

\node{} Throughout this section, all varieties are defined over $\Qbar$. We set up the basic notions needed to define $1$-motives.

\node{} A \emph{lattice} $L \simeq \Z^{\rk L}$ is a free, finitely generated abelian group. An \emph{algebraic torus} $T$ is an algebraic group such that $T \simeq \G_m^{\dim T}$, where $\G_m = \Spec(\Qbar[t,t^{-1}])$ denotes the multiplicative group.

\node{} A \emph{semi-abelian variety} $G$ is an algebraic group that is an extension of an abelian variety $A$ by an algebraic torus $T$, that is, $0 \to T \to G \to A \to 0$ is exact.

\node{} A \emph{marking} on a semi-abelian variety $G$ is a group homomorphism $\varphi \colon L \to G(\Qbar)$ from a lattice $L$ to the group of $\Qbar$-points of $G$.

\node{} A \emph{$1$-motive} $M$ is a two-term complex $[L \overset{\varphi}{\to} G]$, where $L$ is a lattice, $G$ is a semi-abelian variety, and $\varphi$ is a marking.

\node{} A \emph{morphism} of $1$-motives is a morphism of complexes, i.e., a pair of morphisms between the underlying semi-abelian varieties and lattices that commute with the markings. We denote the category of $1$-motives by $\motz$. Morphisms in $\motz$ are denoted $\hom_\Z$ (``integral morphisms''), and endomorphisms by $\End_\Z$.

\node{} The underlying lattice $L$, the semi-abelian variety $G$, the toric part $T$, and the abelian core $A$ are all functorially assigned to a $1$-motive. When discussing a $1$-motive $M$, we will refer to these pieces by the corresponding letters $L,\, T,\, A$.

\node{} To work up to isogeny, we pass to the category $\mot$, the isogeny category of $\motz$, which is an abelian category~\cite{Deligne1974}. Unless specified otherwise, we work with $1$-motives up to isogeny, i.e., in $\mot$. In particular, $\hom = \hom_\Z \otimes \Q$ and $\End = \End_\Z \otimes \Q$ refer to morphisms and endomorphisms in $\mot$ as opposed to $\motz$. If we want to emphasize that we are in $\mot$, we will write $\hom_\Q$ and $\End_\Q$, in particular for abelian varieties considered up to isogeny.

\subsection{The mixed Hodge structure associated to a $1$-motive}

\node{} Let $\mhsz$ denote the category of $(\Z,\Qbar)$-mixed Hodge structures. Its isogeny category is $\mhs$, the (abelian) category of $(\Q,\Qbar)$-mixed Hodge structures. There is a duality functor $\bullet^\vee = \hom_{\mhs}(\bullet,\Q(0))$, which commutes with the duality functors on (bi)filtered vector spaces.

\node{}\label{node:h1_functor} There is an ``integral homology'' functor
\begin{equation}
\H_1(\bullet;\Z) \colon \motz \to \mhsz,
\end{equation}
which induces a fully faithful functor
\begin{equation}
\H_1 \colon \mot \to \mhs,
\end{equation}
and which restricts to the usual homology functor on semi-abelian varieties. The construction of $\H_1$ over the complex numbers is due to~\cite[(10.1.3)]{Deligne1974}; the generalization to arbitrary characteristic zero fields appears in~\cite{Andre2019}. See also~\cite[Prop.~8.17]{HW} for a proof over $\Qbar$.

\node{} Composing $\H_1$ with the duality functor yields a cohomology functor, which we denote by
\begin{equation}
\H^1(\bullet) \coleq \H_1(\bullet)^\vee.
\end{equation}

\node{} We denote by $\H_1^\B(M)$ the filtered $\Q$-vector space underlying $\H_1(M)$, and by $\H_1^\B(M;\Z)$ the underlying filtered lattice. The bifiltered $\Qbar$-vector space underlying $\H_1(M)$ is denoted by $\H_1^\AdR(M)$. Dually, we use $\H^1_\B$ and $\H^1_\AdR$ for the corresponding components of $\H^1$.

\node{Relation to curves.}\label{node:homology_of_motive_of_curve} Let $(C\setminus{}D,E)$ be a relative curve. There is a functorial $1$-motive $J_{C,E}^D$ associated to $(C\setminus{}D,E)$, whose abelian core is the Jacobian $J_C$ of $C$; see~\eqref{def:jac_mot}. There is a natural isomorphism
\begin{equation}\label{eq:h1_mot_to_curve}
\H^1(J_{C,E}^D;\Z) \simeq \H^1(C\setminus{}D,E;\Z),
\end{equation}
see the proof of Lemma~12.9 in~\cite{HW}. In practice, we will use the identification~\eqref{eq:h1_mot_to_curve} and pushout/pullback constructions to represent the cohomology of motives via the cohomology of curves; see~\S\ref{sec:pushpull}.

\subsection{Subgroup theorem for $1$-motives}

\node{} For a $1$-motive $M$, let $\cf(M) = \H_1^\B(M) \otimes \H^1_\AdR(M)$ denote the tensor product of Betti and algebraic de Rham realizations, called the \emph{space of formal periods} of $M$. We define $\cf_1(M) \subset \cf(M)$ to be the cone of rank-$1$ tensors:
\begin{equation}
\cf_1(M) = \set{\gamma \otimes \omega : \gamma \in \H_1^\B(M), \, \omega \in \H^1_\AdR(M)}.
\end{equation}

\node{Space of period relations.} The comparison isomorphism between Betti and de Rham realizations induces a \emph{period pairing}
\begin{equation}
\wp_M \colon \cf(M) \to \C.
\end{equation}
The kernel of $\wp_M$, denoted $\cR(M)$, is the \emph{space of period relations} of $M$. We define the \emph{cone of rank-$1$ relations} as $\cR_1(M) = \cR(M) \cap \cf_1(M)$.

\node{Functoriality of the period map.}\label{node:wp_func} The period pairing is compatible with morphisms of $1$-motives. If $f \colon M_1 \to M_2$ is a $\Q$-morphism, then for $\gamma \in \H_1^\B(M_1)$ and $\omega \in \H^1_\AdR(M_2)$ we have
\begin{equation}
\wp_{M_2}(f_*\gamma \otimes \omega) = \wp_{M_1}(\gamma \otimes f^* \omega),
\end{equation}
where $f_*$ and $f^*$ denote the induced maps on Betti and de Rham realizations, respectively.

\node{}\label{node:ses} Consider a short exact sequence
\begin{equation}\label{eq:ses}
0 \to M' \overset{\iota}{\to} M \overset{\jmath}{\to} M'' \to 0,
\end{equation}
and elements $\gamma' \in \H_1^\B(M')$, $\omega'' \in \H^1_\AdR(M'')$. Then we have
\begin{equation}
\wp_M(\iota_*\gamma' \otimes \jmath^*\omega'') = \wp_{M''}(\jmath_*\iota_*\gamma' \otimes \omega'') = \wp_{M''}(0\otimes \omega'') = 0.
\end{equation}
Therefore, $\iota_*\gamma' \otimes \jmath^*\omega''$ belongs to $\cR_1(M)$.

\node{} Let us denote by $\cR_\ses(M) \subset \cR_1(M)$ the cone of all relations $\iota_*(\gamma') \otimes \jmath^*(\omega'')$ arising from short exact sequences as in~\eqref{node:ses}. The following is one of the key results of~\cite{HW}.

\node{Theorem.}\label{thm:ses}~(\cite[Theorem 9.7]{HW}) There is an equality $\cR_\ses(M) = \cR_1(M)$. \qed

\node{The case of abelian varieties.}\label{node:abelian_strategy} When $M = [0 \to A]$ is an abelian variety, all short exact sequences with $M$ in the middle split up to isogeny. Consequently, short exact sequences involving $A$ in $\mot$ are parametrized by the endomorphism algebra $\End(A)$, a finite-dimensional $\Q$-algebra. By Theorem~\eqref{thm:ses}, the space $\cR_1(A)$ can be computed using $\End(A)$ and its action on homology and cohomology.

\node{} In contrast to abelian varieties, not every short exact sequence of $1$-motives splits up to isogeny. Therefore, the computation of period relations will be more involved. We turn to this problem next.

\subsection{Cartier duality}\label{sec:cartier_duality}

\node{}\label{node:cartier} Cartier duality is a contravariant duality functor $M \mapsto M^\vee$ on $\motz$, satisfying a natural equivalence $M \simeq (M^\vee)^\vee$~\cite[(10.2.11)]{Deligne1974}.

\node{}\label{node:delignes_equivalences} Cartier duality on $1$-motives agrees with a Tate twist of duality on mixed Hodge structures. More precisely, there is a natural equivalence
\begin{equation}
\H_1(M^\vee;\Z) \simeq \Hom_{\mhsz}(M,\Z(1)),
\end{equation}
see~\cite[(10.2.10)]{Deligne1974}.

\node{} The Cartier dual of an algebraic torus $T$ is the motive $[\Xi(T) \to 0]$, where $\Xi(T) = \hom(T,\G_m) \simeq \Z^{\dim T}$ is the \emph{character lattice} of $T$.

\node{} The Cartier dual of a lattice motive $[L \to 0]$ is the algebraic torus $\G_m^L = \hom(L,\G_m) \simeq \G_m^{\rk L}$.

\node{}\label{node:cartier_dual_of_G_easy} If $G$ is a semi-abelian variety with toric part $T \simeq \G_m$ and abelian core $A$, then $G$ can be viewed as a torsor under $\G_m$ over $A$, defining an element of $\H^1(A,\co_A^\times)$, and hence a point $p_G$ on $A^\vee = \Pic^0(A)$. The Cartier dual of $G$ is the $1$-motive
\begin{equation}
G^\vee = [\Z \overset{1 \mapsto p_G}{\to} A^\vee].
\end{equation}

\node{}\label{node:cartier_dual_of_G} In general, the Cartier dual of a semi-abelian variety $G$, extending $A$ by a torus $T$, is the $1$-motive
\begin{equation}
G^\vee = [\Xi(T) \to A^\vee],
\end{equation}
where the marking $\Xi(T) \to A^\vee(\Qbar)$ is defined as follows: for each $\ell \in \Xi(T)$, the pushout $\ell_*G$ is an extension of $A$ by $\G_m$, and the associated point $p_{\ell_*G} \in A^\vee(\Qbar)$ defines the image of $\ell$.

\node{} By~\eqref{node:cartier_dual_of_G}, and the fact that Cartier duality is an equivalence, we see that it establishes a contravariant equivalence between $1$-motives of the second and third kind.

\node{} If $M = [L \to G]$ is a $1$-motive with toric part $T$ and abelian core $A$, then the Cartier dual is a $1$-motive of the form
\begin{equation}
M^\vee = [\Xi(T) \to [L\to A]^\vee],
\end{equation}
where $[L\to A]^\vee$ is the Cartier dual of $[L \to A]$. Moreover, the quotient of $M^\vee$ by its toric part $\G_m^L$ recovers $G^\vee = [\Xi(T) \to A^\vee]$.

\node{Explicit Cartier duality for curves.}\label{node:explicit_duality} Combining~\eqref{node:delignes_equivalences}, Poincaré--Lefschetz duality~\eqref{rem:poincare_lefschetz}, and the fact that $\H_1$ is fully faithful~\eqref{node:h1_functor}, we see that the $1$-motive $J_{C,E}^D$ associated to $(C\setminus{}D,E)$ is Cartier dual to the motive $J_{C,D}^E$ associated to $(C\setminus{}E,D)$. This will allow us to represent duals explicitly.

\subsection{Endomorphism types of motives}

\node{} An abelian variety is called a $1$-motive of the \emph{first kind} ($T=0$ and $L=0$). An abelian variety with a marking is a $1$-motive of the \emph{second kind} ($T=0$). A semi-abelian variety is a $1$-motive of the \emph{third kind} ($L=0$).

\node{Definition.}\label{def:reduced} A $1$-motive $M$ of the second or third kind is \emph{reduced} if and only if the natural map $\End(M) \to \End(A)$ is injective, where $A$ is the abelian core of $M$. A general $1$-motive $M = [L \to G]$ is reduced if and only if its underlying second and third type motives, $M/T = [L \to A]$ and $G$, are reduced.

\node{Definition.} A $1$-motive $M$ with abelian core $A$ is \emph{super-saturated} if the natural map $\End(M) \to \End(A)$ is surjective. If $M$ is reduced and super-saturated, then it is called \emph{saturated}.

\node{Definition.}\label{def:baker} A motive is a \emph{Baker motive} if its abelian core is trivial, $A=0$.

\node{Definition.}\label{def:composite} A motive is \emph{composite} if it is a direct sum of a Baker motive and a saturated motive.

\node{Remark.} Our definition of reduced agrees with that of Huber and Wüstholz~\eqref{lem:reduced_is_lattice_inj}. It is clear that our definition of saturated agrees with theirs.

\node{Lemma.}\label{lem:reduced_is_lattice_inj} A $1$-motive of the second kind $M=[L \to A]$ is reduced if and only if the induced map 
\begin{equation}\label{eq:induced_marking}
  L\otimes_\Z \Q \to A(\Qbar)\otimes_\Z \Q 
\end{equation}
is injective. A $1$-motive of the third kind $G$ is reduced if and only if its Cartier dual $G^\vee = [\Xi(T) \to A^\vee]$ is reduced.

\itnode{Proof.} If~\eqref{eq:induced_marking} is injective, then any $\Q$-endomorphism of $A$ that lifts, lifts uniquely to $L\otimes \Q$. Therefore, injectivity of the $\Q$-marking implies injectivity of $\End(M) \to \End(A)$ for a second kind motive. Conversely, if there is a nontrivial kernel in the $\Q$-marking, then nonzero $\Q$-endomorphisms of the kernel can be extended to $M$ in a way that maps to zero on $A$. For a third kind motive, apply the same argument to the Cartier dual.  
\hfill\qed

\node{Remark.} A motive $M$ can satisfy $\End(M) \simeq \End(A)$ and still be non-reduced if $M$ is not of the second or third kind. For instance, consider a simple abelian variety $A$ with $\End(A) \simeq \Q$, and the motive $[\Z \to \G_m \oplus A]$ where we mark a non-torsion point in each component. This motive has no nontrivial endomorphisms, but the underlying semi-abelian variety $\G_m \oplus A$ is non-reduced.

\subsection{Expected period relations}

\node{} Let $M = [L \to G]$ be a $1$-motive with toric part $T$ and abelian core $A$.

\node{Definition.} Consider the short exact sequence $T \toi M \tos [L \to A]$, and identify $\H_1^\B(T)$ with its image in $\H_1^\B(M)$. We define the \emph{toric relations} as
\begin{equation}
\cR_T(M) \coleq \ker \left( \H_1^\B(T) \otimes \H^1_\AdR(M) \tos \H_1^\B(T) \otimes \H^1_\AdR(T) \overset{\wp_T}{\to} \C \right).
\end{equation}

\node{} Note that $\H_1^\B(T) \otimes \H^1_\AdR([L\to A])$ is contained in $\cR_T(M)$, as it is the kernel of the projection onto $\cf(T)$. The period map $\wp_M$ restricted to $\H_1^\B(T) \otimes \H^1_\AdR(M)$ factors through $\wp_T$ by functoriality~\eqref{node:wp_func} applied to $T \toi M$.

\node{Definition.} Consider the short exact sequence $G \toi M \tos L$, and identify $\H^1_\AdR(L)$ with its image in $\H^1_\AdR(M)$. We define the \emph{lattice relations} as
\begin{equation}
\cR_L(M) \coleq \ker \left( \H_1^\B(M) \otimes \H^1_\AdR(L) \tos \H_1^\B(L) \otimes \H^1_\AdR(L) \overset{\wp_L}{\to} \C \right).
\end{equation}
Observe that $\H_1^\B(G) \otimes \H^1_\AdR(L)$ is contained in $\cR_L(M)$.

\node{Definition.}\label{def:triv_rels} We define the \emph{trivial relations} $\cR_\triv(M)$ as $\cR_\triv(M) \coleq \cR_L(M) + \cR_T(M)$.

\node{Remark.}\label{rem:trivial} The space of trivial relations is straightforward to compute. For the toric relations, the pairing $\wp_T$ is induced from the duality pairing $\Xi(T)^\vee \otimes \Xi(T) \to \Z$, tensored with $\Q\langle \tpi \rangle$ under the canonical identifications~\eqref{node:cohom_of_T_and_xi}. The quotient leading up to $\wp_T$ is simply the quotient by $W_1\H^1_\AdR(M)$.  Similarly, for the lattice relations, we consider the $\Qbar$-tensor of the duality pairing $L \otimes L^\vee \to \Z$ and the quotient by $W_{-1}\H_1^\B(M)$. See also~\eqref{node:nontrivial_weight_pieces}.

\node{Definition.}\label{def:endo_relations} We define the \emph{endo-relations} $\cR_{\End}(M)$ as the subspace
\begin{equation}
\cR_{\End}(M) \coleq \set{\phi_*\gamma \otimes \omega - \gamma \otimes \phi^* \omega : \phi \in \End(M),\; \gamma \in \H_1^\B(M),\; \omega \in \H^1_\AdR(M)},
\end{equation}
which is finite-dimensional and contained in $\cR(M)$ by~\eqref{node:wp_func}.

\node{Definition.}\label{def:exp_rels} We define the \emph{expected relations} of $M$ as
$\cR_e(M) \coleq \cR_{\End}(M) + \cR_\triv(M)$.

\node{Theorem.}\label{thm:expected_relations} If $M$ is a saturated motive, then the space of period relations coincides with the expected relations, that is, $\cR(M) = \cR_e(M)$.

\node{} The proof of this theorem will occupy the rest of this section. We also state the following corollary, which extends the theorem to composite motives; see~\eqref{prop:rels_of_baker_and_sat} for a more explicit statement.

\node{Corollary.}\label{cor:rels_of_Baker_and_sat} If $M$ is a composite motive $M_B \oplus M^{\sat}$, then $\cR(M)$ is generated by the trivial relations, the cross-terms~\eqref{node:crossterms} arising from the direct sum decomposition, the Baker relations $\cR(M_B)$, and the expected relations of $M^{\sat}$.

\itnode{Proof.} This is the content of~\eqref{prop:rels_of_baker_and_sat}, where we substitute $\cR_e(M^{\sat})$ for $\cR(M^{\sat})$ in light of Theorem~\eqref{thm:expected_relations}.
\hfill\qed

\node{Remark.}\label{rem:how_to_use} Recall that the determination of $\cR_\triv(M)$ is straightforward~\eqref{rem:trivial}. We will later show that the period relations of a Baker motive can be computed explicitly~\eqref{lem:baker_period_relations}. The expected relations of a saturated motive can also be computed directly, \emph{provided} that the action of the endomorphism algebra on (co)homology is known. If, further, one  a direct sum decomposition $\H_1(M) = \H_1(M_B) \oplus \H_1(M^{\sat})$, then $\cR(M)$ can be computed using Corollary~\eqref{cor:rels_of_Baker_and_sat}.

\subsection{Filtrations and period matrices}

\node{} Let $M = [L \to G]$ be a $1$-motive with toric part $T$ and abelian core $A$.

\node{}\label{node:motivic_weight_filtration} The \emph{motivic weight filtration} on $M$ is the ascending filtration
\begin{equation}
0 = W_{-3}M \subset W_{-2}M = T \subset W_{-1}M = G \subset W_0M = M.
\end{equation}
It induces the weight filtration on Betti homology, satisfying $\H_1^\B(W_k M) = W_k\H_1^\B(M)$.  
The \emph{motivic weight cofiltration} $Q_{-k-1}M \coloneqq M/W_kM$ gives
\begin{equation}
M \to [L \to A] \to [L \to 0] \to 0,
\end{equation}
and induces the weight filtration on de Rham cohomology, satisfying $\H^1_\AdR(Q_j M) = W_j\H^1_\AdR(M)$.

\node{} A basis $v_1,\dots,v_n$ for a filtered vector space $V$ with filtration $V_\bullet$ is said to \emph{respect the filtration} if for each $i$, the first $\dim V_i$ elements $v_1,\dots,v_{\dim V_i}$ form a basis of $V_i$.

\node{} The short exact sequences
\begin{equation}
0 \to W_kM \to M \to Q_{-k-1}M \to 0
\end{equation}
induce \emph{weight relations} in the period relations:
\begin{equation}
W_k\H_1^\B(M) \otimes W_{-k-1}\H^1_\AdR(M) \subset \cR(M).
\end{equation}

\node{Period matrix and block form.} Choose bases $\gamma_1,\dots,\gamma_m$ for $\H_1^\B(M)$ and $\omega_1,\dots,\omega_m$ for $\H^1_\AdR(M)$, both respecting the weight filtrations.  
The period isomorphism is represented by the \emph{period matrix}
\begin{equation}
\cp_M = \left( \int_{\gamma_j} \omega_i \right)_{i,j=1,\dots,m}.
\end{equation}
By the weight relations, the matrix $\cp_M$ has block form
\begin{equation}
\cp_M = 
\begin{pmatrix}
0 & 0 & \cp_L \\
0 & \cp_A & \cq_{AL} \\
\cp_T & \cq_{TA} & \cq_{TL}
\end{pmatrix},
\end{equation}
where $\cp_T$, $\cp_A$, and $\cp_L$ are period matrices for $T$, $A$, and $L$, respectively.

\node{}\label{node:nontrivial_weight_pieces} The only nontrivial weight relations in the period relations induced by the motivic weight filtration are $\H_1^\B(T) \otimes \H^1_\AdR([L\to A])$ and $\H_1^\B(G) \otimes \H^1_\AdR(L)$, and both are contained in the trivial relations $\cR_\triv(M)$.  In particular, the vanishing of the three $0$ blocks in the period matrix $\cp_M$ follows from the trivial relations.

\node{Toric and lattice periods.} We refer to the entries of $\cp_L$ as \emph{lattice periods} and to the entries of $\cp_T$ as \emph{toric periods}. The matrix $\cp_L$ consists of algebraic numbers, while $\cp_T$ consists of algebraic multiples of $\tpi$.

\subsection{Period relations of Baker motives}

\node{}\label{node:cohom_of_T_and_xi} If $T$ is an algebraic torus with character lattice $\Xi(T)$, then there are canonical isomorphisms
\begin{align}
\H_1^\B(T) &= \Xi(T)^\vee \otimes_\Z \H_1^\B(\G_m) = \Xi(T)^\vee \otimes_\Z \Q, \label{eq:Tbetti} \\ 
\H^1_\AdR(T) &= \Xi(T) \otimes_\Z \H^1_\AdR(\G_m) = \Xi(T) \otimes_\Z \Qbar.  \label{eq:Tadr}
\end{align}
The isomorphism~\eqref{eq:Tadr} is obtained by pulling back $\dd z/z$ from $\G_m$ along characters, and the isomorphism~\eqref{eq:Tbetti} is obtained by pushing forward loops along characters.  
With respect to dual bases for $\Xi(T)$ and $\Xi(T)^\vee$, the period matrix $\cp_T$ is a $\tpi$-multiple of the identity matrix.

\node{} If $L$ is a lattice, viewed as the $1$-motive $M = [L \to 0]$, then there are canonical isomorphisms
\[
\H_1^\B(L) = L \otimes_\Z \Q, \quad \H^1_\AdR(L) = L^\vee \otimes_\Z \Qbar.
\]
With respect to dual bases for $L$ and $L^\vee$, the period matrix $\cp_L$ is the identity matrix.

\node{} If $M = [L \to T]$ is a Baker motive, then with respect to bases compatible with the weight filtrations, the period matrix $\cp_M$ has block form
\begin{equation}\label{eq:baker_periods}
\cp_M = 
\begin{pmatrix}
0 & \cp_L \\
\cp_T & *
\end{pmatrix},
\end{equation}
where $\cp_L$ is a matrix of algebraic numbers, $\cp_T$ is a matrix of algebraic multiples of $\tpi$, and $*$ denotes a matrix whose entries are logarithms of algebraic numbers.

\subsection{Independence statements for periods of saturated motives} 

\node{} If $M$ is saturated and the abelian core $A$ is simple, we say that $M$ is a \emph{simple motive}. Then $E \coloneqq \End(M) \simeq \End(A)$ is a skew field. The weight filtration and cofiltration on $M$, being functorial, are stable under the action of $E$. Thus, the weight filtered $\Q$-vector space $\H_1^\B(M)$ inherits a weight filtration as an $E$-vector space. Let $\gamma_\bullet$ be an $E$-basis for $\H_1^\B(M)$ and $\omega_\bullet$ a $\Qbar$-basis for $\H^1_\AdR(M)$, both respecting the weight filtrations. We define the \emph{reduced period matrix} $\widetilde\cp_M$ as the matrix representing the period pairing with respect to these bases:
\begin{equation}\label{eq:reduced_periods}
\widetilde\cp_M = \left( \int_{\gamma_j} \omega_i \right)_{i,j} =
\begin{pmatrix}
0 & 0 & \widetilde\cp_{L} \\
0 & \widetilde\cp_A & * \\
\widetilde\cp_T & * & *
\end{pmatrix}.
\end{equation}

\node{} We now recall the key technical ingredient behind the ``Dimension estimate'' theorem~\cite[Theorem~15.3]{HW}, adapted to our notation.

\node{Lemma }\cite[Lemma~15.6]{HW}\textbf{.}\label{lem:period_independence} 
Except for the \emph{trivial relations}~\eqref{def:triv_rels}, there are no nonzero $\Qbar$-linear relations among the entries of the reduced period matrices $\widetilde\cp_L$, $\widetilde\cp_A$, $\widetilde\cp_T$, and the $*$ blocks. Here, the trivial relations refer specifically to the relations among entries of $\widetilde\cp_T$, a matrix with entries in $\Qbar\langle \tpi \rangle$, and among entries of $\widetilde\cp_L$, a matrix with entries in $\Qbar$. \hfill\qed

\node{} A stronger statement holds, with an analogous proof; both the statement and the method of proof are indicated in~\cite[Propositions~15.9 and 15.20]{HW}.

\node{Lemma.}\label{lem:no_other_rels}
Let $M_1$ and $M_2$ be $1$-motives with respective abelian cores $A_1$ and $A_2$. Assume that $M_1$ is saturated and that $\hom_\Q(A_1,A_2) = 0$. Then we have
\begin{equation}
\im \wp_{M_1} \cap \im \wp_{M_2} \subset \Qbar\langle 1, \tpi \rangle.
\end{equation}
Moreover, letting $L_1$ and $T_1$ denote the lattice and toric parts of $M_1$, respectively:
\begin{enumerate}
\item If $L_1 = 0$, then $\im \wp_{M_1} \cap \Qbar = 0$.
\item If $T_1 = 0$, then $\im \wp_{M_1} \cap \Qbar\langle \tpi \rangle = 0$.
\end{enumerate}

\itnode{Sketch of proof.}
By Lemmas~15.22 and 15.24 of~\cite{HW}, we may replace $M_2$ by a direct sum of a saturated motive and a Baker motive without changing its abelian core, and only enlarging its space of periods.  
The desired statement follows by adapting the argument of Lemma~15.6 of \emph{loc.cit.}, specifically by examining submotives of a suitably constructed auxiliary motive, as indicated in Propositions~15.9 and 15.20.
\hfill\qed

\node{} Suppose $M$ is composite, that is, a direct sum $M = M_B \oplus M^{\sat}$ where $M_B$ is a Baker motive and $M^{\sat}$ is a saturated motive. The direct sum decomposition induces a decomposition of Betti homology and de Rham cohomology:
\begin{equation}
\H_1^\B(M) = \H_1^\B(M_B) \oplus \H_1^\B(M^{\sat}), \qquad 
\H^1_\AdR(M) = \H^1_\AdR(M_B) \oplus \H^1_\AdR(M^{\sat}).
\end{equation}
Thus, the space of formal periods $\calF(M) = \H_1^\B(M) \otimes \H^1_\AdR(M)$ decomposes naturally as the direct sum of four components:
\begin{equation}\label{eq:direct_sum_periods}
\calF(M) = \calF(M_B) \bigoplus \calF(M^{\sat}) \bigoplus 
\big( \H_1^\B(M_B) \otimes \H^1_\AdR(M^{\sat}) \big) \bigoplus
\big( \H_1^\B(M^{\sat}) \otimes \H^1_\AdR(M_B) \big).
\end{equation}

\node{}\label{node:crossterms} The period relations $\cR(M_B)$ and $\cR(M^{\sat})$ are embedded into $\cR(M)$ via the identifications $\calF(M_B) \subset \calF(M)$ and $\calF(M^{\sat}) \subset \calF(M)$. Furthermore, the \emph{cross-terms} $\H_1^\B(M_B) \otimes \H^1_\AdR(M^{\sat})$ and $\H_1^\B(M^{\sat}) \otimes \H^1_\AdR(M_B)$ are contained in $\cR(M)$: they pair a homology class from one summand with a cohomology class from the other, and thus vanish under the period pairing.

\node{Proposition.}\label{prop:rels_of_baker_and_sat} We have an equality of period relations:
\begin{equation*}
\cR(M) = \cR_\triv(M) + \cR(M_B) + \cR(M^{\sat}) + 
\H_1^\B(M_B) \otimes \H^1_\AdR(M^{\sat}) +
\H_1^\B(M^{\sat}) \otimes \H^1_\AdR(M_B).
\end{equation*}

\itnode{Proof.} The right-hand side is contained in $\cR(M)$ by construction.  
The equality asserts that there are no further relations among the periods of $M_B$ and $M^{\sat}$.  
This follows from Lemma~\eqref{lem:no_other_rels}, since $M^{\sat}$ is saturated and the abelian core of $M_B$ is zero.
\hfill\qed

\subsection{Period relations of a power of a simple motive}\label{sec:simple_powers}

\node{Proposition.}\label{prop:simple} If $M$ is a saturated motive with simple abelian core, then the space of period relations and the expected relations coincide, i.e., $\cR(M) = \cR_e(M)$.
\itnode{Proof.}
In a suitable basis, the period map $\wp_M \colon \cf(M) \to \C$ is represented by the vector of entries of the period matrix $\cp_M$. Reducing modulo the endo-relations, we can consider $\widetilde\wp_M \colon \cf(M)/\cR_{\End}(M) \to \C$ which is represented by the vector of entries of $\widetilde \cp_M$ as given in~\eqref{eq:reduced_periods}. By Lemma~\eqref{lem:period_independence}, the only remaining relations are coming from the zero blocks and the relations between $\widetilde \cp_T$ and $\widetilde \cp_L$. These remaining relations are contained in the trivial relations. 
\hfill\qed

\node{} \label{node:almost_simple} Let $A = B^n$ where $B$ is a simple abelian variety and $n\ge 1$. For the skew field $E = \End(B)$, the endomorphism algebra $\End(A)$ is isomorphic to the algebra $\Mat_n(E)$ of $n \times n$ matrices over $E$. Let $m_{ij} \in \End(A)$ denote the $n\times n$ matrix with a $1$ on the $(i,j)$-th entry and zeros everywhere else. Equivalently, let $p_i$ be the projection onto the $i$-th factor, $\iota_i$ the identification of $B$ onto the $i$-th component, so that $m_{ij} = \iota_i \circ p_j$. 

\node{} 
\label{node:almost_simple_motive}
If $M$ is a saturated $1$-motive with abelian core $A = B^n$ as in~\eqref{node:almost_simple}, the identification $\End(M) \simeq \End(A)$ furnishes $M$ with the endomorphisms algebra $\Mat_n(E)$. Consequently, $M$ must be of the form $(M')^n$ where $M'$ has abelian core $B$. The operators $p_i$, $\iota_i$, $m_{ij}$ lift to $M$ in the natural way.

\node{}
Let $M = (M')^n$ as in~\eqref{node:almost_simple_motive}. Then, we can choose bases so that the period matrix of $M$ is block diagonal,
\begin{equation}\label{eq:almost_simple_block}
  \cp_M = \left( 
  \begin{array}{cccc}
    \cp_{M'}       & 0        & \dots  & 0 \\
    0              & \cp_{M'} & \dots  & 0 \\
    &          & \ddots &   \\
    0              & 0        &        & \cp_{M'}
  \end{array}
  \right) .
\end{equation}

\node{Proposition.}\label{prop:intermediate} For $M$ as above, the space of period relations coincides with the space of expected relations.
\itnode{Proof.}
  We will use the elements $m_{ij} \in \End(M) \simeq \Mat_n(E)$ from~\eqref{node:almost_simple}. Proposition~\eqref{prop:simple} implies that the entries of $E \cdot m_{11}$ give all the non-trivial relations amongst the entries of the first copy of $\cp_{M'}$ in~\eqref{eq:almost_simple_block}. There can be no other relations except those suggested by the repeated blocks. We will show that these obvious relations can all be accounted by the endo-relations. %

  We will first show that endomorphisms give enough relations to deduce that the $i$-th diagonal block is equal to the $1$-st diagonal block. Let $\gamma \in \H_1^\B(M')$ and $\omega \in \H^1_\AdR(M')$. We will use the identities $\iota_i = m_{i1} \circ \iota_1$, $p_i = p_1 \circ m_{1i}$, $m_{1i} \circ m_{i1} \circ \iota_1 = \iota_1$ to deduce the following identities modulo the endomorphism relations:
  \begin{equation}
    \iota_{i,*}\gamma \otimes p_i^* \omega = m_{i1,*}\iota_{1,*}\gamma \otimes m_{1i}^*p_1^*\omega
    \equiv m_{1i,*}m_{i1,*}\iota_{1,*}\gamma\otimes p_1^* \omega 
    = \iota_{1,*} \gamma \otimes p_1^* \omega.
  \end{equation}

  To see that the $(i,j)$-th block is zero for $i \neq j$, use the identities $m_{jj} \circ \iota_i = 0$ and $p_j \circ m_{jj} = p_j$. Then, modulo endomorphism relations, we have
  \begin{equation}
    \iota_{i,*} \gamma \otimes p_j^* \omega = \iota_{i,*} \gamma \otimes m_{jj}^* p_j^* \omega  
     \equiv m_{jj,*}\iota_{i,*} \gamma \otimes p_j^* \omega 
     = 0.
  \end{equation}
  The trivial relations account for the remaining relations between the lattice and toric periods. 
\hfill\qed

\subsection{Period relations for saturated motives}\label{sec:rels_of_sat_mots}

\node{} We will decompose a saturated motive as indicated in~\cite[\S15.2.2]{HW} and read off the relations in terms of endomorphisms.

\node{}\label{node:general_ab_var}
  Let $A$ be a non-zero abelian variety. Then $A$ decomposes (up to isogeny) as 
  \begin{equation}
    A = B_1^{n_1} \oplus \dots \oplus B_m^{n_m}
  \end{equation}
  where $B_i$'s are simple, $n_i \ge 1$, and $\hom_\Q(B_i,B_j) = 0$ for $i \neq j$. Let $E_i = \End(B_i)$ so that
  \begin{equation}\label{eq:endo_decomposition}
  \End(A) \simeq \Mat_{n_1}(E_1) \oplus \dots \oplus \Mat_{n_m}(E_m).
\end{equation}

\node{}\label{node:decompose_saturated}
  Let $M$ be a non-zero saturated $1$-motive with abelian core $A$. Necessarily, $A$ is non-zero and decomposes as in~\eqref{node:general_ab_var}. Since $\End(M) \simeq \End(A)$, $M$ has a similar structure
  \begin{equation}
    M = (M_1')^{n_1} \oplus \dots \oplus (M_m')^{n_m}, 
  \end{equation}
  with the abelian core of $M_i'$ being simple $B_i$. Let us write $M_i = (M_i')^{n_i}$ so that $M_i$ belongs to the case considered in~\S\ref{sec:simple_powers}.

\node{}
  Let $\jmath_i \colon M_i \toi M$ be the inclusion of the $i$-th component into $M$ and $q_i \colon M \tos M_i$ be the projection onto the $i$-th component. For each $i,j = 1,\dots,m$ we define the following element of $\End(M)$:  
  \begin{equation}
    \xm_{ij} = \begin{cases}
      0 & i \neq j \\
      \jmath_i \circ q_i 
    \end{cases}.
  \end{equation}
  We have $q_j \circ \xm_{jj} = q_j$ and $\xm_{jj} \circ \jmath_i = 0$ if $i \neq j$. 

\node{Remark.} 
    Since $\hom(B_i,B_j) = 0$ for $i \neq j$ the same holds for the saturated motives $M_i$ and $M_j$. The maps $\xm_{ij} = 0$ for $i \neq j$ correspond to projecting onto the $j$-th factor and then mapping as $0$ to the $i$-th component.

\node{}
  The period matrix for $M$ is of the form
  \begin{equation}\label{eq:saturated_block}
    \cp_M = \left( 
    \begin{array}{cccc}
      \cp_{M_1} & 0         & \dots  & 0 \\
      0         & \cp_{M_2} & \dots  & 0 \\
                &           & \ddots & \\
      0         & 0         &        & \cp_{M_m}
    \end{array}
    \right). 
  \end{equation}

\itnode{Proof of Theorem~\ref{thm:expected_relations}.}
    All relations within the blocks $\cp_{M_i}$ are contained in the expected relations via Proposition~\eqref{prop:intermediate}. 

    The zero blocks in~\eqref{eq:saturated_block} can be explained by the endomorphism relations. The proof of this fact is similar to the argument at the end of Proposition~\eqref{prop:intermediate}, but this time we use the identities $q_j \circ \xm_{jj} = q_j$ and $\xm_{jj} \circ \jmath_i = 0$ if $i \neq j$. 

    There are no further relations between the periods in distinct blocks $\cp_{M_i}$ and $\cp_{M_j}$, $i \neq j$, except for the trivial relations, by Lemma~\eqref{lem:no_other_rels}.
\hfill\qed

\newpage
\section{Push-pull Jacobian motives}\label{sec:pushpull}

\node{} In this section, we describe a mild form of pushout and pullback constructions on $1$-motives, sufficient to realize a given $1$-motive as the push-pull of a composite motive. This reduction allows us to apply Theorem~\eqref{thm:expected_relations} to the composite and compute period relations for the original motive via~\eqref{prop:pushpull_relations}. We carry curves, motives, and their mixed Hodge structures effectively along these constructions without losing track of periods or relations.

\node{} To a curve $(C\setminus{}D,E)$ we associate a Jacobian motive~\eqref{def:jac_mot}, a $1$-motive with abelian core $J_C$, the Jacobian of $C$. Our push-pull constructions only change the lattice and the torus around this core.

\node{} There are three results in this section. The first is an explicit description of the ``points'' of the push-pull of a Jacobian motive in terms of restricted divisor classes on a curve~\eqref{prop:points_of_pushpull}. The second is an explicit description, in terms of the MHS of a curve, of the MHS of the push-pull Jacobian motive~\eqref{prop:MHS_of_pushpull}. The third result states that if a motive $M_1$ is realized as the push-pull of another motive $M_2$ then the period relations $\cR(M_1)$ of the former can be computed from the period relations $\cR(M_2)$ of the latter by linear algebra on (co)homology groups~\eqref{prop:pushpull_relations}.

\subsection{Jacobian motive of a punctured relative curve}

\node{} Let $C/\Qbar$ be a smooth proper curve and let $D, E \subset C(\Qbar)$ be disjoint finite sets. The group of invertible rational functions on $C$ are $\kappa(C)^\times$ and we will write $\kappa^\times$ for the invertible locally constant functions on $C$.

\node{} We define the multiplicative subgroup
\begin{equation}
  \kappa(C)^\times_D \coleq \set{f \in \kappa(C)^\times : f|_D \in \G_m^D} \subset \kappa(C)^\times
\end{equation}
of rational functions without zeros or poles on $D$ and the subgroup
\begin{equation}
  \kappa(C)^\times_{D,1} \coleq \set{f \in \kappa(C)^\times : f|_D = 1} \subset \kappa(C)^\times_D
\end{equation}
of functions constant $1$ along $D$.

\node{} The group of divisors (over $\Qbar$) on $C\setminus{}D$ is denoted by $\Div(C\setminus{}D)$. Let $\Div^0(C\setminus{}D)$ be the subgroup of divisors on $C\setminus{}D$ of multidegree $0$ on $C$, i.e., of degree $0$ on each component.

\node{}\label{node:LS} For a set $S \subset C(\Qbar)$, denote by $L_S \subset \Div^0(C)$ the lattice of multidegree $0$ divisors supported on $S$.

\node{} There exists~\cite{Serre1988} a semi-abelian variety $J_C^D$, the \emph{Jacobian of $(C\setminus{}D)$}, defined over $\Qbar$ whose $\Qbar$-points coincide with the quotient
\begin{equation}
  \frac{\Div^0(C\setminus{}D)}{\set{\divv(f) : f \in \kappa(C)^\times_{D,1}}}.
\end{equation}
Note the restricted use of rational functions in the quotient. See~\eqref{lem:shape_of_JD} for the realization of $J_C^D$ as a toric extension of the Jacobian of $C$.

\node{} The \emph{Abel--Jacobi map} is the quotient map 
\begin{equation}
  \AJ \colon \Div^0(C\setminus{}D) \to J_C^D(\Qbar) : \xi \mapsto \xi \mod \kappa(C)^\times_{D,1}.
\end{equation}

\node{Definition.}\label{def:jac_mot}  The \emph{Jacobian motive of $(C\setminus{}D,E)$} is $\JJ \coleq [L_E \overset{\AJ}{\to} J_C^D]$, where the the marking is given by restricting the Abel--Jacobi map $\AJ$ to $L_E \subset \Div^0(C\setminus{}D)$, divisors of multidegree $0$ supported on $E$~\eqref{node:LS}.

\node{} When either $D$ or $E$ is empty, we will drop it from notation. In particular, $J_C$ is the Jacobian of the (smooth proper) curve $C$. 

\node{Remark.} We considered it would be jarring to refer to $J_C^D$ and $\JJ$ of $J_C$ as Albanese varieties or complexes but, technically, that is how we defined them. The usual Jacobian $J_C$ is both the Albanese and the $\Pic^0$ of $C$. On the other hand, $J_C^D$ is the Albanese of $C\setminus{}D$ but the $\Pic^0$ of the ``nodal curve'' $(C,D)$, where $D$ is contracted to a single node with transversal branches~\cite{Serre1988}. The complex $[L_D \to J_C]$ is the Albanese of the ``nodal'' curve $(C,D)$~\cite{BS99} and $\Pic^0$ of $C\setminus{}D$, morally standing for the quotient of $J_C$ by $L_D$. More generally, the Jacobian motive $\JJ$ is best viewed as the Albanese of the punctured ``nodal'' curve $(C\setminus{}D,E)$, although it is also the $\Pic^0$ of $(C\setminus{}E,D)$.

\subsection{Push-pull constructions}

\node{} Let $M = [L \to G]$ be a $1$-motive with toric part $T$ and abelian core $A$.

\node{} Let $\psi \colon Y \to L$ be a morphism of lattices. Then $\psi^*M \to M$, the \emph{$\psi$-pullback of $M$}, is defined by the pullback diagram
\begin{equation}
\begin{tikzcd}
	G & M & L \\
	G & {\psi^*M} & Y
	\arrow[hook, from=1-1, to=1-2]
	\arrow[two heads, from=1-2, to=1-3]
	\arrow["{\operatorname{id}_G}", from=2-1, to=1-1]
	\arrow[hook, from=2-1, to=2-2]
	\arrow[from=2-2, to=1-2]
	\arrow["\ulcorner"{anchor=center, pos=0.125, rotate=90}, draw=none, from=2-2, to=1-3]
	\arrow[two heads, from=2-2, to=2-3]
	\arrow["\psi"', from=2-3, to=1-3]
\end{tikzcd}
\end{equation}
We will often write $\psi^*M$ for the pullback, leaving the map $\psi^*M \to M$ (which includes $Y \to L$) implicit. 

\node{}\label{node:pullback_as_kernel} Let $\lambda \colon M \to L = M/G$ denote the projection map. The pullback $\psi^*M$ is canonically isomorphic to the kernel of the map $(\lambda - \psi) \colon M \oplus Y \to L$. This is a general construction valid in any abelian category.

\node{} Let $\chi \colon X \to \Xi(T)$ a morphism of lattices. Recall that the Cartier dual $M^\vee$ of $M$ has lattice part $\Xi(T)$. The $\chi$-pullback of $M^\vee$ is $\chi^*M^\vee \to M^\vee$. Taking Cartier duals gives the \emph{$\chi$-pushout of $M$}, $M \to \chi_*M$, where $\chi_*M \coleq (\chi^*M^\vee)^\vee$. As with pullbacks, we will often write $\chi_*M$ for the pushout, leaving the map implicit.

\node{} An alternative description of the $\chi$-pushout is as follows. The map $\chi$ induces the map $\widehat\chi \colon T \to \G_m^X$ of tori, where $\G_m^X = \hom(X,\G_m)$ is an algebraic torus with dimension equal to the rank of $X$. Then the $\chi$-pushout of $M$ is the pushout of $M$ via $\widehat\chi$, 
\begin{equation}
\begin{tikzcd}
	T & M & {[L\to A]} \\
	{\mathbb{G}_m^X} & {\chi_*M} & {[L\to A]}
	\arrow[hook, from=1-1, to=1-2]
	\arrow["{\widehat\chi}"', from=1-1, to=2-1]
	\arrow[two heads, from=1-2, to=1-3]
	\arrow[from=1-2, to=2-2]
	\arrow["{\operatorname{id}}", from=1-3, to=2-3]
	\arrow[hook, from=2-1, to=2-2]
	\arrow["\ulcorner"{anchor=center, pos=0.125, rotate=180}, draw=none, from=2-2, to=1-1]
	\arrow[two heads, from=2-2, to=2-3]
\end{tikzcd}
\end{equation}

\node{}\label{node:pushout_as_cokernel} Let $\tau \colon T \to M$ be the inclusion map. The pushout $\chi_*M$ is the cokernel of the map $(\widehat\chi, - \tau) \colon T \to \G_m^X \oplus M$. This is a general construction valid in any abelian category. 

\node{Definition.} Given morphisms of lattices $\chi \colon X \to \Xi(T)$ and $\psi \colon Y \to L$ then the \emph{push-pull} of $M$ with respect to $(\chi,\psi)$ is $\chi_*\psi^* M$. 

\node{} The universal properties of the pushout and pullback give a canonical isomorphism $\chi_*\psi^*M \simeq \psi^*\chi_*M$. We will implicitly identify the two.%

\subsection{Points of a push-pull Jacobian motive}

\node{} Consider two morphisms of lattices $\chi \colon X \to L_D$ and $\psi \colon Y \to L_E$. We will describe the \emph{push-pull Jacobian motive} $\ppJ \coleq \chi_*\psi^*\JJ$ using divisors on $C$, see Proposition~\eqref{prop:points_of_pushpull}. To be more precise, we will describe the points of the underlying semi-abelian variety and the marking, but this distinction will be blurred.

\node{} To make the pushout clearer, we need to give an alternative construction of the Jacobian motive of $(C\setminus{}D,E)$.

\node{} For any $f \in \kappa(C)^\times_D$, the restriction of $f$ to $D$ defines an element $f|_D \in \G_m^D = \hom(D,\G_m)$. For any divisor $\xi = \sum a_i p_i$ supported on $D$, define $f(\xi) = \prod f(p_i)^{a_i}$. In particular, this defines an element $f|_D \in \G_m^{L_D} = \hom(L_D,\G_m)$. Note that the locally constant functions $\kappa^\times$ map to $1 \in \G_m^{L_D}$ since $L_D$ consists of multidegree $0$ divisors.

\node{}\label{node:kappa_map} We have a map $\kappa(C)^\times_D \to \G_m^{L_D} \oplus \Div^0(C\setminus{}D) : f \mapsto (f|_D,\divv(f))$, which restricts to the map $\kappa(C)^\times_{D,1} \to 1 \oplus \Div^0(C\setminus{}D)$, where $1$ stands for the trivial group. In particular, by definition,
\begin{equation}
\JJ = \frac{1 \oplus \Div^0(C\setminus{}D)}{\kappa(C)^{\times}_{D,1}}.
\end{equation}

\node{Lemma.}\label{lem:surj_to_JC} There is a natural isomorphism
\begin{equation}
  \frac{\Div^0(C\setminus{}D)}{\kappa(C)^\times_D} \isoto J_C
\end{equation}
\itnode{Proof.}
Every divisor on $C$ is linearly equivalent to a divisor on $C\setminus{}D$. 
\hfill\qed

\node{Lemma.}\label{lem:shape_of_JD} The Jacobian of $(C\setminus{}D)$ fits into an exact sequence
\begin{equation}
  0 \to \G_m^{L_D} \to J_C^D \to J_C \to 0.
\end{equation}
\itnode{Proof.}
Surjectivity is by Lemma~\ref{lem:surj_to_JC}. The kernel of $J_C^D \to J_C$ is the quotient
\begin{equation}
  \frac{\set{\divv(f) : f \in \kappa(C)^\times_D}}{\set{\divv(f) : f \in \kappa(C)^\times_{D,1}}} %
\end{equation}
which is easily seen to be the torus $\G_m^{L_D}=\G_m^D/\kappa^\times$. 
\hfill\qed

\node{Lemma.}\label{lem:JCD_alternative} The inclusion 
\begin{equation}
  1 \oplus \Div^0(C\setminus{}D) \toi \G_m^{L_D} \oplus \Div^0(C\setminus{}D)
\end{equation}
induces an isomorphism
\begin{equation}
  J_C^D \isoto \frac{\G_m^{L_D}\oplus \Div^0(C\setminus{}D)}{\kappa(C)^\times_D}.
\end{equation}
\itnode{Proof.}
Both sides surject onto $J_C$ with kernel $\G_m^{L_D}$. The proof of Lemma~\eqref{lem:shape_of_JD} makes it clear that the map induced from the left to the right induces the identity on $J_C$ and $\G_m^{L_D}$. Now use the five lemma.
\hfill\qed

\node{}\label{node:alternative_marking} Lemma~\eqref{lem:JCD_alternative} gives an alternate construction of the Jacobian motive $J_{C,D}^E = [L_E \overset{\AJ}{\to} J_C]$. Using the isomorphism described in Lemma~\eqref{lem:JCD_alternative} we see that the marking for the alternative description of $J_C^D$ is $\xi \mapsto (1,\xi) \mod \kappa(C)^\times_{D,1}$. %

\node{}\label{node:lattice_pullback_of_funcs} The lattice morphism $\chi \colon X \to L_D$ induces a map $\kappa(C)^\times_D \to \G_m^X : f \mapsto f|_D \circ \chi$. Therefore, modifying~\eqref{node:kappa_map} we get a map
\begin{equation}
  \kappa(C)^\times_D \to \G_m^X \oplus \Div^0(C\setminus{}D) : f \mapsto (f|_D \circ \chi,-\divv(f)).
\end{equation}

\node{Lemma.}\label{lem:pushout_points} Let $\chi \colon X \to L_D$ be a morphism of lattices. The pushout $\chi_*J_C^D$ is naturally isomorphic to
\begin{equation}
  \frac{\G_m^X \oplus \Div^0(C\setminus{}D)}{\kappa(C)^\times_D}.
\end{equation}
\itnode{Proof.}
This is clear from the alternative construction of $J_C^D$ in Lemma~\eqref{lem:JCD_alternative}. To elaborate, it is convenient to view (the points of) both $J_C^D$ and $\chi_*J_C^D$ as pushouts of 
\begin{equation}
  0 \to \kappa(C)^\times_D/\kappa^\times \to \Div(C\setminus{}D) \to J_C \to 0
\end{equation}
via the restriction maps to $\G_m^{L_D}$ and $\G_m^X$ respectively.
\hfill\qed

\node{Proposition.}\label{prop:points_of_pushpull} Let $\chi \colon X \to L_D$ and $\psi \colon Y \to L_E$ be morphisms of lattices. The push-pull Jacobian motive $\ppJ$ is given by
\begin{equation}
 Y \to \frac{\G_m^X \oplus \Div^0(C\setminus{}D)}{\kappa(C)^\times_D} : y \mapsto (1,\psi(y)) \mod \kappa(C)^\times_D.
\end{equation}
\itnode{Proof.}
This combines Lemma~\eqref{lem:pushout_points} and the observation made in~\eqref{node:alternative_marking}.
\hfill\qed

\subsection{A moving lemma}

\node{Definition.}\label{def:divisor_class} 
For $\xi \in \Div^0(C\setminus{}D)$, let $[\xi]_D = \{\xi + \divv(f) : f \in \kappa(C)_{D,1}\}$ denote the \emph{$D$-linear equivalence class of $\xi$}. This equivalence class represents a point of $J_C^D$. When we have a lattice morphism $\chi \colon X \to L_D$ and $a \in \G_m^X$, then we let $[a,\xi]_\chi = \{(a-f|_D\circ\chi, \xi + \divv(f))\}$ be the \emph{$\chi$-linear equivalance class of $(a,\xi)$}, which defines a point of $J_C^\chi$. When $a=1$ is the identity, we write $[\xi]_\chi \coleq [1,\xi]_\chi$, making the two notations compatible when $\chi=\id_{L_D}$.

\node{} For $D \subset D'$ have a surjection $J_C^{D'} \tos J_C^{D}$. This implies that any $D$-linear equivalence class $[\xi]_D$ can be represented by a divisor $\xi'$ supported outside of $D'$. We need to make the determination of $\xi'$ from $\xi$ effective.

\node{Lemma.}\label{lem:move} Let $D \subset D' \subset C(\Qbar)$ be finite sets. Then any divisor $\xi \in \Div(C \setminus D)$ can be effectively moved away from $D'$ within its $D$-linear equivalence class.

\itnode{Proof.} We may assume that $\xi$ is supported only on $D' \setminus D$, leaving the rest in place. Write $\xi = \xi_+ - \xi_-$ where $\xi_+,\, \xi_-$ are effective and have disjoint support. Let $\xi_-^\epsilon$ be obtained by increasing each coefficient in $\xi_-$ by $1$.

Let $\widetilde{D} = D' \setminus (D \cup \supp(\xi))$. Fix a divisor $\xp$ of positive degree on every component of $C$, and choose $n \gg 0$ such that
\begin{equation}
  h^1(n\xp - \xi_-^\epsilon - \widetilde{D} - D) = 0.
\end{equation}
Consider the exact sequence
\begin{equation}
  0 \to \co_C(n\xp - \xi_-^\epsilon - \widetilde{D} - D) \to \co_C(n\xp + \xi) \to \co_{\xi_+} \oplus \co_{\supp(\xi_-)} \oplus \co_{\widetilde{D} \cup D} \to 0.
\end{equation}
By our choice of $n$, the global sections of the middle term surject onto those of the right term. Hence, there exists $f \in \kappa(C)^\times$ such that
\begin{equation}
  \divv(f) + \xi + n\xp \ge 0
\end{equation}
and
\begin{equation}
  \divv(f) - \xi_-^\epsilon - \widetilde{D} - D + n\xp \not\ge 0.
\end{equation}
In particular, $f$ has just enough zeros to cancel $\xi_-$, but no extra zeros on $\widetilde{D}$ or $D$. Surjectivity onto $\co_{\xi_+}$ allows us to arrange poles cancelling $\xi_+$. Finally, surjectivity to $\co_D$ ensures we can arrange $f|_D = 1_D$.

Thus, $\xi + \divv(f)$ is $D$-linearly equivalent to $\xi$ and supported away from $D'$.
\hfill\qed

\node{Remark.} This construction gives an explicit algorithm: given $\xi \in \Div(C \setminus D)$, it produces $f \in \kappa(C)_{D}^\times$ such that $\xi' \coleq \xi + \divv(f) \in \Div^0(C \setminus D')$. It follows that, for any lattice morphism $\chi \colon X \to L_D$, the divisor class $[\xi]_\chi$ can be represented by $[f|_D\circ \chi,\xi']_\chi$. 

\subsection{Homology and cohomology of a push-pull Jacobian motive}

\node{} Consider two morphisms of lattices $\chi \colon X \to L_D$ and $\psi \colon Y \to L_E$. We wish to describe the (co)homology of the push-pull motive $\ppJ \coleq \chi_*\psi^*\JJ$.

\node{} Recall that we have a canonical isomorphism $\H^1(\JJ) \simeq \H^1(C\setminus{}D,E)$~\eqref{eq:h1_mot_to_curve}. We will combine our explicit representation of $\H^1(C\setminus{}D,E)$~\eqref{node:representation_of_MHS}, with the push-pull constructions below to give an explicit representation of $\H^1(\ppJ)$.

\node{} The $\H_1,\, \H^1$ functors from $\mot$ to $\mhs$ are fully-faithful and exact. Therefore, the homology and cohomology of the push-pull motives can be described by making use of the descriptions of the pullback as a kernel~\eqref{node:pullback_as_kernel} and the pushout as a cokernel~\eqref{node:pushout_as_cokernel}.

\node{} For a lattice $L$ we will write $L_\Q \coleq L \otimes_\Z \Q$ and $L_\Qbar \coleq L \otimes_\Z \Qbar$. Similarly, use subscripts for the $\Q$ or $\Qbar$ tensor of maps between lattices.

\node{} We will need to refer to parts of the usual exact sequence involving the (co)homology of $(C\setminus{}D,E)$. The quotient map to the lattice $\lambda \colon \JJ \tos L_E$ and the inclusion of the torus $\tau \colon \G_m^{L_D} \toi \JJ$ induce the following maps 
\begin{align}
  \partial \coleq \lambda_\B &\colon \H_1^\B(C\setminus{}D,E) \to L_{E,\Q}, \\
  \Gy \coleq \tau_\B &\colon L_{D,\Q}^\vee \to \H_1^\B(C\setminus{}D,E), \\
  \delta \coleq \lambda^\adr &\colon L_{E,\Qbar}^\vee \to \H^1_\adr(C\setminus{}D,E), \\
  \Res \coleq \tau^\adr &\colon \H^1_\adr(C\setminus{}D,E) \to L_{D,\Qbar}.
\end{align}
The map $\partial \coleq \lambda_\B$ is the usual boundary map on chains and $\Res \coleq \tau^\adr$ is the residue map on differential forms. The map $\Gy\coleq\tau_\B$ is the Gysin map, putting circles around the punctures $D$ of $C$. The map $\delta \coleq \lambda^\adr$ is the connecting morphism in cohomology, corresponding to the inclusion of constants on $E$, see~\eqref{eq:ses_for_differentials} and~\eqref{eq:split_cohom}, after the obvious identification $L_{E,\Qbar}^\vee = \Qbar_0^E$~\eqref{node:qbarE}.

\node{Lemma.} For a morphism $\chi \colon X \to L_D$ of lattices, we have
\begin{equation}
  \H^1_\adr(J_C^\chi) = \set{(\eta,x) \in \H^1_\adr(C\setminus{}D) \oplus X_\Qbar : \Res(\eta) = \chi(x)}.
\end{equation}
\itnode{Proof.}
Applying $\H^1_\adr$ to $\G_m^{L_D} \toi J_C^D \tos J_C$ and making use of the identification~\eqref{eq:h1_mot_to_curve} we get 
\begin{equation}
  0 \to \H^1_\adr(C) \to \H^1_\adr(C\setminus{}D) \overset{\Res}{\to} L_D \otimes \Qbar \to 0.
\end{equation}
The pullback of this exact sequence via $\chi_\Qbar$ will give the cohomology of $J_C^\chi$ by functoriality of $\H^1$ and~\eqref{node:pullback_as_kernel}.
\hfill\qed

\node{} The general case can be proven in the same manner. We will simply state the result. 

\node{Proposition.}\label{prop:MHS_of_pushpull} Let $\chi \colon X \to L_D$ and $\psi \colon Y \to L_E$ be morphisms of lattices. The mixed Hodge structre on the push-pull Jacobian motive $\ppJ$ is given by:
\begin{align}
  \H^1_\adr(\ppJ) &= \frac{\set{(s,\eta,x) \in Y^\vee_\Qbar \oplus \H^1_\adr(C\setminus{}D,E) \oplus X_\Qbar : \Res(\eta) = \chi(x)}}{\im(-\psi_\Q^\vee,\delta,0)}, \\
  \H_1^\B(\ppJ) &= \frac{\set{(t,\gamma,y) \in X^\vee_\Q \oplus \H_1^\B(C\setminus{}D,E) \oplus Y_\Q : \partial \gamma = \psi(y)}}{\im(-\chi_\Q^\vee,\Gy,0)},
\end{align}
and, by denoting the class of triplets by square brackets, the period pairing 
\begin{equation}
  \wp_{\ppJ}([s,\eta,x] \otimes [t,\gamma,y]) = y(s) + \int_\gamma \eta + \tpip \, t(x),
\end{equation}
where we used the integration pairing~\eqref{node:period_pairing}. The weight and Hodge filtrations are deduced immediately from those of the (co)homology of~$(C\setminus{}D,E)$. \hfill \qed

\subsection{Pushing and pulling period relations}

\node{} We are going to show in this subsection that if $M$ is the push-pull of a motive $M'$ then $\cR(M)$ can be computed from $\cR(M')$ using just linear algebra on (co)homology spaces. 

\node{}\label{node:pullback_relation_setup} Suppose $M = \psi^* M' \overset{f}{\to} M'$ is a pullback by a lattice map $\psi$. Let $I^\B \coleq \im f_*^\B \subset \H_1^\B(M')$ and $I_\adr \coleq \im f^*_\adr \subset \H^1_\adr(M)$ denote the images of the induced maps on (co)homology. Choose a splitting of the surjections 
\begin{equation}
  \H_1^\B(M) \tos I^\B, \quad \H^1_\adr(M') \tos I_\adr,
\end{equation}
and denote the lifts of the images by $\widetilde I^\B \subset \H_1^\B(M)$ and $\widetilde I_\adr \subset  \H^1_\adr(M')$. This gives an isomorphism
\begin{equation}
  \sigma \colon I^\B \otimes \widetilde I_\adr \isoto \widetilde I^\B \otimes I_\adr.
\end{equation}

\node{Lemma.}\label{lem:pullback_relations} %
We have $\cR(M) = \cR_\triv(M) + \ker(f_*^\B)\otimes \im(f^*_\adr) + \sigma\left(\cR(M') \cap (I^\B \otimes \widetilde I_\adr) \right)$.
\itnode{Proof.}
Let $K=\ker(f) \subset M$ be the kernel of $f$. Notice that $K$ is the lattice $[\ker \psi \to 0]$ and, therefore, has periods in $\Qbar$. Let $I = \im(f) \subset M'$ be the image of $f$. The lattice part of $I$ is $\im(\psi)$ and if we pick a splitting of the lattice of $M$ as $\ker(\psi) \oplus \im(\psi)$ then we can embed $I$ into $M$ by identifying $I$ with the submotive of $M$ obtained from the pullback via $\im(\psi) \toi \ker(\psi) \oplus \im(\psi)$. Call this submotive $\oI \subset M$. We now have a decomposition $M = K \oplus \oI$ which allows us to view $\cf(\oI)$ and $\cf(K)$ as subspaces of $\cf(M)$.

The period matrix $\cp_M$ of $M$ is of the block diagonal form
\begin{equation}\label{eq:block_M_IK}
  \cp_M = 
  \begin{blockarray}{ccc}
                         & \ker(f_*^\B) & \H_1^\B(\oI) \\
    \begin{block}{c(cc)}
    \H^1_\AdR(K)         & \cp_K        & 0  \\
    \ker(f^*_\adr)       & 0            & \cp_I \\
    \end{block}
  \end{blockarray}\vspace{-2ex}
\end{equation}
The lower left zero block in the matrix $\cp_M$ corresponds to $\ker(f_*^\B)\otimes \im(f^*_\adr) \subset \cR(M)$. Let $M = [L \to G]$.

The space $\H_1^\B(K) \simeq \ker(f_*^\B) \subset \H_1^\B(M)$ is a lift of $\ker(\psi)_\Q \subset \H_1^\B(L)$ via the quotient $\H_1^\B(M) \to \H_1^\B(L)$. Modulo $\H_1^\B(G)\otimes \H^1_\AdR(L) \subset \cR_\triv(M)$, the space $\ker(f_*^\B) \otimes \H^1_\AdR(L) \subset \cf(M)$ is independent of the lift and the period pairing on $\ker(f_*^\B) \otimes \H^1_\AdR(L)$ factors through $\wp_L$. Therefore, $(\ker(f_*^\B) \otimes \H^1_\AdR(L)) \cap \cR(M)$ is contained in the trivial relations.

The submotive $\oI \subset M$ contains the underlying semi-abelian variety $G = W_{-1}M$, regardless of the splitting chosen for $\im\psi$. Therefore, the projection $M \to K$ factors through $M \to L$. In particular, $\H^1_\AdR(K) \subset \H^1_\AdR(L) \subset \H^1_\AdR(M)$, regardless of the splitting. 

Moreover, modulo $\ker(f_*^\B)\otimes \im(f^*_\adr) \subset \cR(M)$, the space $\ker(f_*^\B) \otimes \H^1_\AdR(K) \subset \cf(M)$ is independent of our choice of splitting.

All combined, we see that the space of relations arising from the entries of $\cp_K$ are contained in $\cR_\triv(M) + \ker(f_*^\B)\otimes \im(f^*_\adr)$. This explains all relations in the first column of~\eqref{eq:block_M_IK}.

Every zero entry in the upper right zero block in $\cp_M$ is coming from trivial relations: Since, $G \subset \oI$ and $\H^1_\AdR(K) \subset \H^1_\AdR(L)$ the weight relation $\H_1^\B(G) \otimes \H^1_\AdR(L)$ explains a portion of this zero block. The splitting $L = \ker(\psi) \oplus \im(\psi)$ identifies $\H^1_\AdR(K)$ with $\im(\psi)^\perp_\Qbar \subset \H^1_\AdR(L) = L_\Qbar$. Modulo $\H_1^\B(G)$, $\H_1^\B(\oI)$ is identified with $\im(\psi)_\Q \subset L_\Q$. This orthogonal pairing explains the rest of the zero block.

Furthermore, since $\cp_K$ consists of algebraic numbers and algebraic periods of any motive (modulo relations) must come from its lattice part, we conclude that any relations between the periods of $\cp_\oI$ and $\cp_K$ must lie in the space $\cR(\oI) + \cR_\triv(M) \subset \cR(M)$.  %

So far, we proved 
\begin{equation}\label{eq:so_far}
\cR(M) = \cR_{\triv}(M) + \ker(f_*^\B)\otimes \im(f^*_\adr) + \cR(\oI). 
\end{equation}
Observe, $\H^1_\AdR(\oI)= \im f^*_\AdR = I_\AdR$ is independent of the splitting. As for $\H_1^\B(\oI)$, it is a lift of $I_\B$. The other lift $\widetilde I_B$ will differ by elements of $\ker(f_*^\B)$ so
\begin{equation}
  \widetilde I^\B \otimes I_\AdR \equiv \H_1^\B(\oI) \otimes I_\AdR \mod \ker(f_*^\B) \otimes \im(f^*_\AdR).
\end{equation}
Therefore, in the expression~\eqref{eq:so_far}, we can replace $\cR(\oI)$ with $\cR(M) \cap (\widetilde I^\B \otimes I_\AdR)$, which is equal to
\begin{equation}
 \ker(\wp_M|_{I^\B \otimes \widetilde I_\adr}).
\end{equation}

Via the isomorphism $\widetilde I^\B \isoto I^\B$ we can identify the period pairing $\wp_M$ restricted to $\widetilde I^\B \otimes I_\AdR$ with the period pairing of $\oI \simeq I$. Of course, the period pairing on $I \subset M'$ is obtained analogously by the period pairing of $M'$ restricted to $I^\B \otimes \widetilde I_\AdR$. The kernel $\cR(\oI)$ of the period pairing on $\oI$ can now be identified with with 
\begin{equation}
  \cR(M') \cap (I^\B \otimes \widetilde I_\adr).
\end{equation}
The map $\sigma$ makes the identifications explicit.
\hfill\qed

\node{} Suppose $M' \overset{g}{\to} \chi^* M' = M$ is a pushout by a map $\chi$ character lattices. Let $I^\B \coleq \im f_*^\B \subset \H_1^\B(M)$ and $I_\adr \coleq \im f^*_\adr \subset \H^1_\AdR(M')$. Once again, choose lifts $\widetilde I^\B \subset \H_1^\B(M')$ and $\widetilde I_\adr \subset \H^1_\AdR(M)$ of $I^\B$ and $I_\adr$ respectively. This gives an isomorphism 
\begin{equation}
  \sigma \colon \widetilde I^\B \otimes I_\adr \isoto I^\B \otimes \widetilde I_\adr.
\end{equation}

\node{Lemma.}\label{lem:pushout_relations} We have $\cR(M) = \cR_\triv(M) + \im(f_*^\B) \otimes \ker(f^*_\adr) + \sigma\left(\cR(M') \cap (\widetilde I^\B \otimes I_\adr)\right)$.
\itnode{Proof.}
One can replicate the proof of Lemma~\eqref{lem:pullback_relations}, this time making reference to the toric periods and toric relations instead. Reference to algebraic periods must be replaced by algebraic multiples of $\tpi$. 
\hfill\qed

\node{Proposition.}\label{prop:pushpull_relations} If $M$ is the push-pull of a motive $M'$ then explicit linear algebra on the (co)homology spaces computes $\cR(M)$ from $\cR(M')$.
\itnode{Proof.}
Let $M$ be the push-pull of $M'$ via the pair of lattice morphisms $(\chi,\psi)$. Let $M'' = \psi^*M' \to M'$ be the $\psi$-pullback. Compute $\cR(M'')$ from $\cR(M')$ via Lemma~\eqref{lem:pullback_relations}. Now, $M'' \to \chi_*M'' = M$ so that $\cR(M)$ can be computed from $\cR(M'')$ via Lemma~\eqref{lem:pushout_relations}.

The statement of the referenced lemmas require only the computation of kernels and images of vector space maps as well as making arbitrary choices to lift quotients of vector spaces. Recall that the trivial relations are also computed by pure linear algebra~\eqref{rem:trivial}. 
\hfill\qed

\newpage
\section{Correspondences}\label{sec:corr}

\node{} In this section, we recall the classical theory of correspondences acting on Jacobian varieties of smooth proper curves, and the representation of the endomorphism algebra $\End_\Z(J_C)$ by correspondences. We also recall that the action of a correspondence on divisors and homology can be computed effectively. We extend these constructions to Jacobian motives associated to marked and punctured curves.

\node{} Correspondences will be crucial for representing the action of $\End_\Z(J_C)$ on Jacobian motives explicitly. In particular, they allow us to compute the endo-relations introduced in~\eqref{def:endo_relations}. They also form the basis for the supersaturation process in~\S\ref{sec:supsat}. 

\node{} We will also recall the practical methods recently developed to compute a representation of $\End_\Z(J_C)$~\cite{CMSV18,CLV21}. 

\subsection{Correspondences on Jacobian varieties}

\node{} Let $C_1, \, C_2$ be smooth proper irreducible curves over $\Qbar$. 

\node{Correspondences.} An \emph{irreducible correspondence from $C_1$ to $C_2$} is a triple $(C, f_1, f_2)$ where $C$ is an irreducible curve and $f_i \colon C \to C_i$, $i=1,2$, are finite morphisms. A finite formal $\Z$-linear combination $\calC = \sum a_j (C^{(j)},f_1^{(j)},f_2^{(j)})$ of irreducible correspondences is a \emph{correspondence}. The additive group of correspondences is $\Corr(C_1,C_2)$.

\node{Remark.} A correspondence from $C_1$ to $C_2$ defines a divisor in $C_1\times C_2$ with no horizontal or vertical components. Conversely, such a divisor will yield a correspondence by taking the sum of its resolved irreducible components together with the two projection maps. These two are \emph{functionally} equivalent, the latter being the more common way to define correspondences. We chose a definition that is more compatible with our representations of curves as function fields.

\node{Representing a correspondence.}\label{node:represent_correspondence} An irreducible correspondence is represented by the induced inclusions of function fields $\kappa(C_i) \to \kappa(C)$ for $i=1,2$. A correspondence is represented by the formal $\Z$-linear combination of representations of its irreducible components $(C^{(j)},f_1^{(j)},f_2^{(j)})$.

\node{Transpose of a correspondence.} The \emph{transpose} of an irreducible correspondence $(C,f_1,f_2)$ is the irreducible correspondence $(C,f_2,f_1)$. This operation is extended linearly to correspondences $\Corr(C_1,C_2) \isoto \Corr(C_2,C_1)$ and we denote the transpose of $\calC$ by $\calC^t$. 

\node{Action on divisors.} An irreducible correspondence $\calC = (C,f_1,f_2)$ induces a morphism $\calC_* \colon \Div(C_1) \to \Div(C_2)$ defined by $\xi \mapsto f_{2,*}f_1^*\xi$, namely, by pulling back divisors via $f_1$ and pushing forward via $f_2$. By taking $\Z$-linear combinations, any correspondence $\calC$ induces a morphism $\calC_* \colon \Div(C_1) \to \Div(C_2)$. We will also write $\calC^*$ to denote the transpose $\calC^t_* \colon \Div(C_2) \to \Div(C_1)$.

\node{Linear equivalence.} For $\calC = (C,f_1,f_2)$ we have $\calC_*\divv(g) = \divv(\Nm_{f_2}(f_1^*(g)))$. In particular, any correspondence $\calC$ induces an action on divisor classes which we will denote by $[\calC_*] \colon J_{C_1} \to J_{C_2}$.

\node{Action on points is effective.}\label{node:correspond_on_pts} Given a representation of $\calC$ and $\xi \in \Div(C_1)$ represented as a formal sum of points, the computation of $\calC_*(\xi) \in \Div(C_2)$ is effective since the pullback and pushforward maps on divisors are effective. 

\node{}\label{node:correspond_on_H1} An irreducible correspondence $\calC = (C,f_1,f_2)$ induces a morphism of MHS $ \calC_* \coleq f_{2,*} \circ f_1^! \colon \H_1(C_1) \to \H_1(C_2)$. This operation extends linearly to any correspondence. %

\node{Lemma.}\label{lem:correspondence_action_on_MHS} Given bases for the Betti and de Rham realizations of $\H_1(C_i)$ the action of any correspondence $\calC_*$ can be effectively computed as matrices in the given bases.
\itnode{Proof.} %
 First, compute bases for the homology and cohomology of $C$ as in \S\ref{sec:adr} and \S\ref{sec:betti}. Then apply Propositions~\eqref{prop:coh_pullback_matrix},~\eqref{prop:coh_transfer_matrix},~\eqref{prop:effective_hom_pushforward},~\eqref{prop:effective_hom_transfer} to compute the pullback, pushforward, and transfer maps. 
\hfill\qed

\node{} The group $\hom_\Z(J_{C_1}, J_{C_2}) \simeq \Z^\rho$ is torsion-free and of finite rank. It is classical that each morphism $J_{C_1} \to J_{C_2}$ arises from a correspondence, hence the natural map $\Corr(C_1, C_2) \to \hom_\Z(J_{C_1}, J_{C_2})$ is surjective.

\node{}\label{node:rep_hom} A \emph{representation} of $\hom_\Z(J_{C_1}, J_{C_2})$ consists of correspondences $\calC_1, \dots, \calC_\rho \in \Corr(C_1, C_2)$ such that the induced maps $\mu_i = [\calC_{i,*}]$ form a $\Z$-basis for $\hom_\Z(J_{C_1}, J_{C_2})$.

\node{} If $C_1$ and $C_2$ are reducible, with irreducible components $\{C_{1,j}\}$ and $\{C_{2,j}\}$ respectively, then $\Corr(C_1, C_2)$ is defined as the direct sum of correspondences between components: $\Corr(C_1, C_2) = \bigoplus_{i,j} \Corr(C_{1,i}, C_{2,j})$. The group $\hom_\Z(J_{C_1}, J_{C_2})$ admits the analogous decomposition and its representation is defined as in~\eqref{node:rep_hom}.

\node{}\label{node:reducible_C_correspondences} Let $C = \coprod_i C_i$ be a smooth proper curve over $\Qbar$ with $C_i$'s the irreducible components of $C$. We view $\Corr(C) \coleq \Corr(C, C) = \bigoplus \Corr(C_i, C_j)$ as a matrix of correspondences. Similarly, the endomorphism ring decomposes as $\End_\Z(J_C) = \bigoplus \hom_\Z(J_{C_i}, J_{C_j})$.

\node{} A \emph{representation} of $\End_\Z(J_C)$ consists of a basis of correspondences $\calC_1,\dots,\calC_\rho \in \Corr(C)$ as in~\eqref{node:rep_hom} \emph{and} a tensor $A \in \Z^{\rho \times \rho \times \rho}$ such that $[\calC_{i,*}] \cdot [\calC_{j,*}] = \sum_{k} A_{ij}^k [\calC_{k,*}] \in \End_\Z(J_C)$. 

\node{Remark.} In other words, we represent $\End_\Z(J_C)$ by an abstract $\Z$-algebra together with a $\Z$-linear map to $\Corr(C)$ such that the composition with $\Corr(C) \to \End_\Z(J_C)$ yields an isomorphism of $\Z$-algebras.

\node{Effectivity.}\label{node:correspondence_is_effective} Given a representation of $\End_\Z(J_C)$, its action on $J_C(\Qbar)$---with points represented by divisor classes---and on $\H_1(C)$ can be computed effectively, using~\eqref{node:correspond_on_pts} and~\eqref{node:correspond_on_H1}.

\subsection{Computing a representation of the endomorphism algebra of a Jacobian}

\node{} Costa, Mascot, Sijsling, and Voight~\cite{CMSV18} give an algorithm to compute a representation of $\Hom_\Z(J_1, J_2)$ from the equations of two irreducible curves $C_1$ and $C_2$, \emph{provided} the rank $a = \rk \Hom(J_1, J_2)$ is known. When $a$ is unknown, the method becomes an algorithm only under the Mumford--Tate conjecture~\cite{CLV21}. The procedure is implemented and performs well in practice.

\node{} While the algorithm of~\cite{CMSV18} is formulated for computing $\End_\Z(J_C)$ when $C$ is irreducible over $\Qbar$, it also applies to $\Hom_\Z(J_{C_1}, J_{C_2})$ for irreducible curves $C_1$ and $C_2$. For a possibly reducible curve $C$, use the identification in~\eqref{node:reducible_C_correspondences} to reduce to the irreducible case and compute a representation of $\End_\Z(J_C)$.

\node{Theorem.}~(\cite{CMSV18, CLV21}) There exists a (practical, implemented) semi-algorithm to compute $\End_\Z(J_C)$ from the equations of a smooth proper curve $C/\Qbar$. This becomes an algorithm assuming the Mumford--Tate conjecture, or unconditionally if $\rk \End_\Z(J_C)$ is known. \hfill\qed

\node{Remark.} In practice, the semi-algorithm is expected to terminate even when the rank is not known in advance. When it does, the output is provably correct. A practical upper bound for the rank, that we expect to be sharp, can be obtained from the crystalline cohomology of a prime reduction of the surface $C \times C$; see the final section of~\cite{CS21}.

\node{} Although the semi-algorithm of~\cite{CMSV18} generally suffices in practice, for \emph{theoretical decidability} one may use a complete algorithm due to Lombardo~\cite[\S3]{Lombardo2018}, which computes the full endomorphism algebra of an abelian variety. While impractical, it is unconditional. For our purposes, we require only the rank and need not recover Betti or de Rham realizations. See also~\cite[\S7.4]{AL24} for another impractical but complete method. 

\node{} Lombardo's algorithm takes as input equations for an abelian variety $A$ with respect to a projective embedding, together with equations for the subvariety $\Gamma_A \subset A^3$ representing the graph of the group operation $(x, y, x \cdot_A y^{-1})$. It outputs a description of the endomorphism algebra $\End_\Z(A)$, including its rank.

\node{Theorem.}~(\cite{Lombardo2018, CMSV18, CLV21})\label{thm:endo} There exists an unconditional algorithm that takes as input the equations defining a smooth proper curve $C/\Qbar$ and returns a representation of $\End_\Z(J_C)$.
\itnode{Proof.}
Given an irreducible curve $C/\Qbar$, we can compute equations for $J_C$ and its group law by following the approach of~\cite{Chow1954}. For irreducible curves $C_1,\dots,C_k$, we compute projective embeddings of the $J_{C_i}$ and embed their product $\bigoplus J_{C_i}$ into projective space via the Segre embedding. Given equations for a smooth proper curve $C/\Qbar$, this procedure yields equations for $J_C$ and its group operations. We then apply Lombardo's algorithm to compute the rank of $\End_\Z(J_C)$, and the algorithm of~\cite{CMSV18} to compute a representation of the endomorphism ring.
\hfill\qed

\node{Remark.}\label{rem:semisimple}
The algorithm of~\cite{CMSV18} computes a representation of the integral endomorphism algebra $\ce = \End_{\Z}(J_C)$ together with the decomposition of its rational span
\begin{equation}
\label{eq:decompose_ce}
\ce_\Q \coleq \ce \otimes_{\Z}\Q \cong B_1 \times \dots \times B_r
\end{equation}
as a product of simple $\Q$-algebras $B_i$, each isomorphic to a matrix algebra $\Mat_{n_i}(D_i)$ over a division algebra $D_i$ with center a number field; see~\cite[Rem.~7.2.13]{CMSV18}.

Given a finitely generated $\ce$-module $X$ (e.g., a lattice with a $\ce$-action), this decomposition induces a canonical and effective splitting
\begin{equation}
X_{\Q} \coleq X \otimes_{\Z}\Q = X_{1,\Q} \oplus \dots \oplus X_{r,\Q},
\end{equation}
where each $X_{i,\Q}$ is a $B_i$-module. Intersecting with $X$, we obtain sublattices $X_i \coleq X \cap X_{i,\Q}$ of finite index in $X_{i,\Q}$, giving a corresponding decomposition up to isogeny:
\begin{equation}
X \sim X_1 \oplus \dots \oplus X_r.
\end{equation}
This splitting is explicitly computable from the data returned by~\cite{CMSV18}.

\subsection{Induced morphisms between Jacobian motives}\label{sec:induced_morphism_for_jac_motives}

\node{}\label{node:extending_morphism} Suppose $p \colon C \to C'$ is a morphism of smooth proper curves over $\Qbar$, and let $D, E \subset C(\Qbar)$ and $D', E' \subset C'(\Qbar)$ be finite disjoint sets. Then $p$ induces a morphism $p \colon (C \setminus D, E) \to (C' \setminus D', E')$ if and only if $p(E) \subset E'$ and $p^{-1}(D') \subset D$.

\node{} The morphism $p$ induces a pushforward $p_* \colon \Div^0(C \setminus D) \to \Div^0(C' \setminus D')$, which restricts to a morphism of lattices $L_E \to L_{E'}$.  For $g \in \kappa(C)_{D,1}^\times$, we have $\Nm_p(g) \in \kappa(C')_{D',1}^\times$ and $p_*(\divv(g)) = \divv(\Nm_p(g))$. Thus, we have a well defined map
\begin{equation}\label{eq:pushforward_on_albanese}
  p_* \colon \frac{\Div^0(C \setminus D)}{\kappa(C)_{D,1}^\times} \to \frac{\Div^0(C' \setminus D')}{\kappa(C')_{D',1}^\times}
\end{equation}
which respects the markings $L_E \to J_C^D$ and $L_{E'} \to J_{C'}^{D'}$. %

\node{} We denote the induced map on Jacobian motives by $p_* \colon J_{C,E}^D \to J_{C',E'}^{D'}$. At the level of mixed Hodge structures, this corresponds to
\begin{equation}
  p_* \colon \H_1(J_{C,E}^D) = \H_1(C \setminus D, E) \to \H_1(J_{C',E'}^{D'}) = \H_1(C' \setminus D', E').
\end{equation}

\node{}\label{node:dual_of_morphism} The Cartier dual of the map $p_*$ is a morphism $p^! \colon J_{C',D'}^{E'} \to J_{C,D}^E$. On points, this is induced by the pullback map on divisors
\begin{equation}
  p^* \colon \Div^0(C' \setminus E') \to \Div^0(C \setminus E),
\end{equation}
which restricts to a lattice morphism $L_{D'} \to L_D$. Observe that, for $g\in \kappa(C')_{E',1}$ we have $p^*\divv(g) = \divv(g\circ p)$ and $g\circ p \in \kappa(C)_{E,1}$ because $p(E) \subset E'$. Therefore, linear equivalence is preserved. 

\node{Remark.} To prove that the Cartier dual of $p^*$ behaves as we described on points, recall that $J_{C,D}^E$ can be described as a Picard group on $(C \setminus E, D)$ and the pullback map on line bundles is given by divisorial pullback.

\node{}\label{node:cartier_dual_of_pushforward} On homology, the corresponding map is the transfer map
\begin{equation}
  p^! \colon \H_1(C' \setminus E', D') \to \H_1(C \setminus E, D).
\end{equation}

\node{Remark.} To prove this statement, note that $p^!$ is the pullback on cohomology of $(C' \setminus D', E')$, identified with homology of $(C' \setminus E', D')$ via Poincaré--Lefschetz duality. Recall that Cartier duality on motives correspond to the Poincaré--Lefschetz duality on $\H_1$.

\node{Effectivity.}\label{node:correspondence_is_still_effective} Both the action on the points and the action on $\H_1$ for $p^!$ and $p_*$ is effective. The chain of arguments essentially follows that of~\eqref{node:correspondence_is_effective}.

\subsection{Correspondences on Jacobian motives}

\node{} For $i=1,2$ let $C_i$ be smooth proper irreducible curves over $\Qbar$ and $D_i, E_i \subset C_i(\Qbar)$ be finite disjoint sets.

\node{Definition.} An \emph{irreducible correspondence} from $(C_1\setminus{}D_1,E_1)$ to $(C_2\setminus{}D_2,E_2)$ is a correspondace $(C,f_1,f_2)$ from $C_1$ to $C_2$ such that $f_1(f_2^{-1}D_2) \subset D_1$ and $f_2(f_1^{-1}E_1) \subset E_2$. A correspondence is a formal $\Z$-linear combination of irreducible correspondences.

\node{} Given an irreducible correspondence $\calC = (C,f_1,f_2)$ as above, let $D=f_2^{-1}(D_2)$ and $E = f_1^{-1}(E_1)$. Then, recalling~\eqref{node:extending_morphism}, we have morphisms $f_1 \colon (C\setminus{}E,D) \to (C_1 \setminus{}E_1,D_1)$ and $f_2 \colon (C \setminus{}D,E) \to (C_2\setminus{}D_2,E_2)$. Note the flip in the role of punctures and markings for $f_1$.

\node{}\label{node:correspondence_as_a_map_general_case} Consequently, applying the constructions in~\S\ref{sec:induced_morphism_for_jac_motives}, we get maps
\begin{equation}
  J_{C_1,E_1}^{D_1} \overset{f_1^!}{\to} J_{C,E}^D \overset{f_{2,*}}{\to} J_{C_2,E_2}^{D_2},
\end{equation}
whose composition we will denote by 
\begin{equation}
  [\calC_*] = f_{2,*}\circ f_1^! \colon J_{C_1,E_1}^{D_1} \to J_{C_2,E_2}^{D_2}.
\end{equation}
This definition extends additively to formal $\Z$-linear combinations of irreducible correspondences.

\node{} Using~\S\ref{sec:induced_morphism_for_jac_motives}, we see that the map on divisor classes is defined at the level of divisors by $q_*p^* \colon \Div^0(C_1\setminus{}D_1) \to \Div^0(C_2\setminus{}D_2)$ and this restricts to $L_{E_1} \to L_{E_2}$. 

\node{} Using~\S\ref{sec:induced_morphism_for_jac_motives}, we see that the induced map on mixed Hodge structures is $q_* p^! \colon \H_1(C_1\setminus{}D_1) \to \H_1(C_2\setminus{}D_2)$ which is effective by~\eqref{node:correspondence_is_still_effective}.

\node{}\label{node:correspondence_is_still_effective_again} Therefore, given a correspondence $\calC$ from $(C_1\setminus{}D_1,E_1)$ to $(C_2\setminus{}D_2,E_2)$ we can effectively compute the action on divisors representing the map $[\calC_*]\colon J_{C_1,E_1}^{D_1} \to J_{C_2,E_2}^{D_2}$ and the action on $\H_1$.

\node{Lemma.}\label{lem:Dprime} Let $\cc_1,\dots,\cc_k$ be a tuple of correspondences on $C$ and $D \subset C(\Qbar)$ a finite subset. Then we can effectively compute a finite set $D' \subset C(\Qbar)$ containing $D$ such that, for all $i=1,\dots,\rho$, $\cc_i$ induces a correspondence from $(C,D)$ to $(C,D')$, in particular $\cc_{i,*} \colon L_D \to L_{D'}$. Dually, the transpose correspondences $\cc_i^t$ induce correspondences from $C\setminus{}D'$ to $C\setminus{}D$, in particular, $\cc_i^* \colon \Div^0(C\setminus{}D') \to \Div^0(C\setminus{}D)$.
\itnode{Proof.}
Write $\cc_i = \sum_{j} a_{ij} (C_{ij},p_{ij},q_{ij})$ as a sum of irreducible correspondences. Let $D_{ij} \subset C_{ij}$ be the preimage $p_{ij}^{-1}(D)$. Let $D'$ be the union $D\bigcup_{ij}q_{ij}(D_{ij})$. Then, for each $i,j$, $(C_{ij},p_{ij},q_{ij})$ is a correspondence from $(C,D)$ to $(C,D')$. Therefore, $\cc_i$ is a correspondence from $(C,D)$ to $(C,D')$. The dual statement is clear.
\hfill\qed

\node{Remark.} It is tempting to choose a basis $\xi_1,\dots,\xi_d$ for $L_D$ and define $D'$ to be the union of the supports of $\cc_{i,*}\xi_j$. Although this ensures $\cc_{i,*} \colon L_D \to L_{D'}$, each irreducible component of $\cc_{i,*}$ may not induce such a map and then the construction of the action on $\H_1$ requires an extra step.

\newpage
\section{Explicit morphisms and decompositions of push-pull Jacobian motives}

\node{} In this section, we describe the operations on push-pull Jacobian motives that can be carried out explicitly, both at the level of points and mixed Hodge structures.

\node{} We begin by introducing the notions of point-explicit and MHS-explicit morphisms~\eqref{node:explicit_map}. In this context, explicit means that the action on divisor classes and on $\H_1$ can be computed effectively by linear algebra and divisor arithmetic on curves, without invoking higher-dimensional techniques. We also develop techniques to determine more suitable representatives for a given push-pull Jacobian motive in its isogeny class.

\node{} We explain how to compute the kernel of the Abel--Jacobi map from a lattice~\eqref{node:torsion_algo_eff}, and build on this to explicitly split off Baker components from push-pull Jacobian motives~\eqref{lem:explicit_isogeny_second_kind},~\eqref{lem:explicit_isogeny_third_type}.

\node{} We also explain how to compute the period relations of a Baker motive~\eqref{lem:baker_period_relations}.

\subsection{Explicit representations of morphisms}

\node{Representing motives.} In this work, we represent a subclass of motives by push-pull Jacobian motives. Not all motives admit such a description; in fact, not even all abelian varieties are Jacobians. Nevertheless, the subclass we consider is evidently rich enough to represent any $1$-period and, as we will show, is enough to capture all relations.

\node{Representing a push-pull.} To represent a motive $M = [Y \to [X \to A]^\vee]$ as a push-pull of another motive $M' = [Y' \to [X' \to A]^\vee]$ is to represent $M'$ and lattice morphisms $\chi \colon X \to X'$ and $\psi \colon Y \to Y'$ such that $M = \psi^*\chi_*M'$. If $M$ already has another representation, then it is implied that the isomorphism $M \simeq\psi^*\chi_*M'$ (or its inverse) is explicit, at least at the level of $\H_1$.

\node{Different representations, same motive.} Let $D \subset D' \subset C(\Qbar)$ be finite sets. The motive $J_C^D$ is represented as a Jacobian motive, while it is also isomorphic to the pushout $\iota_*J_C^{D'} = J_C^{\iota}$ via the inclusion of lattices $\iota \colon L_D \to L_{D'}$. Although isomorphic, these motives have different representations, affecting how points and Betti or de Rham realizations are described. For example, points of $J_C^D$ are represented by divisors on $C \setminus D$, whereas points of $\iota_*J_C^{D'}$ by divisors on $C \setminus D'$. Such distinctions matter when applying morphisms or moving divisors away from problematic loci.

\node{Representing morphisms.} We will need to single out representations of certain morphisms between push-pull Jacobian motives. We aim for an effective---even, efficient---computation of the action of the map on ``points'' and on $\H_1$. Prime examples are those induced from lattice pushout and pullback maps as well as correspondences. The map induced by a correspondence in the sense of~\eqref{node:correspondence_as_a_map_general_case} is represented by the correspondence itself as in~\eqref{node:represent_correspondence}. 

\node{Representing pullback and pushout maps.} Consider the push-pull Jacobian motive $M = J_{C,\psi}^\chi$ where $\chi \colon X \to L_D$ and $\psi \colon Y \to L_E$. A morphism of lattices $\mu \colon X' \to X$ \emph{represents} the pushout map $J_{C,\psi}^{\chi} \to \mu_*J_{C,\psi}^{\chi} =J_{C,\psi}^{\chi\circ \mu} $. Similarly, a lattice morphism $\nu \colon Y' \to Y$ \emph{represents} the pullback map $\nu^*J_{C,\psi}^\chi = J_{C,\psi\circ\nu}^\chi \to J_{C,\psi}^\chi$. 

\node{Points of a motive.} Let $M = [L \to G]$ be a $1$-motive. By the \emph{points of $M$}, we mean the complex $L \to G(\Qbar)$, or, up to isogeny, the complex $L \otimes \Q \to G(\Qbar) \otimes \Q$. For push-pull Jacobian motives, both the group $G(\Qbar)$ and the map $L \to G(\Qbar)$ can be described explicitly in terms of divisor classes~\eqref{prop:points_of_pushpull}.

\node{Point-explicit morphisms.} A morphism $M_1 \to M_2$, or rather its representation, between two push-pull Jacobian motives is said to be \emph{point-explicit} if the induced map $M_1(\Qbar) \to M_2(\Qbar)$ is effective. In particular, given a divisor representing a divisor class on $M_1$, one can compute a divisor representing its image in $M_2$. We further require that the operations involved reduce to simple linear algebra or to divisor arithmetic on curves and do not rely on higher dimensional algebraic geometry. For instance, correspondences induce point-explicit morphisms~\eqref{sec:induced_morphism_for_jac_motives}. 

\node{Pushout and pullback are point-explicit.} Recall that we represent points of a push-pull Jacobian motive by the class of a tuple $(a, \xi)$, where $a \in \G_m^X$ and $\xi \in \Div^0(C \setminus D)$ for some data $X, C, D$. A pushout morphism induced by a lattice map $\chi \colon X' \to X$ acts on $a$ via the dual map $\widehat{\chi} \colon \G_m^X \to \G_m^{X'}$, while leaving $\xi$ unchanged. That is, the point $[a, \xi]$ is sent to $[\widehat{\chi}(a), \xi]$. The marking is adjusted accordingly, so pushout maps are point-explicit. Pullback maps induced by lattice morphisms are evidently point-explicit, as they only act non-trivially on the lattice part of the points. 

\node{MHS of a motive.} For a motive $M$, the associated mixed Hodge structure is given by $\H_1(M) = (\H_1^\B(M), \H^1_\AdR(M)^\vee, \comp)$. When $M$ is a push-pull Jacobian motive, we have explicit representatives for bases of both underlying vector spaces~\eqref{prop:MHS_of_pushpull}.

\node{MHS-explicit morphisms.} A morphism $M_1 \to M_2$ between two push-pull Jacobian motives is \emph{MHS-explicit} if the induced map between their mixed Hodge structures $\H_1(M_1) \to \H_1(M_2)$ is effective.  Once again, we require that the operations involved reduce to simple linear algebra or to divisor arithmetic on curves and do not rely on higher dimensional algebraic geometry. For instance, correspondences are MHS-explicit~\eqref{node:correspondence_is_still_effective_again}.

\node{} Note that, given basis elements for the Betti and de Rham realizations of $M_1$, one can compute under an MHS-explicit map their images in the corresponding realizations of $M_2$. This means we have explicit linear maps: matrices in chosen bases that describe the induced morphisms on each realization. 

\node{Pullbacks and pushouts are MHS-explicit.} This is clear from the representation of the MHS given in~\eqref{prop:MHS_of_pushpull}. Pullback and pushout maps will act as the identity on the central component $\H_1(C\setminus{}D,E)$ and act naturally via the lattice maps on the outer components.

\node{Explicit map.}\label{node:explicit_map} A morphism between two push-pull Jacobian motives is said to be \emph{explicit} if it is both point-explicit and MHS-explicit. That is, we can compute its action on divisor classes representing points, and describe its induced morphism on mixed Hodge structures via explicit matrices. We established above that pullback and pushout morphisms as well as correspondences are explicit.

\node{Period effective maps.} An MHS-explicit map $M_1 \to M_2$ will be called \emph{period effective} if $\cR(M_1)$ can be determined from $\cR(M_2)$ by applying prescribed linear algebra operations. Trivially, MHS-explicit isomorphisms are period effective. We also proved that pushouts and pullbacks are period effective~\eqref{prop:pushpull_relations}. 

\subsection{Isogenies, direct sum decompositions, and making lattice maps injective}

\node{} Our primary concern is with period relations and, therefore, with the MHS $\H_1$. However, we cannot compute effectively with the comparison isomorphism and we resort to numerical approximations. The motives are used to remedy this defect: their points allow us to recover the exactness lost in working with the approximate comparison isomorphism. In this sense, motives serve as an effective model for the underlying $\H_1$. Naturally, we wish to allow ourselves to work up to isogeny, as this leaves $\H_1$ unchanged.

\node{Effectively interchanging isogenous motives.} Suppose $M_1 \to M_2$ is an MHS-explicit isogeny, so that the induced map $\H_1(M_1) \to \H_1(M_2)$ is an effective isomorphism. Then we obtain an identification of the space of period relations: $\cR(M_1) \isoto \cR(M_2)$. Thus, for our purposes, we may pass between $M_1$ and $M_2$ interchangeably.

\node{Direct sum decompositions.}\label{node:explicit_split} An \emph{explicit direct sum decomposition} of a motive $M_3$ consists of explicit morphisms $M_1 \to M_3$ and $M_2 \to M_3$ (or $M_3 \to M_1$ and $M_3 \to M_2$) such that the induced map $M_1 \oplus M_2 \to M_3$ (or $M_3 \to M_1 \oplus M_2$) is an isomorphism. In any case, we obtain an explicit identification $\H_1(M_3) = \H_1(M_1) \oplus \H_1(M_2)$. Furthermore, we get a pointwise isomorphisms $M_1(\Qbar) \oplus M_2(\Qbar) \isoto M_3(\Qbar)$, or in the other direction. We call this giving an explicit direct sum decomposition $M_3 \simeq M_1 \oplus M_2$.

\node{Isogeny direct sum decomposition.}\label{node:explicit_isogeny_split} An \emph{explicit isogeny direct sum decomposition} of a motive $M_3$ into $M_1$ and $M_2$ consists of explicit isogenies $M_1 \sim M_1'$, $M_2 \sim M_2'$, and explicit morphisms $M_3 \to M_1'$ and $M_3 \to M_2'$ such that $M_3 \to M_1' \oplus M_2'$ is an isogeny. This gives an explicit identification $\H_1(M_3) = \H_1(M_1) \oplus \H_1(M_2)$ and an explicit, pointwise, map $M_3(\Qbar) \otimes \Q \isoto (M_1(\Qbar) \oplus M_2(\Qbar))\otimes \Q$. As in~\eqref{node:explicit_split}, we also allow for maps in the other direction. We will refer to this as an explicit isogeny direct sum decomposition and denote it by $M_3 \sim M_1 \oplus M_2$. 

\node{Disjoint union decomposition.} Suppose $(C \setminus D, E)$ is the disjoint union of $(C_i \setminus D_i, E_i)$'s. Then the associated Jacobian motive $J_{C,E}^D$ naturally admits an explicit direct sum decomposition into the $J_{C_i,E_i}^{D_i}$'s. This decomposition is induced by functoriality and is explicit on both points and mixed Hodge structures.

\node{Making lattice maps injective.}\label{node:make_lattice_injective}
Let $J_{C,\chi}$ be a pullback Jacobian motive with $\chi \colon X \to L_D$. Then we may effectively construct a disjoint union $C'$ of $C$ with finitely many copies of $\p$, a finite set $D' \subset C'(\Qbar)$, and an \emph{injective} lattice map $\chi' \colon X \to L_{D'}$, such that there exists an explicit isomorphism $J_{C,\chi} \cong J_{C',\chi'}$.  This gives a model in which the lattice map is injective---a seemingly minor step that yields significant expository simplification.

\itnode{Proof.}
Let $\chi \colon X \to L_D$ and let $K = \ker \chi$. Choose a basis $e_1,\dots,e_k$ for $K$, and a complement $X' \subset X$ such that $X = K \oplus X'$. Let $C'$ be the disjoint union of $C$ and $k$ copies of $\p$, and let $D'$ be the union of $D$ with $\{0,\infty\}$ on each added copy.

Define $\chi' \colon X \to L_{D'}$ by sending $e_i \mapsto \infty - 0$ on the $i$-th $\p$, and restricting to $\chi$ on $X'$. Then $\chi'$ is injective by construction.

Observe that $J_{C,\chi} \cong [K \to 0] \oplus [X' \to J_C]$ is an explicit decomposition because the kernel of the lattice map splits from the representation of points and of $\H_1$ by the pullback construction. On the other hand, $J_{C',\chi'} \cong [\Z \to 0]^{\oplus k} \oplus [X' \to J_C]$ by the disjoint union decomposition of $(C',D')$. The two motives are identified via identification of summands.
\hfill\qed

\node{Remark.}\label{rem:injective_lattice} The construction in~\eqref{node:make_lattice_injective} is effective and leaves $J_C$ unchanged. We will invoke this operation by saying that we may assume $\chi$ is injective, without renaming $C$, $D$, or $\chi$. 

\subsection{Computing the kernel of an Abel--Jacobi map}\label{sec:computing_AJ_kernel}

\node{} Consider $J_{C,D} = [L_D \overset{\AJ}{\to} J_C]$.  Let $K_D \subset L_D$ be the kernel of the Abel--Jacobi map $\AJ \colon L_D \to J(\Qbar)$, these are the divisors supported on $D$ which are linearly equivalent to $0$. We will give an algorithm here to compute $K_D$. Let $T_D \subset L_D$ be the kernel of $L_D \to J(\Qbar) \otimes_\Z \Q$, these are the divisors supported on $D$ whose linear equivalence classes are torsion.

\node{Saturating a lattice.} Once we know $K_D \subset L_D$, we can saturate it in $L_D$ to compute $T_D$. If we choose bases and represent the inclusion $K_D \hookrightarrow L_D$ by an integer matrix, the saturation of the image can be computed from the Smith normal form of the matrix. 

\node{Preparation.} Without loss of generality, we may assume $C$ is irreducible. Fix a basis $\xi_1,\dots,\xi_n$ for $L_D$, these are divisors of degree $0$ on $C$ supported on $D$. Let $\kk \subset \Qbar$ be a number field over which $C$ and each point of $D$ can be defined. In particular, the divisors $\xi_1,\dots,\xi_n$ are defined over $\kk$.

\node{The Néron--Tate height pairing.}\label{node:height_pairing} The algorithm of~\cite{BHM19} allows us to compute the Néron--Tate height pairing $\langle [\xi_i], [\xi_j] \rangle$ of the divisor classes to arbitrary precision. By performing the relevant period computations with certified error bounds, we may compute complex balls containing each N\'eron--Tate height with desired precision. The N\'eron--Tate height pairing is positive definite on $J_C(\kk) \otimes \R$. The matrix $\left(\langle [\xi_i], [\xi_j] \rangle\right)_{i,j=1}^n$ gives the pullback of the height pairing to $L_D \otimes \R$. The rank of this pullback matrix is the rank of $L_D \otimes \Q \to J_C(\Q)$ despite the real coefficients, since the map is defined over $\Q$.

\node{Testing membership in the kernel.}\label{node:test_kernel}
Given an element $\xi \in L_D$, $\AJ(\xi) = 0$ in $J_C$ if and only if $\xi$ is linearly equivalent to zero. This can be tested effectively by computing the Riemann--Roch space of $\xi$. In particular, given a sublattice $K' \subset L_D$, performing this test on a basis for $K'$ we can decide whether $K' \subset \ker \AJ$.

\node{An efficient guess for the kernel.}\label{node:guess_AJ_kernel}
We can use lattice reduction~\cite{LLL} to efficiently guess a sublattice $K' \subset L_D$ likely to lie in $K_D$. Using the filtrations and the comparison isomorphism on $\H_1(C,D)$ we have an approximation of the map
\begin{equation}
  \H_1^\B(C,D;\Z) \to (F^1\H^1_\AdR(C,D))^\vee = (F^1\H^1_\AdR(C))^\vee
\end{equation}
Use LLL to guess the integral kernel of the map. Project the putative kernel to $L_D = \H_1^\B(C,D)/\H_1^\B(C)$ to formulate the guess $K'$. This guess $K'$ can be refined by increasing the working precision on the period matrix.

\node{Computing the kernel.}\label{node:torsion_algo_eff}
Let $K= 0$. We start by formulating a guess~\eqref{node:guess_AJ_kernel} of a kernel $K' \subset L_D$ and check~\eqref{node:test_kernel} if a basis for $K'$ consists of divisors linearly equivalent to $0$. Update $K$ by adding the basis elements of $K'$ that were linearly equivalent to $0$. Increase precision of the period matrix used for the guess and repeat. This gives an non-decreasing sequence of lattices $K$ in $K_D$. 

In parallel, we will compute a non-decreasing sequence of lower bounds on the rank of $\AJ \otimes \R$ on $L_D \otimes \R$. To find the lower bounds, we simply check which minors of the height pairing~\eqref{node:height_pairing} on $L \otimes \R$ do not contain zero using higher and higher precision. 

Whenever our lower bound matches the corank of $K$, we compute the saturation $T$ of $K$, find finitely many elements in $T$ dominating $T/K$ and check which of these finitely elements are linearly equivalent to $0$. Then $K_D$ is the span of $K$ with these new elements.

The algorithm terminates because LLL will eventually find a lattice of full rank in $K_D$ and our approximation of the minor of the correct corank will eventually stop containing $0$.

\subsection{Splitting Baker components, part I}

\node{} We will now explicitly split the kernel $K_D$ of the Abel--Jacobi map from the Jacobian motive $J_{C,D} = [L_D \to J_C]$. This means giving a pullback Jacobian $J_{C,\aleph_D}$, $\aleph_D \colon N_D \to L_D$, and an explicit isogeny direct sum decomposition $J_{C,D} \sim [K_D \to 0] \oplus J_{C,\aleph_D}$, recall~\eqref{node:explicit_isogeny_split}.

\node{} We first need an explicit map from $[K_D \to 0]$ into $J_{C,D}$, meaning it must be both point-explicit and MHS-explicit~\eqref{node:explicit_map}.

\node{} At the level of points, the only non-trivial part of the inclusion of the motive $[K_D \to 0]$ into $J_{C,D}$ is the inclusion of lattices $K_D \subset L_D$. Having computed $K_D$ explicitly in $L_D$ in~\eqref{node:torsion_algo_eff}, there is nothing to prove here.

\node{Lemma.} The injection $[K_D \to 0] \to J_{C,D} = [L_D \to J_C]$ is MHS-explicit.
\itnode{Proof.}
Let $e_1,\dots,e_k$ be a basis for the kernel $K_D \subset L_D$. For each $e_i \in K_D$, compute a rational function $f_i \in \kappa(C)^\times$ such that $e_i = \divv(f_i)$.

Recall that the Jacobian motive of $(\mathbb{P}^1, \{0,\infty\})$ satisfies $J_{\mathbb{P}^1,\{0,\infty\}} \cong [\mathbb{Z}\to 0]$. Let $\Gamma_{f_i}$ denote the graph of the rational function $f_i \colon C\to \p$. Consider the transpose correspondence $\Gamma_{f_i}^t$ from $(\p,\{0,\infty\})$ to $(C,D)$. It acts~\eqref{node:dual_of_morphism} by pulling back the unique divisor class $[\infty-0]$ to $\divv(f_i)$, thus mapping the generator $1$ from $J_{\p,\{0,\infty\}}$ to $e_i$ in $J_{C,D}$.

Identify the motive $[K_D \to 0]$ with the Jacobian motive of a disjoint union of $k$ copies of $(\mathbb{P}^1,\{0,\infty\})$. Define a correspondence $\Gamma \subset \left(\coprod_{i=1}^k \mathbb{P}^1\right) \times C$ as the disjoint union of graphs $\Gamma_{f_i} \subset \mathbb{P}^1 \times C$. The transpose correspondence $\Gamma^t$ defines a morphism from the Jacobian motive $J_{\mathbb{P}^1,\{0,\infty\}}^{\oplus k}$ to the Jacobian motive $J_{C,D}$, representing the injection $[K_D \to 0] \toi J_{C,D}$.

This map is MHS-explicit since morphisms represented by correspondances have this property~\eqref{node:correspondence_is_still_effective_again}.
\hfill\qed

\node{Remark (Explicit interpretation of the vanishing condition (14.3)).}\label{node:reinterpret_143} In the monograph~\cite[p.137]{HW}, a certain vanishing condition ``(14.3)'' is introduced for rational differentials of the second type. Explicitly, given a rational differential $\omega$ on $C$ of the second type, and a submotive $\iota \colon [\Z \to 0] \toi J_{C,D}, \quad 1 \mapsto \divv(f), \quad f\in \kappa(C)^\times$, condition (14.3) asks for the vanishing of $\iota^*([\omega])$, where $[\omega] \in \H^1_\AdR(C,D)$ is the cohomology class of $\omega$. The statement of Lemma~\eqref{lem:explicit_isogeny_second_kind} implies that this condition can be explicitly checked by arithmetic on curves. We will spell this out here, as it may be of independent interest.

Following the proof of the lemma above, we need the Cartier dual~\eqref{node:dual_of_morphism} of the pushforward map given by $f \colon C\setminus{}D \to \p\setminus{}\{0,\infty\}$. On mixed Hodge structures the induced map is $f^! \colon \H_1(\mathbb{P}^1,\{0,\infty\}) \to \H_1(C,D)$, see~\eqref{node:cartier_dual_of_pushforward}. 

On algebraic de Rham cohomology, this corresponds to the transfer map,
\begin{equation}
f_! \colon \H^1_\AdR(C,D) \to \H^1_\AdR(\mathbb{P}^1,\{0,\infty\}),
\end{equation}
which evaluates to the $f$-trace on differentials~\eqref{node:coh_transfer} mapping $[\omega]$ to $[\tr_f \omega]$. Note that the trace of $\omega$ has no residual poles and must be exact on $\p$.

We have an identification $\H^1_\AdR(\mathbb{P}^1,\{0,\infty\}) \simeq \Qbar^{\{0,\infty\}}/\Qbar \simeq \Qbar$, which sends $\dd h \mapsto h(0)-h(\infty)$, see~\eqref{node:df_to_fE}. Thus, the condition $\iota^*([\omega])=0$ is explicitly checked by computing the rational function $h$ such that $\tr_f\omega = \dd h$ and evaluating the difference $h(0)-h(\infty)$.

\node{Representing the quotient.}\label{node:rep_quot} The quotient $J_{C,D}/[K_D \to 0]$ admits an explicit representation as a pullback motive, which we describe here. This representation is not unique, but any two representations differ only by an isogeny. Let $K_D \subset T_D$ be as computed previously, and define $n \in \N$ to be the least common multiple of the orders of elements in the quotient $T_D/K_D$. Define $L_{nD} \coleq n \cdot L_D \subset L_D$ and let $J_{C,nD}$ denote the pullback associated to $L_{nD} \toi L_D$. By construction, $L_{nD} \cap T_D = L_{nD} \cap K_D = K_{nD}$, the kernel of the Abel–Jacobi map on $L_{nD}$. Define $N_D \coleq L_{nD}/K_{nD}$, which naturally maps to $J_C$. The resulting motive $M\coleq[N_D \to J_C] = J_{C,nD}/[K_{nD} \to 0]$ is isogenous to the original quotient $J_{C,D}/[K_D \to 0]$. First note that for any splitting $\aleph_D \colon N_D \to L_{nD} \subset L_D$, the map $N_D \to J_C$ will factor through $\aleph_D$, hence $M \simeq J_{C,\aleph_D}$. Moreover, the natural inclusion $J_{C,\aleph_D} \toi J_{C,nD}$ followed by the quotient map induces an explicit isomorphism
\begin{equation}
  \H_1(J_{C,\aleph_D}) \isoto \H_1(J_{C,nD})/\H_1([K_{nD} \to 0]) = \H_1(M),
\end{equation}
as well as the identity map on the points $J_{C,\aleph_D}(\Qbar) \isoto M(\Qbar)$. To comply with our choice of representation, we will write $J_{C,\aleph_D}$ for $M$, however, we could have slightly extended our representation to allow for quotients of lattices and cover $M$ since the choice of $\aleph_D$ plays no role here.

\node{Lemma.}\label{lem:explicit_isogeny_second_kind} We have an explicit isogeny direct sum decomposition $J_{C,D} \sim [K_D \to 0] \oplus J_{C,\aleph_{D}}$. In fact, we can describe the decomposition maps in either direction
\itnode{Proof.}
The explicit maps $J_{C,\aleph_D} \to J_{C,D}$ and $[K_D \to 0] \to J_{C,D}$ give an explicit isogeny direct sum decomposition of $J_{C,D}$.

We can also realize the isogeny in the other direction. First, replace $J_{C,D}$ with the explicitly isogenous representative $J_{C,nD}$. The lattice $L_{nD}$ splits via $\aleph_D$ into $K_{nD} \oplus N_D$. Projecting onto $[K_{nD}\to 0]$ and onto the quotient $J_{C,nD}/[K_{nD}\to 0] = [N_D \to J_C] \simeq J_{C,\aleph_D}$ recovers the decomposition in the other direction.
\hfill\qed

\subsection{Splitting Baker components, part II}

\node{} We continue with the notation of the previous part and consider $J_C^D$. This time we will give an explicit isogeny direct sum decomposition $J_C^D \sim \G_m^{K_D} \oplus J_{C,\aleph_D}$.

\node{Lemma.}\label{lem:explicit_surjection}
The canonical surjection from $J_C^D$ to $\G_m^{K_D}$ is explicit.
\itnode{Proof.}
The proof is dual to that of~\eqref{lem:explicit_isogeny_second_kind}, using $\p \setminus \{0,\infty\}$ in place of $(\p,\{0,\infty\})$. We repeat a part of the argument here to clarify the action on points.

Let $K_D = \langle \divv(f_1), \dots, \divv(f_k) \rangle$ for some $f_i \in \kappa(C)^\times$, and let $\Gamma_{f_i}$ denote the graph of $f_i \colon C \setminus D \to \p \setminus\{0,\infty\}$. Since $J_{\p \setminus\{0,\infty\}} = \G_m$, we identify $\G_m^{K_D} \simeq (J_{\p}^{\{0,\infty\}})^{\oplus k}$ as the Jacobian of $k$ disjoint copies of $\p \setminus\{0,\infty\}$.

Let $\Gamma = \coprod \Gamma_{f_i}$ be the disjoint union of these graphs. Then $\Gamma$ defines a correspondence from $C \setminus D$ to $\coprod (\p \setminus\{0,\infty\})$ inducing the desired map $J_C^D \to \G_m^{K_D}$. The induced map is explicit by virtue of being induced by a correspondance.
\hfill\qed

\node{Action on points made explicit.} \label{node:point_proj_to_torus}
The Abel--Jacobi map $\AJ \colon \Div^0(\p \setminus \{0,\infty\}) \to J_{\p}^{ \{0,\infty\}} \simeq \G_m $ is the one mapping $[1:a]-[1:b]$ to $a/b$. Equivalently, any divisor $\divv h$ is mapped to $h(\infty)/h(0) \in \G_m$. Given a divisor $\xi$ representing a point of $J_C^D$ the point $\AJ(f_{i,*}(\xi)) \in \G_m$ is the $i$-th coordinate of the map $J_C^D \tos \G_m^{K_D}$ for $K_D = \langle \divv(f_1), \dots, \divv(f_k) \rangle$.

\node{Lemma.}\label{lem:explicit_isogeny_third_type} We have an explicit isogeny direct sum decomposition $J_C^D \sim \G_m^{K_D} \oplus J_C^{\aleph_D}$. We can describe the decomposition maps in either direction.
\itnode{Proof.}
The map $J_C^D \tos \G_m^{K_D}$ is made explicit in~\eqref{lem:explicit_surjection}. The map $J_C^D \tos J_C^{\aleph_D}$ is explicit by pushout construction. This gives the explicit isogeny direct sum decomposition in one direction. To get the isogeny in the other direction repeat the argument in~\eqref{lem:explicit_isogeny_second_kind} by taking duals.
\hfill\qed

\node{}\label{node:including_aleph} We explicate the inclusion $J_C^{\aleph_D} \toi J_C^D$ for a representative $[\xi]_{\aleph_D}$, where $\xi \in \Div^0(C\setminus{}D)$. Use the explicit map $\pi \colon J_C^D \to \G_m^{K_D}$ described in~\eqref{node:point_proj_to_torus} to compute $\pi([\xi]_D) \in \G_m^{K_D}$. Then, use the splitting $\iota\colon \G_m^{K_D} \toi \G_m^{L_D} \toi J_C^D$ induced from $L_D \sim K_D \oplus N_D$ and compute $[\xi']_D \coleq [\xi]_D - \iota(\pi([\xi]_D))$ which necessarily lies in $\ker \pi \sim J_C^{\aleph_D}$. 

\subsection{Splitting Baker components, part III}

\node{} We continue with the notation of the previous parts. We will now split the Baker components from $J_{C,\chi}$ and $J_C^{\chi}$ for a lattice map $\chi \colon X \to L_D$.

\node{}\label{node:K_chi} Let $K_\chi \subset X$ be the pullback of $K_D \subset L_D$, i.e., the kernel of $X \to J_C$. Since $K_D \subset L_D$ was computed, $K_\chi$ is computable.

\node{}\label{node:ignore_torsion_in_chi} Scaling $\chi$ by an integer corresponds to passing to an isogenous representative of $J_{C,\chi}$. Then, taking $n>0$ as in~\eqref{node:rep_quot}, work with $n \chi \colon X \to L_{nD}$ so the pullback of torsion divisors classes and $0$ divisor classes agree and give $K_\chi$. Thus $N_\chi \coleq X/K_\chi \to J_C$ is well defined. Choose a lift ${\aleph_\chi} \colon N_\chi \to L_{nD} \subset L_D$ to represent $[N_\chi \to J_C]$ as the pullback $J_{C,{\aleph_\chi}}$. By construction $ N_\chi \to J_C(\Qbar)\otimes \Q$ is injective, thus the motives $J_{C,\aleph_\chi}$ and $J_C^{\aleph_\chi}$ are reduced~\eqref{lem:reduced_is_lattice_inj}.

\node{} Recall~\eqref{node:make_lattice_injective} that by a small modification to $C,D,\chi$ we may assume that $\chi \colon X \to L_D$ is injective. We make this assumption now.

\node{Proposition.}\label{prop:explicit_isogeny_with_chi} We have an explicit isogeny direct sum decompositions $J_{C,\chi} \sim [K_\chi \to 0] \oplus J_{C,{\aleph_\chi}}$ and $J_C^\chi \sim \G_m^{K_\chi} \oplus J_C^{\aleph_\chi}$. We can describe the decomposition maps in either direction. Observe that $J_{C,{\aleph_\chi}}$ is reduced. 
\itnode{Proof.}
The maps $[K_\chi \to 0] \to J_{C,\chi}$ and $J_{C,{\aleph_\chi}} \to J_{C,D}$ are trivially point-explicit. To determine the decomposition at the level of $\H_1$ observe that we have an inclusion, having made $\chi$ injective, $\H_1(J_{C,\chi}) \toi \H_1(J_{C,D})$. We may now use the decomposition on the right hand side from~\eqref{lem:explicit_isogeny_second_kind}. For the maps in the reverse direction repeat the arguments there.

For the dual decomposition $J_C^\chi \sim \G_m^{K_\chi} \oplus J_C^{\aleph_\chi}$ we make use of~\eqref{lem:explicit_isogeny_third_type}. For instance, injectivity of $\chi$ implies $J_C^D \tos J_C^\chi$ and $K_\chi \toi K_D$ implies $\G_m^{K_D} \tos \G_m^{K_\chi}$. To describe the map $J_C^\chi \to \G_m^{K_\chi}$ lift the divisor to $J_C^D$, tautological at the level of representatives, map to $\G_m^{K_D}$ and then down to $\G_m^{K_\chi}$. On the map of $\H_1$, use the splitting $\H_1(\G_m^{K_D}) \toi \H_1(J_C^D)$ followed by the surjection $\H_1(J_C^D) \tos \H_1(J_C^\chi)$. The other map $J_C^\chi \tos J_C^{\aleph_\chi}$ is induced by pushforward. For the maps in the reverse direction we repeat the argument in~\eqref{lem:explicit_isogeny_second_kind} and use~\eqref{node:including_aleph}. 
\hfill\qed

\subsection{Splitting Baker components, part IV}

\node{} Given a push-pull motive $M = J_{C,\psi}^{\chi}$ with $\chi \colon X \to L_D$ and $\psi \colon Y \to L_E$, we will give an algorithm that produces a push-pull Jacobian motive $M'$, an explicit isogeny direct sum decompoition of $M'$ into a Baker component and a reduced push-pull Jacobian, and a realization of $M$ as a push-pull of $M'$. 

\node{Remark.} Recall that period relations of $M$ can be computed from that of $M'$~\eqref{prop:pushpull_relations} and the periods of $M'$ can be computed once we can split it into a direct sum of a Baker and a reduced motive~\eqref{cor:rels_of_Baker_and_sat}.

\node{} We make use of the notation in the previous parts. Anything defined for $\chi$, e.g., $\aleph_\chi$, will also be defined for $\psi$. Modify the representative as in~\eqref{node:make_lattice_injective} to assume both $\chi$ and $\psi$ are injective. Also, pass to an isogenous representative, as in~\eqref{node:make_lattice_injective}, so as to assume that neither $\chi$ or $\psi$ hit non-trivial torsion divisors.

\node{}\label{node:chi_torsion}  Consider the map $L_E \to J_C^{\chi}$. Compute $K_E \subset L_E$ and observe that $K_E \to \G_m^{X}$ can be made explicit: write $K_E = \langle \divv(f_1),\dots,\divv(f_k) \rangle$ and recall~\eqref{node:lattice_pullback_of_funcs} that $K_E \to \G_m^X$ is described by $\divv(f_i) \mapsto (f_i)|_D\circ\chi$. Compute the kernel $K^\chi_E \coleq \ker (K_E \to \G_m^X)$.  By scaling $\psi$ again, we assume $K^\chi_\psi \coleq \psi^{-1}K^\chi_E$ is saturated in $Y$. Let $N_\psi^\chi \coleq Y/K_\psi^\chi$ and choose a splitting $\aleph_\psi^\chi \colon N_\psi^\chi \to Y$. If $Y=L_E$ and $\psi$ was the identity (prior to scaling) then we write $N_E^\chi$ and $\aleph_E^\psi$ instead.

\node{Lemma.}\label{lem:split_kchipsi} The map $[K^\chi_\psi \to 0] \to J_{C,\psi}^\chi$ can be made explicit. In particular, we have an explicit direct sum decomposition $J_{C,\psi}^\chi \simeq [K^\chi_\psi \to 0] \oplus J_{C,\aleph_\psi^\chi}^\chi$.
\itnode{Proof.}
The action on the points is a triviality. The difficulty is in computing the image under $\H_1$. We can choose a basis of $K^\chi_\psi$ consisting of divisors of functions $f \in \kappa(C)_D$. Therefore, it will be sufficient to fix one such $f$ with $[\divv(f)]_\chi = 0$ with support in $Y \toi L_E$ and to determine the corresponding image of $[\Z \to 0] \to J_{C,\psi}^\chi$ under $\H_1$.

Since $\divv(f)$ is supported on $E$, we can view it as a map $f \colon (C\setminus{}E,D) \to (\p\setminus{}\{0,\infty\},f(D))$. The pullback induces an explicit map $f^* \colon J_{\p,\{0,\infty\}}^{f(D)} \to J_{C,E}^D$. Here $J_{\p,\{0,\infty\}}^{f(D)} = [\Z\langle \infty - 0 \rangle \to J_{\p}^{f(D)}]$ and $f^*$ maps $\infty -0$ to $\divv(f)$.

Note that $\chi \to L_D \overset{f_*}{\to} L_{f(D)}$ gives a pushout $J_{\p,\{0,\infty\}}^{f_*\circ\chi}$ of the Baker motive $J_{\p,\{0,\infty\}}^{f(D)}$ which maps to $J_{C,E}^\chi$. Since $[\divv(f)]_\chi = 0$, the pushout has a split $[\Z \to 0]$ submotive. Making this splitting explicit by way of correspondences is possible but cumbersome. In practice, it is much easier to put the period matrix into block form, by change of bases in the de Rham and Betti realizations, using the fact that we can work effectively with the Baker periods~\eqref{lem:baker_period_entries}. Once $\H_1([\Z \to 0]) \toi \H_1(J_{\p,\{0,\infty\}}^{f_*\circ\chi})$ is computed, map it to $\H_1(J_{C,E}^\chi)$ using $f$, and then use the pullback via $\psi$ to determine the image in $\H_1(J_{C,\psi}^\chi) \subset \H_1(J_{\p,\{0,\infty\}}^{f_*\circ\chi})$ (here we use $\divv(f) \in Y \toi L_E$).

The map $J_{C,\aleph_\psi^\chi}^\chi \toi J_{C,\psi}^\chi$ is explicit and completes the direct sum decomposition together with $[K^\chi_\psi \to 0] \toi J_{C,\psi}^\chi$.
\hfill\qed

\node{Lemma.}\label{lem:pullback_remove_Baker} We can effectively realize $J_{C,\psi}^\chi$, up to isogeny, as a pullback of a motive $M' = [Y \to \G_m^{K_\chi}] \oplus J_{C,\psi}^{\aleph_\chi}$ given as a direct sum. Note that the submotive $J_C^{\aleph_\chi}$ is reduced.
\itnode{Proof.}
Explicitly split $J_C^\chi \sim \G_m^{K_\chi} \oplus J_C^{\aleph_\chi}$ as in~\eqref{prop:explicit_isogeny_with_chi}. Since this isogeny is point explicit, we may map our marking $Y \to J_C^\chi$ to the two components, at least after scaling $\psi$. The pullback via the diagonal $Y \to Y\oplus Y$ recovers $J_{C,\psi}^\chi$.
\hfill\qed

\node{Lemma.}\label{lem:pushout_remove_Baker} We can effectively realize $J_{C,\psi}^\chi$, up to isogeny, as a pushout of a motive $M' = [K_\psi \to \G_m^X] \oplus J_{C,\aleph_\psi}^\chi$. Note that the quotient motive $J_{C,\aleph_\psi}$ is reduced.
\itnode{Proof.}
The construction of $M'$ is, by now, clear. Pushout $M'$ via the diagonal $X \to X\oplus X$ to recover $J_{C,\psi}^\chi$ up to isogeny.
\hfill\qed

\node{Proposition.}\label{prop:split_baker} Given $J_{C,\psi}^\chi$, we can effectively construct a Baker motive $M_B$ and a \emph{reduced} push-pull Jacobian motive $M^{\red}$ together with lattice maps realizing $J_{C,\psi}^\chi$ as a push-pull of $M_B \oplus M^{\red}$, up to isogeny. 
\itnode{Proof.}
We note that $\aleph_\chi$ and $\aleph_\psi$ give reduced motives~\eqref{node:ignore_torsion_in_chi}. Now apply the two lemmas~\eqref{lem:pullback_remove_Baker} and~\eqref{lem:pushout_remove_Baker}, one after the other, removing the Baker components at each step to end up with $M^{\red} = J_{C,\aleph_\psi}^{\aleph_\chi}$.
\hfill\qed

\subsection{Period relations of a Baker motive}

\node{} Let $M_B$ be a Baker motive. Just like any other motive in our computational scheme, it is represented as a push-pull Jacobian motive. We will give an algorithm to compute the period relations $\cR(M_B) \subset \H_1^\B(M_B)\otimes \H^1_\AdR(M_B)$.

\node{} Let $C = \p$ and $D,E \subset C(\Qbar)$ be disjoint finite sets. Denote the points of $E$ by $e_i$'s and $D$ by $d_i$'s. To simplify notation, we will assume $d_0 \in D$ is $\infty$ and identify the rest of the points in $E, D$ with algebraic numbers $\Qbar$. For each $i,j$ fix a branch of the logarithm $\log(e_j - d_i)$. We will say that a period of $J_{\p,E}^D$ is in \emph{standard form} if it is expressed as an explicit $\Qbar$-linear combination of $1$, $\tpi$, and $\log(e_i-d_j)$'s.

\node{Lemma.}\label{lem:baker_period_entries} 
Given a representation of $\H_1(\p\setminus{}D,E)$, we can effectively express each entry in its period matrix in standard form.

\itnode{Proof.} 
We are given representative elements for bases for $\H^1_\AdR(\p\setminus{}D,E)$ and $\H_1^\B(\p\setminus{}D,E)$. We fix one of each: a representative form $(a,\omega)$ and $1$-chain $\gamma$. We will write the period $\int_\gamma (a,\omega) = a(\partial \gamma) + \int_\gamma \omega$ in standard form. The first term is an explicit algebraic multiple of $1$ so we focus on the second term. 

Write $\omega = P(x)/Q(x)\, \dd x$, and decompose this into a sum of elementary fractions:
\begin{equation}
\omega = \sum_i \frac{P_i(x)}{(x - d_i)^{n_i}}\, \dd x,
\end{equation}
with $n_i \geq 0$, $\deg P_i < n_i$, and $P_i$ relatively prime to $(x - d_i)$. If $n_i \neq 1$, then each summand admits an elementary antiderivative, and the corresponding integral can be computed by evaluating this antiderivative at the algebraic endpoints of $\gamma$ giving an explicit algebraic number.

When $n_i = -1$, $P_i$ is a constant, the residue of $\omega$ at $c_i$, call it $r_i \in \Qbar$. The local antiderivatives are branches of $r_i \log(x - d_i)$. The integral of this term over $\gamma$ will differ from the evaluation of the logarithm on $\partial \gamma$ by an \emph{integer} multiple of $r_i\tpi$. This integer multiple can be explicitly and rigorously determined, for instance by numerical approximation of the integral and of the logarithms with rigorous error bounds. The evaluation of $r_i \log(x - d_i)$ on $\partial \gamma$ is an explicit $\Qbar$-linear combination of $\log(e_j-d_i)$'s.
\hfill\qed

\node{Lemma.}\label{lem:baker_period_relations}
Let $M_B$ be a Baker motive. Then the space $\cR(M_B)$ of $\Qbar$-linear relations between the entries of its period matrix can be computed effectively.

\itnode{Proof.} 
Since a Baker motive has trivial abelian part, we may write it as a push-pull of $J_{C,E}^D$, with $C = \coprod_{i=1}^k \Pp^1$ and $D,E \subset C(\Qbar)$ disjoint finite subsets. By~\eqref{prop:pushpull_relations}, the period relations of $M_B$ are determined by those of $J_{C,E}^D$, and thus we may assume $M_B = J_{C,E}^D$. The Jacobian motive $J_{C,E}^D$ is a direct sum of irreducible Jacobian motives with the underlying curve $C = \p$. We may, therefore, apply Lemma~\eqref{lem:baker_period_entries} to write the entries of the period matrix of $J_{C,E}^D$ in standard form.

Let $N = \dim \H_1^\B = \dim \H^1_\AdR$ and identify the space $\H_1^\B \otimes \H^1_\AdR$ with $\Qbar^{N^2}$ via the representing bases. The entries of the period matrix represents the period pairing in our coordinates $\Qbar^{N^2} \to \C$. The kernel $\cR(M_B)$ is the space of relations between the entries of the period matrix. 

Using Baker's theorem~\cite{Baker1966}, see also~\cite[Theorem 1.05]{HW}, any $\Qbar$-linear relation among logarithms of algebraic numbers lies in the $\Qbar$-span of $\Q$-linear relations among the same logarithms, which are in turn spanned by the multiplicative relations among the arguments. Moreover, a non-zero $\Qbar$-linear combination of logarithms of algebraic numbers is never algebraic. Therefore, computing the desired space of relations reduces to finding the lattice of multiplicative relations between non-zero algebraic numbers. 

To find multiplicative relations between $a_1,\dots,a_n \in \Qbar^\times \subset \C$, construct the abstract number field $K = \Q[x]/p(x)$ and elements $k_1,\dots,k_n \in K$ such that $K$ has an embedding that maps $(k_1,\dots,k_n)$ to $(a_1,\dots,a_n)$. A multiplicative relation between $a_i$'s hold if and only if the same relation holds between the $k_i$'s in $K$. Find a finite set of places $S$ of $K$ such that $S$ contains all Archimedean places and the $k_i$'s are invertible outside of $S$. Now use Drichlet's $S$-unit theorem, or rather its proof using logarithmic embeddings~\cite[Theorem~3.12]{Narkiewicz2004}, together with a day-and-night algorithm analogous to the one given in~\S\ref{sec:computing_AJ_kernel} to find the lattice of all multiplicative relations between the $k_i$'s.
\hfill\qed

\newpage
\section{Supersaturation}\label{sec:supsat}

\node{}\label{node:supsat} In this section, we will describe an algorithm that effectively supersaturates a given Jacobian motive $J_{C,E}^D$. By \emph{effective supersaturation} we mean that the algorithm takes a Jacobian motive $M = J_{C,E}^D$ and outputs the following: 
\begin{enumerate}
  \item A supersaturated motive $M^\ss = J_{C,\psi'}^{\chi'}$ represented as a push-pull Jacobian motive. 
  \item Lattice morphisms that realize $M$ as a push-pull of $M^\ss$. 
  \item An explicit isogeny direct sum decomposition of $M^{\ss} \sim M_B \oplus M^{\sat}$ into the direct sum of a Baker motive and a \emph{saturated} push-pull Jacobian motive $M^{\sat}=J_{C,\psi}^\chi$.
  \item The action of $\End_\Q(J_C)$ on $\H_1(M^\sat)$.
\end{enumerate}

\node{}\label{node:rels_of_M_from_Mss} Note that period relations of $M$ can be determined effectively from those of $M^\ss$ by~\eqref{prop:pushpull_relations}. The isogeny splitting gives an effective isomorphism $\H_1(M^{\ss}) \simeq \H_1(M_B) \oplus \H_1(M^{\sat})$. Since $M_B$ is Baker and $M^{\sat}$ is saturated, we can compute the period relations of $M^{\ss}$ from those of $M_B$ and $M^{\sat}$~\eqref{cor:rels_of_Baker_and_sat}. The period relations of $M_B$ can be computed directly~\eqref{lem:baker_period_relations}. The period relations of $M^{\sat}$ can be computed since it is saturated and we know the endomorphism action~\eqref{thm:expected_relations}. Therefore, once the supersaturation is complete, we can compute the period relations of $M$. %

\node{Remark.} In principle, during the algorithm, we change $C$ by adding disjoint copies of $\p$ to $C$, see~\eqref{node:make_lattice_injective}. However, this does not change $J_C$ and we will make this operation without changing the name of $C$.

\node{}\label{node:basis_of_correspondences} Throughout, we fix a representation of $\ce \coleq \End_\Z(J_C)$. This can be computed~\eqref{thm:endo} and involves correspondences $\cc_1,\dots,\cc_\rho$ over $C$ such that $\mu_1=[\cc_{1,*}],\dots,\mu_\rho = [\cc_{\rho,*}]$ is a $\Z$-linear basis for $\ce$. We will take $\cc_1 = \Delta_C$ to be the diagonal, inducing the identity on divisors. We also have the integer tensor $(a_{ij}^k)$ representing the product, $\mu_i\cdot\mu_j = \sum_k a_{ij}^k\mu_k$. 

\subsection{Supersaturation of Jacobian motives of the second kind}\label{sec:supsat_second_kind}

\node{} We begin by giving an algorithm to supersaturate a Jacobian motive $J_{C,D} = [L_D \to J_C]$ of the second kind. Let $X' = \ce \otimes L_D$, and define the supersaturated motive as $M^{\ss} = [X' \to J_C]$ via the map $\mu \otimes \xi \mapsto \mu([\xi])$. The motive $J_{C,D}$ is the pullback of $M^{\ss}$ via $L_D \overset{1 \otimes \id}{\too} \ce \otimes L_D$. The motive $M^{\ss}$ coincides with the one considered in~\cite[\S 15]{HW}. We note that $M^{\ss}$ can be explicitly represented as a pullback Jacobian motive, and while the induced action on the points of $J_C$ is straightforward, the corresponding action of $\ce$ on $\H_1(M^{\ss})$ expressed via the push-pull representation is significantly more involved.

\node{}\label{node:chi_primes} Determine~\eqref{lem:Dprime} finite subsets $D', D'' \subset C(\Qbar)$ such that our chosen basis of correspondences~\eqref{node:basis_of_correspondences} induce correspondences from $(C,D)$ to $(C,D')$ and from $(C,D')$ to $(C,D'')$. Explicitly, we obtain lattice maps:
\begin{equation}\label{eq:role_of_Dprimes}
\chi' \colon \ce \otimes L_D \to L_{D'} : \mu_i \otimes \xi \mapsto \cc_{i,*}(\xi), \quad \text{and} \quad
\chi'' \colon \ce \otimes L_{D'} \to L_{D''} : \mu_i \otimes \xi \mapsto \cc_{i,*}(\xi).
\end{equation}
Because $\cc_1=\Delta_C$, we have $D \subset D' \subset D''$. Let $\iota_{D',D''} \colon L_{D'} \toi L_{D''}$ be the natural inclusion map.

\node{Pullback Jacobian representation.}\label{node:pushpull_second_kind} By construction, the map $X' \to J_C$ factors through $\chi' \colon X' \to L_{D'}$. Consequently, we represent the motive $M^{\ss}$ as the pullback Jacobian motive $\chi^{\prime,*} J_{C,D'}$.

\node{Action on points.} The induced action on the points of $M^{\ss}$ is clear. The algebra $\ce$ acts on the lattice $\ce \otimes L_D$ through left multiplication, which is effective. For the points of $J_C(\Qbar)$, the action of each element of $\ce$ is explicitly realized via the chosen basis of correspondences.

\node{Action on cohomology.} The goal of the rest of the subsection is to provide an algorithm explicitly computing the action of $\ce$ on $\H_1(M^{\ss}) = \chi^{\prime,*}\H_1(J_{C,D'})$.

\node{} To simplify exposition, we first ensure the lattice maps $\chi'$ and $\chi''$ are injective~\eqref{node:make_lattice_injective}, and scale these maps appropriately to avoid torsion divisors~\eqref{node:ignore_torsion_in_chi}.

\node{}\label{node:split_K_X} Define $K \subset \ce \otimes L_D$ as the pullback of $T_{D'}$ (equivalently, of $K_{D'}$). Since the $\ce$-action maps torsion points of $J_C$ to torsion points, the sublattice $K$ must be $\ce$-invariant. Use the semisimplicity of $\ce$, and use~\eqref{rem:semisimple}, to identify a $\ce$-invariant complementary sublattice $X \subset X'$ such that $K \oplus X$ has full rank in $X'$. We denote by $\chi \colon X \toi L_{D'}$ the restriction of $\chi'$ to $X$, noting this restriction remains injective.

\node{} Through the inclusion $K \oplus X \toi X'$, the original motive $M^{\ss}$ explicitly decomposes, up to isogeny, into the direct sum $[K \to 0] \oplus J_{C,\chi}$, where $J_{C,\chi} = [X \to J_C]$, using~\eqref{prop:explicit_isogeny_with_chi}. Since the $\ce$-action on $[K \to 0]$ and its corresponding $\H_1$ is explicit, the $\ce$-action on $\H_1(M^{\ss}) = \H_1([K \to 0]) \oplus \H_1(J_{C,\chi})$ can be deduced from that of $\H_1(J_{C,\chi})$, which remains to be determined.

\node{}\label{node:commute_markings_1}
Our next goal is to realize the $\ce$‑action on $J_{C,\chi}$ by \emph{correspondences} from $(C,D')$ to $(C,D'')$. However, the obvious diagram of lattices 
\begin{equation}\label{eq:does_not_commute}
  \begin{tikzcd}[row sep=1.2em,column sep=3em]
      \ce \otimes X \ar[r,"\text{multiply}"] \ar[d,"{\mathrm{id}\,\otimes\,\chi}"'] 
      & X \ar[d,"{\iota_{D',D''} \circ \chi}"] \\
      \ce \otimes L_{D'} \ar[r,"\chi''"'] 
        & L_{D''}
  \end{tikzcd}
\end{equation}
\emph{fails to commute}.  Indeed, for basis elements \( \mu_i,\mu_j\in\ce \text{ and } \xi\in X \toi \ce \otimes L_{D'}\) the upper path is induced by
\begin{equation}
  \mu_i\otimes\mu_j\otimes\xi \;\xmapsto{}\; (\mu_i\mu_j)\otimes\xi \;=\; \sum_k a_{ij}^k\,\mu_k\otimes\xi \;\xmapsto{\;\chi''\;} \sum_k a_{ij}^k\,\cc_k\!\cdot\!\xi,
\end{equation}
whereas the lower path is induced by \( \mu_i\otimes \mu_j \otimes \xi \mapsto \cc_i\!\cdot\!\cc_j\!\cdot\!\xi. \)
These two need not agree but, as the $\ce$‑action on $J_C$ is well‑defined, the difference \( \cc_i\!\cdot\!\cc_j\!\cdot\!\xi - \sum_k a_{ij}^k\,\cc_k\!\cdot\!\xi \) is \emph{linearly equivalent to $0$} in $L_{D''}$.

\node{}\label{node:commute_markings_2}
Recall $\chi \colon X \toi L_{D'} \toi L_{D''}$ avoids both $K_{D'}$ and $K_{D''}$ (as well as torsion divisors). Use~\eqref{lem:explicit_isogeny_second_kind} to explicitly isogeny split $J_{C,D''} \sim [K_{D''} \to 0] \oplus J_{C,\aleph_{D''}}$ where $\aleph_{D''} \colon N_{D''} \to L_{D''}$ is defined as in~\eqref{node:rep_quot}. We may scale $\chi$ to assume $X \to L_{D''}$ factors through $n_{D''} L_{D''}$ and hence $X \to J_C$ factors through $\aleph_{D''}$. Since $N_{D''}$ is a quotient of $n_{D''} L_{D''}$ by classes linearly equivalent to $0$, the diagram~\eqref{eq:does_not_commute} \emph{will commute} when $L_{D''}$ is replaced by $N_{D''}$ (and $L_{D'}$ is replaced by a multiple). 

\node{} In particular, fixing the diagram~\eqref{eq:does_not_commute} as in~\eqref{node:commute_markings_2} and including the Jacobian components, we get the following diagram (defined up to isogeny), which commutes at the level of points:
\begin{equation}\label{eq:final_comm_diag}
  \begin{tikzcd}[row sep=1.2em, column sep=4em]
    \ce \otimes {[X \to J_C]} \ar[r] \ar[d] 
      & {[X \to J_C]} \ar[d] \\
    \ce \otimes {[L_{D'} \to J_C]} \ar[r] 
      & {[N_{D''} \to J_C]}.
  \end{tikzcd}
\end{equation}
We will now realize the bottom arrow using correspondences.

\node{}\label{eq:coprod_C} Let $\widetilde{C} = \coprod_{\rho} C$ denote the disjoint union of $\rho$ copies of $C$, indexed by the basis elements $\mu_i \in \ce_{\Z}$. Define $\widetilde{D}' \subset \widetilde{C}(\Qbar)$ as the disjoint union of $\rho$ copies of $D'$. Our choice of basis for $\ce_{\Z}$ identifies the bottom left motive in~\eqref{eq:final_comm_diag} with a Jacobian motive: 
\begin{equation}
  \ce_{\Z}\otimes [L_{D'} \to J_C] \simeq [L_{\widetilde{D}'} \to J_{\widetilde{C}}].
\end{equation} 
To describe the bottom map~\eqref{eq:final_comm_diag} first consider the correspondence $\cc \coleq \coprod \cc_i$ which induces a correspondence from $(\widetilde{C},\widetilde{D}')$ to $(C,D'')$ and hence an explicit map $\cc_{*}\colon J_{\widetilde{C},\widetilde{D}'} \to J_{C,D''}$. %

\node{} We modify diagram~\eqref{eq:final_comm_diag} to take into account the correspondence to $J_{C,D''}$ followed by the explicit surjection coming from the splitting from Lemma~\eqref{lem:explicit_isogeny_second_kind}:
\begin{equation}\label{eq:final_comm_diag_split}
  \begin{tikzcd}[row sep=1.4em, column sep=4em]
    \ce \otimes {[X \to J_C]} \ar[r] \ar[d,hook] & {[X \to J_C]} \ar[rd,hook]\\
    \ce \otimes {L_{D'} \to J_C} = {[L_{\widetilde{D}'} \to J_{\widetilde{C}}]} \ar[r,"{\cc_{*}}"] & {[L_{D''} \to J_C]} \ar[r,two heads] & {[N_{D''} \to J_C]}.
  \end{tikzcd}
\end{equation}
With this arrangement, each map in the diagram is explicitly realized either by correspondences, lattice pullback, or Baker splitting, making the induced maps on $\H_1$'s effective. The commutativity of the diagram follows from that of~\eqref{eq:final_comm_diag}.

\node{} We recall that the diagram does not commute at $J_{C,D''}$ but only at $J_{C,\aleph_{D''}}$. Nevertheless, since $X \toi L_{D'}$ and $X \toi N_{D''}$ are injective, the downward arrows in~\eqref{eq:final_comm_diag_split} induce injections $\ce\otimes\H_1(J_{C,\chi}) \toi \H_1(J_{\widetilde{C},\widetilde{D}'})$ and $\H_1(J_{C,\chi}) \toi \H_1(J_{C,\aleph_{D''}})$. Therefore, we can compute the action $\ce\otimes\H_1(J_{C,\chi}) \to \H_1(J_{C,\chi})$.

\subsection{Supersaturation of Jacobian motives of the third kind}\label{sec:supsat_third_kind}

\node{} This time we give an algorithm to supersaturate a Jacobian motive $J_C^D = [L_D \to J_C]^\vee$ of the third kind. 

\node{Remark.} On could imagine supersaturating $J_{C,D}$ and then passing to the duals, however, we need to be able to do both on the ``same side'' to handle the general case. Moreover, we have not given an algorithm to compute the duality pairing $\H^1_\AdR(C\setminus{}D) \otimes \H^1_\AdR(C,D) \to \C$, therefore we do not use duality to transfer period relations.

\node{} Let $X' = L_D \otimes \ce$ and define the supersaturated motive as $M^{\ss} = [X' \to J_C]^\vee$ where the map is $\xi \otimes \mu \mapsto [\xi]\cdot \mu$. The map $L_D \overset{\id \otimes 1}{\too} \to L_D \otimes \ce$ realizes $M$ as a pushout of $M^{\ss}$. We use the right action on the lattice so that when we dualize we get the left action of $\ce$ on the Jacobian.

\node{}\label{node:third_type_Dprimes} Determine~\eqref{lem:Dprime} finite subsets $D', D'' \subset C(\Qbar)$ such that our chosen basis of correspondences~\eqref{node:basis_of_correspondences} induce correspondences from $C\setminus{}D$ to $C\setminus{}D'$ and from $C\setminus{}D'$ to $C\setminus{}D''$. Note that we have induced maps $\chi' \colon L_D \otimes \ce \to L_{D'}$ and $\chi'' \colon L_{D'} \otimes \ce \to L_{D''}$ defined as in~\eqref{eq:role_of_Dprimes}. 

\node{Pushout Jacobian representation.}\label{node:pushpull_third_kind} Once again, the map $X' \to J_C$ factors through $\chi' \colon X' \to L_{D'}$ and, hence, the motive $M^{\ss}$ is represented as the pushout motive $\chi'_* J_C^{D'}$.

\node{The action.} It remains to compute the action of $\ce$ on the points and on $\H_1$ of $M^{\ss}$. We do so following the dual construction for the second kind motives. We will give a quick summary on how build the dual diagram.

\node{} Once again, we may assume $\chi'$ and $\chi''$ are injective and avoid torsion divisors. We determine~\eqref{prop:explicit_isogeny_with_chi} an explicit isogeny splitting $M^{\ss} \sim \G_m^K \oplus J_C^\chi$ for a $\ce$-equivariant isogeny splitting $X' \sim K \oplus X$ where $K$ is the kernel of $X' \to J_C(\Qbar)$ and $\xi = \xi'|_X$. Since the identification of the points $M^{\ss}(\Qbar)\otimes \Q$ with the isogeny splitting is effective as well the splitting $\H_1(M^{\ss}) \simeq \H_1(\G_m^K) \oplus J_C^\chi$, we reduce the problem to identifying the $\ce$-action on each component separately. The $\ce$-action on $\G_m^K$ is obvious, so we reduce to identifying the $\ce$ action on $J_C^\chi$.

\node{Remark.} The right action $K \otimes \ce \to K$ is induced from multiplication on the right by the identification $K \subset L_D \otimes \ce$. Dualizing, we get a map $\G_m^K \to \ce^\vee \otimes \G_m^K$ so that tensoring with $\ce$ and then contracting $\ce \otimes \ce^\vee \to \Z$ gives the desired action $\ce \otimes \G_m^K \to \G_m^K$. 

\node{Dual maps.} The map $J_{C,D'} \otimes \ce \to J_{C,D''} : \xi \otimes \mu_i \mapsto \xi \cdot \cc_i = \cc_i^*\xi$ dualizes to give $J_C^{D''} \to \ce^\vee \otimes J_C^{D'} \simeq (J_C^{D'})^{\oplus \rho}$ which at the level of points is given by
\begin{equation}
  \xi \mapsto \left( \mu_i \mapsto \cc_{i,*}\xi \right) = (\cc_{1,*}\xi, \dots, \cc_{\rho,*}\xi).
\end{equation}
This can be seen by identifying $\ce \otimes J_{C,D'}$ with the Jacobian of a curve and describing the morphism $J_{C,D'} \otimes \ce \to J_{C,D''}$ by a correspondence as in~\eqref{eq:coprod_C} and observing that taking duals uses the transpose action of the correspondence.

\node{Dualize the action.}\label{node:commute_comarkings} Repeat the argument in the previous section to arrive at the \emph{commutative} diagram~\eqref{eq:final_comm_diag_split} on the side of lattices but using the right action of $\ce$. Taking duals we get the following \emph{commutative} diagram,
\begin{equation}
  \begin{tikzcd}[row sep=1.4em, column sep=4em]
    \ce^\vee \otimes {J_C^\chi} & {J_C^\chi} \ar[l] \\
    \ce^\vee \otimes {J_C^{D'}} \ar[u,two heads] & J_C^{D''} \ar[l]  & J_C^{\aleph_{D''}} \ar[l, hook] \ar[ul,two heads].
  \end{tikzcd}
\end{equation}
We now tensor by $\ce$ and then contract with $\ce \otimes \ce^\vee \to \Z$.

\node{} Once again $\widetilde{C} = \coprod_{\rho} C$ is the disjoint union of $\rho$ copies of $C$. This time, let $\widetilde{D}'' \subset \widetilde{C}(\Qbar)$ be the disjoint union of $\rho$ copies of $D''$. Using our basis $\ce \simeq \Z^\rho$, we identify $\ce \otimes J_C^{D''} \simeq (J_C^{D''})^{\oplus\rho}$ with $J_{\widetilde{C}}^{\widetilde{D}''}$.

\node{} We now arrive at the commutative diagram, which makes the action on points and on $\H_1$ explicit:
\begin{equation}\label{eq:final_diagram_third_kind}
  \begin{tikzcd}[row sep=1.4em, column sep=4em]
    & \ce \otimes J_C^\chi \ar[r] & J_C^\chi \\
    \ce \otimes J_C^{\aleph_{D''}} \ar[ur,two heads] \ar[r,hook] & \ce \otimes J_C^{D''} = J_{\widetilde{C}}^{\widetilde{D}''}  \ar[r,"{\cc_{*}}"] & J_C^{D'} \ar[u,two heads]
  \end{tikzcd}
\end{equation}

\node{Effectivity.} The inclusion map $J_C^{\aleph_{D''}} \toi J_C^{D''}$ is deduced from the explicit splitting~\eqref{lem:explicit_isogeny_third_type}. The vertical maps are induced by lattice pushouts and the map $\cc_*$ is induced by a correspondence. Therefore, all maps except the top horizontal map are explicit. Lifting to $J_C^{\aleph_{D''}}$ then going through $\cc_*$ and down to $J_C^\chi$ gives effective maps on points $\ce \otimes J_C^\chi \to J_C^\chi$ and on homology $\ce \otimes \H_1(J_C^\chi) \to \H_1(J_C^\chi)$.

\node{Action on points.}\label{rem:lifting_points} For future use, we spell out the action on the points. 
Since $J_C^\chi$ is realized as a pushout of $J_C^{D'}$, a point $p \in J_C^\chi$ is represented by a divisor $\xi \in \Div^0(C\setminus{}D')$ so that $p=[\xi]_\chi$. Move~\eqref{lem:move} the divisor $\xi$ away from $D''$ in its $D'$-linear equivalence class, getting $\xi' \in [\xi]_\chi$ which lifts to $[\xi']_{D''} \in J_C^{D''}$. 
The pushout map gives a point $[\xi']_{\aleph_{D''}} \in J_C^{\aleph_{D''}}$, this is a lift of $p$ along the diagonal map in~\eqref{eq:final_diagram_third_kind}. Now, use the Baker splitting~\eqref{node:including_aleph} to apply the (non-trivial!) inclusion $J_C^{\aleph_{D''}} \to J_C^{D''}$ to the point $[\xi']_{\aleph_{D''}}$ and get $[\xi'']_{D''}$. This operation moves $\xi'$ by linear equivalence but changes its $D''$-linear equivalence class. The point $\mu_i \otimes [\xi'']_{D''}$ gets mapped to $[\cc_{i,*}\cdot \xi'']_{D'}$ via the bottom horizontal map and then down to $[\cc_{i,*}\cdot\xi'']_{\chi}$. Hence $\mu_i \otimes [\xi]_\chi \mapsto [\cc_{i,*}\cdot \xi'']_{\chi}$, which extends linearly to $\ce \otimes J_C^\chi(\Qbar) \to J_C^\chi(\Qbar)$. 

\node{Further notes on the action on points.}\label{rem:further_lifting_points} Note that $\xi'' = \xi' + \divv(u)$ for some $u\in \kappa(C)_{D''}$ and thus $\cc_{i,*}\cdot\xi'' = \cc_{i,*}\cdot\xi' + \divv(v)$ where $\divv(v) = \cc_{i,*}\cdot\divv(u)$ and $v \in \kappa(C)_{D'}$. Later, we will want to represent the point $\mu_i \cdot [\xi]_\chi$ by a divisor in the image of $J_C^{\aleph_{D'}} \toi J_C^{D'}$. Applying the Baker splitting again, compute $\xi''' = \xi'' + \divv(w)$ for $w\in \kappa(C)_{D'}$ so that $\mu_i\cdot[\xi]_\chi = [\xi''']_\chi$ and $[\xi''']_{D'}$ is in the image of $J_C^{\aleph_{D'}}$.

Summarizing, we have $\mu_i \cdot [\xi]_\chi = [\cc_i\cdot \xi' + \divv(f)]$ where $f \in \kappa(C)_{D'}$ depends only on $i$ and $\xi'$. If $\xi$ is already distinct from $D''$ then we can take $\xi = \xi'$. If, $\xi$ already defines a point in the image of $J_C^{\aleph_{D''}}$ then we may take $\divv(f)=0$, $\mu_i \cdot [\xi]_\chi = [\cc_i \cdot \xi]$. In the general case, we may choose $f$ so that $[\cc_i\cdot \xi' + \divv(f)]_{D'}$ is in the image of $J_C^{\aleph_{D'}}$.

\subsection{The general case}

\node{} We now begin the supersaturation process for the motive $M = J_{C,D}^E$. Construct finite subsets $D' \subset D'' \subset C(\overline{\Q})$ as in~\eqref{node:third_type_Dprimes}. Extend this chain by constructing a further set $D'''$ such that the $\cc_i$'s induce correspondances from $(C \setminus D''')$ to $(C \setminus D'')$.

\node{Moving $E$ away from $D'''$.} We may assume without loss of generality that $E \cap D''' = \emptyset$. Indeed, we can move each divisor in $L_E$ away from $D'''$ within its $D$-linear equivalence class, as follows. For each $\xi \in L_E$, we can compute a divisor $\xi' \in [\xi]_D$ disjoint from $D'''$ using~\eqref{lem:move}. Move a basis of divisors of $L_E$ away from $D'''$ and let $\widetilde{E} \subset C \setminus D''$ denote the joint support of the new divisors. This induces a map $L_E \to L_{\widetilde{E}}$ compatible with the Abel--Jacobi maps to $J_C^D$. Hence $J_{C,E}^D$ is realized as a pullback of $J_{C,\widetilde{E}}^D$. Supersaturating the latter supersaturates the former, and we may therefore replace $E$ by $\widetilde{E}$ from now on.

\node{} Supersaturate the underlying motive of the third kind $J_C^D$ following \S\ref{sec:supsat_second_kind}. We continue with the notation there: $\chi, \chi', \chi'', K, X$. This process realizes $J_C^D$, up to isogeny, as a pushout of $\G_m^K \oplus J_C^\chi$. Note that $J_C^\chi$ inherits a $\ce$-action compatible with the usual action on $J_C$.

\node{} Recall that, after a minor adjustment~\eqref{node:make_lattice_injective}, we assume $J_C^{D'}$ dominates the supersaturation of $J_C^D$. Thus, the marking $L_E \to J_C^D$ lifts to $J_C^{D'}$ (and, in fact, to $J_C^{D'''}$), yielding a marking $L_E \to J_C^\chi \oplus \G_m^K$. We thereby realize $J_{C,E}^D$ as a push-pull of $[L_E \to J_C^\chi \oplus \G_m^K]$.

\node{} Use the two projection maps to construct the motive $[L_E \to J_C^\chi] \oplus [L_E \to \G_m^K]$. The diagonal map $L_E \to L_E^2$ realizes $[L_E \to J_C^\chi \oplus \G_m^K]$ as a pullback of $[L_E \to J_C^\chi] \oplus [L_E \to \G_m^K]$. We set aside the Baker component $[L_E \to \G_m^K]$. What remains is to supersaturate $J_{C,E}^\chi = [L_E \to J_C^\chi]$ to complete the supersaturation of $M$.

\node{} Let $M^{\ss} = [\ce \otimes L_E \to J_C^\chi]$, with map $\mu \otimes \xi \mapsto \mu \cdot [\xi]_\chi$. This motive $M^{\ss}$ is clearly supersaturated. Moreover, $J_{C,E}^{\chi}$ is the pullback of $M^{\ss}$ along $L_E \xrightarrow{1 \otimes \id} \ce \otimes L_E$. The remainder of this section aims to make explicit the $\ce$-action on points and (co)homology of $M^{\ss}$.

\node{Jacobian representation.}\label{node:jac_rep_general_case} Recall the pointwise action~\eqref{rem:lifting_points} on $J_C^\chi$. We will use $D'''$ and $D''$ instead of $D''$ and $D'$ for the construction. Using the basis of correspondences $\mu_i = [\cc_i]$ for $\ce$ and a basis $\xi_j$ of $L_E$, we have:
\begin{equation}
\mu_i \cdot [\xi_j]_\chi = [\cc_i \cdot \xi_j + \divv(f_{ij})]_\chi,
\end{equation}
where $f_{ij} \in \kappa(C)_{D''}$, since $\xi_j$ need not be moved away from $D'''$~\eqref{rem:further_lifting_points}. Moreover, by choosing $f_{ij}$'s appropriately, we may assume $\xi_{ij} \coleq \cc_i \cdot \xi_j + \divv(f_{ij})$ defines a $D''$-linear equivalence class in the image of $J_C^{\aleph_{D''}} \toi J_C^{D''}$.

Let $E' \subset C(\overline{\Q}) \setminus D''$ be a finite set supporting the divisors $\xi_{ij}$ for all $i,j$. Then, we define $\psi' \colon \ce \otimes L_E \to L_{E'} : \mu_i \otimes \xi_j \mapsto \xi_{ij}$ to realize $M^{\ss}$ as the pullback $J_{C,\psi'}^\chi$ of $J_{C,E'}^\chi$.

\node{Further decomposing the supersaturated motive.}\label{node:psi_general_case} Once again apply~\eqref{node:make_lattice_injective} to ensure $\psi'$ is injective and avoids torsion~\eqref{node:ignore_torsion_in_chi}. We now follow the argument in~\eqref{node:split_K_X} to decompose $\ce \otimes L_E$ up to isogeny as $K_{\psi'} \oplus Y$ into two $\ce$-invariant pieces: $K_{\psi'}$ is the kernel of $\ce \otimes L_E \overset{\psi'}{\to} L_{E'} \to J_C(\overline{\Q})$ and $Y$ is a $\ce$-invariant complement. Let $\psi \colon Y \to L_{E'}$ be the restriction of $\psi'$ to $Y$.

\node{} Observe that we can realize $J_{C,\psi'}^\chi$ effectively as a pushout of the direct sum $[K_{\psi'} \to \G_m^X] \oplus J_{C,\psi}^\chi$, where the pushout is by the diagonal $X \to X \oplus X$. Putting aside the Baker component, we will consider $J_{C,\psi}^\chi$, which is not only supersaturated but also reduced and hence saturated.

\node{The action figure.} It remains to make the $\ce$-action on $J_{C,\psi}^\chi$ explicit. At this point, we only need the action on $\H_1$ but we will also describe the action on points. This will be done via Figure~\ref{fig:action}. The arrow~\labelarr{action} represents the $\ce$-action we are after and we will show how to navigate the explicit arrows in this diagram to evaluate~\circledref{action} by a composition of explicit maps. We recall the argument in~\eqref{node:commute_markings_1} to emphasize the need to follow the outer arrows in the diagram to achieve commutativity.

\begin{figure}[ht]
  \centering
\begin{tikzcd}
	& {\mathcal{E} \otimes [Y \to J_C^\chi]} && {[Y \to J_C^\chi]} \\
	{\mathcal{E} \otimes [Y \to J_C^{\aleph_{D''}}]} & {\mathcal{E} \otimes [Y \to J_C^{D''}]} && {[L_{E''} \to J_C^\chi]} & {[N^\chi_{E''} \to J_C^\chi]} \\
	& {\mathcal{E} \otimes [L_{E'} \to J_C^{D''}]} && {[L_{E''} \to J_C^{D'}]}
  \arrow[from=1-2, to=1-4, "\circledref{action}"]
	\arrow[hook', from=1-4, to=2-4, "\circledref{rhs:inject}"]
	\arrow[hook', from=2-5, to=2-4,   "\circledref{rhs:pullback}"']
	\arrow[shift right=3, two heads, from=2-4, to=2-5,  "\circledref{rhs:split}"']
	\arrow[hook', from=1-4, to=2-5, "\circledref{rhs:embed}"]
	\arrow[two heads, from=3-4, to=2-4, "\circledref{rhs:surject}"]
	\arrow[two heads, from=2-2, to=1-2,  "\circledref{lhs:coinject}"]
	\arrow[hook, from=2-1, to=2-2, "\circledref{lhs:split}"]
	\arrow[two heads, from=2-1, to=1-2, "\circledref{lhs:coembed}"]
	\arrow[shift left=3, two heads, from=2-2, to=2-1, "\circledref{lhs:pushout}"]
	\arrow[hook, from=2-2, to=3-2, "\circledref{lhs:cosurject}"]
	\arrow[from=3-2, to=3-4, "\circledref{correspond}"]
\end{tikzcd}
\caption{Making the action on a supersaturated motive effective.}
  \label{fig:action}
\end{figure}

\node{The right hand side.} We will describe the arrows appearing on the right hand side of Figure~\ref{fig:action}. Apply~\eqref{lem:Dprime} to construct a finite set $E'' \subset C(\Qbar)\setminus{}D'$ which supports the image of $E' \subset C(\Qbar)\setminus{}D''$ under our basis of correspondences. The inclusion $Y \toi L_{E'} \toi L_{E''}$ induces the map~\labelarr{rhs:inject}. 

\node{} We adopt the notation from~\eqref{node:chi_torsion} to decompose $L_{E''}$ into $K_{E''}^\chi \oplus N_{E''}^\chi$ up to isogeny. The arrow~\labelarr{rhs:pullback} is trivially explicit, since it is the canonical inclusion defined by the pullback induced by the (non-canonical) splitting map $\aleph_{E''}^\chi \colon N_{E''}^\chi \to L_{E''}$. The projection~\labelarr{rhs:split} is not trivial but it is made explicit by~\eqref{lem:split_kchipsi}. The composition of~\circledref{rhs:inject} and~\circledref{rhs:split} gives the map~\labelarr{rhs:embed}, which is an injection because the composition $Y \toi L_{E'} \toi L_{E''} \tos N_{E''}^\chi$ is injective by choice of $Y$.

\node{} The surjection~\labelarr{rhs:surject} is simply the pushout representation of $J_{C,E''}^\chi$ and is trivially explicit. 

\node{The left hand side.} By our choice of basis for $\ce$ we make the identification $\ce = \Z^\rho$ and therefore, for any motive $A$, make the identification $\ce \otimes A = A^{\oplus \rho}$. For a map $f\colon A \to B$, we can define a map $\ce \otimes A \to \ce \otimes B$ by simply acting as identity on $\ce$, i.e., taking $f^{\oplus \rho}$. All the maps appearing solely on the left hand side of Figure~\ref{fig:action} are of this form and we will describe the respective maps on the components.

\node{} The map~\labelarr{lhs:coinject} is the pushout corresponding to $X \toi L_{D'} \toi L_{D''}$. The pushout $J_C^{D''} \tos J_C^{\aleph_{D''}}$ is explicitly split~\eqref{lem:explicit_isogeny_third_type} and the map~\labelarr{lhs:split} is the splitting inclusion. Recall that the map~\labelarr{lhs:coembed} is surjective, since $\chi$ avoids torsion. Our construction of $\psi' \colon \ce \otimes L_E \to L_{E'}$ in~\eqref{node:jac_rep_general_case} is such that we lifted our marking $Y\to J_C^{D''}$ through~\circledref{lhs:coembed}. 

\node{} The map~\labelarr{lhs:pushout} is used indirectly: in constructing $\psi'$, we first lifted the marking to $J_C^{D''}$ naively, mapped down to $J_C^{\aleph_{D''}}$, and then took the image marking in $J_C^{D''}$. The naive lift (through~\circledref{lhs:coinject}) and the one we constructed (through~\circledref{lhs:coembed}) are different markings on $J_C^{D''}$. 

\node{} The inclusion~\labelarr{lhs:cosurject} is the natural pullback map via $\psi\colon Y \to L_{E'}$. 

\node{Applying the correspondence.}\label{node:apply_correspondence} Let $\widetilde{C} = \coprod_{i=1}^\rho C$. Define $\widetilde{D}''$, $\widetilde{E}'$ by taking a copy of $D''$ and $E'$ per component of $\widetilde{C}$. Then, using our basis for $\ce$, we identify $\ce \otimes [L_{E'} \to J_C^{D''}]$ with $J_{\widetilde{C},\widetilde{E}'}^{\widetilde{D}''}$. Define the correspondence $\widetilde{\cc}$ to be ``disjoint union'' of the correspondences $\widetilde{\cc}_i$ for $i=1,\dots,\rho$ from $(\widetilde{C}\setminus{}\widetilde{D}'',\widetilde{E}')$ to $(C\setminus{}D',E'')$. This correspondence defines the map~\labelarr{correspond}, and is the backbone for the diagram in Figure~\ref{fig:action}.

\node{Commutativity.}\label{node:commutativity} The \emph{path of correspondence} is the map attained by composing \circledref{lhs:split}, \circledref{lhs:cosurject}, \circledref{correspond}, \circledref{rhs:surject}, \circledref{rhs:split}. The \emph{path of action} is the map attained by composing \circledref{lhs:coembed}, \circledref{action}, \circledref{rhs:embed}. We claim that these two maps commute. Commutativity at the level of semi-abelian groups was established in~\eqref{node:commute_comarkings}. At the level of markings, we appeal to the commutativity of the diagram on the semi-abelian varieties---that is, on divisor classes---and proceed as in~\eqref{node:commute_markings_1} and~\eqref{node:commute_markings_2}.

We now expand on the proof of commutativity at the level of markings. Observe that $\mu_i \otimes y \in \ce \otimes Y$ is mapped to $\cc_i \cdot \psi(y) \in L_{E''}$ via the path of correspondence, whereas it is mapped to $\psi(\mu_i \cdot y) \in L_{E'} \subset L_{E''}$ via the path of action. Note, however, that $\mu_i \cdot [\psi(y)]\chi = [\cc_i \cdot \psi(y)]\chi$, since $\psi(y)$ defines a divisor lying in the image of $J_C^{\aleph_{D''}}$~\eqref{rem:further_lifting_points}. On the other hand, the $\ce$-action on $Y$ and $J_C^\chi$ are compatible by construction, so we also have $\mu_i \cdot [\psi(y)]\chi = [\psi(\mu_i \cdot y)]\chi$. It follows that $\cc_i \cdot \psi(y) - \psi(\mu_i \cdot y) \in K_{D''}^{\chi}$, and hence the difference vanishes in $N_{D''}^\chi$, establishing the desired commutativity.

\node{Effectivity.} 
We wish to evaluate the map~\circledref{action} on points and on $\H_1$. We established the commutativity~\eqref{node:commutativity} of the path of action, the path involving~\circledref{action}, with the path of correspondence, the path involving~\circledref{correspond}. By lifting via~\circledref{lhs:coembed} we may follow the path of correspondence and then, identifying $J_{C,\psi}^\chi$ with its image via~\circledref{rhs:embed}, we get the desired action. We claim that this operation is effective at the level of points and on $\H_1$. Most of the maps are trivial at the level of representatives because they represent pushout or pullback maps. The non-trivial maps are the splitting maps~\circledref{lhs:split} and~\circledref{rhs:split}. These splittings have been made explicit at the level of points and at the level of $\H_1$ when splitting Baker components. The map~\circledref{correspond} is effective as it is induced by a correspondence~\eqref{node:apply_correspondence}.

\node{} This completes the description and proof of the algorithm that effectively supersaturates $J_{C,D}^E$.

\newpage
\section{Classifying first-order autonomous differential equations}\label{sec:autonomous_ode}

\node{} Let $P(x,y) \in \Qbar[x,y]$ be an irreducible polynomial. The autonomous first-order differential equation $P(u(t), u'(t)) = 0$ falls into one of four classes~\eqref{node:classes}. As noted in~\cite[p.1655]{NPT22}, no algorithm is known to classify such equations in full generality. We resolve this problem by providing a classification algorithm based on our computational model for the mixed Hodge structure of curves.

\node{} Let $C$ be the smooth proper curve with function field $\Qbar[x,y]/P$ and consider the rational differential $\omega = \dd x / y$ on $C$. For initial conditions $(u(0), u'(0))$ avoiding singularities and vertical tangents of the affine plane curve $Z(P) \subset \C^2$, the solution $u(t)$ to $P(u,u')=0$ satisfies, and is determined by, the equation
\begin{equation}
t = \int_{u(0)}^{u(t)} \frac{\dd x}{y(x)}, \quad \text{with } P(x,y(x)) = 0,\ \abs{t} \ll 1.
\end{equation}
For this reason, the differential equation $P(u,u')=0$ is classified by the isomorphism type of the pair $(C,\dd x/y)$.

\node{Classification.}\label{node:classes} Let $C/\Qbar$ be a smooth proper irreducible curve and $\omega$ a nonzero rational differential on $C$. Then $(C,\omega)$ falls into one of the following classes:
\begin{enumerate}
\item $\omega$ is exact,
\item $\omega = h^*(\dd z/z)$ for a non-constant $h \colon C \to \p$,
\item $\omega = h^*\eta$ for a non-constant $h \colon C \to C'$ with $C'$ an elliptic curve and $\eta$ regular on $C'$,
\item none of the above.
\end{enumerate}
The first two cases are effectively decidable using Riemann--Roch. Deciding between the last two cases is the problem.

\node{} We can significantly refine the classification above. Observe that the second and third clauses detect whether $\omega$ is the pullback of a translation invariant form from a semi-abelian variety of dimension $1$. Such forms are always of the third kind, see~\eqref{node:pullback_thirdkind}, and differentials of the third kind are unique in their cohomology class. The first clause can thus be viewed as asserting that the cohomology class of $\omega$ is a pullback from a semi-abelian variety of dimension $0$. We generalise the classification at the level of cohomology classes, and specialise to pullbacks of translation invariant regular forms in~\eqref{node:detect_pullback}.

\node{} Let $D \subset C(\Qbar)$ be the residual polar locus of $\omega$, and let $p \in C(\Qbar) \setminus D$. The Abel--Jacobi map $\AJ_p \colon C \setminus D \to J_C^D$ identifies $\H^1_\AdR(C\setminus{}D,E)$ with $\H^1_\AdR(J_C^D)$. 

\node{Refined classification.} There exists a minimal (up to isogeny) semi-abelian quotient $f\colon J_C^D \tos G$ and a translation invariant differential $\eta$ on $G$ such that $\omega \equiv \AJ_p^* f^* \eta$ in $\H^1(C \setminus D)$. Let $a$ and $b$ be the dimensions of the toric and abelian parts of $G$, respectively. We define the \emph{refined type} of $\omega$ to be the pair $(a,b)$. 

\node{Comparing the classifications.} The original classification corresponds to the cases $(0,0)$, $(1,0)$ with $\omega$ of the third kind, $(0,1)$ with $\omega$ of the third kind, and $a + b > 1$. Detecting whether $\omega$ is of the third kind is comparatively easy, see~\eqref{node:detect_third_kind} below, therefore we consider the refined type $(a,b)$ to be the essential piece of information.

\node{Refined classification algorithm.}\label{node:classification_algo}
Given $(C,\omega)$ for a smooth proper irreducible curve $C/\Qbar$ and a rational differential $\omega$ on $C$, compute its refined type as follows:
\begin{enumerate}
\item Compute a representation for $\H_1(C \setminus D)$.
\item Compute coordinates for $[\omega] \in \H^1_\AdR(C \setminus D)$.
\item \label{item:kernel} Compute the kernel $K \coleq [\omega]^{\perp} \subset \H_1^\B(C\setminus{}D)$ under the period pairing. 
\item Compute $W_{-2}K \coleq K \cap W_{-2} \H_1^\B(C \setminus D)$.
\item Return the tuple $(\dim_\Q W_{-2}K,\, \dim_\Q K/W_{-2}K)$.
\end{enumerate}

\node{} The required subroutines are described in the main algorithm~\eqref{alg:main_algorithm}. In step~\eqref{item:kernel} of~\eqref{node:classification_algo}, compute the kernel of the period pairing and then restrict to $\H_1^\B(C\setminus{}D) \otimes \Qbar \cdot [\omega]$. We now prove correctness.

\node{Theorem.} Algorithm~\eqref{node:classification_algo} terminates and returns the correct refined type of $\omega$.
\itnode{Proof.}
Termination of each step follows from the subroutines cited in~\eqref{alg:main_algorithm}. For correctness, recall the identification $\H_1(J_C^D) = \H_1(C \setminus D)$. By Wüstholz's analytic subgroup theorem~\cite{Wustholz1989}, the orthogonal $K$ is the Betti realization of a semi-abelian subvariety $G' \subset J_C^D$, and $[\omega]$ is a pullback from the quotient $G \coleq J_C^D/G'$. As the kernel is the largest possible, the quotient is the smallest possible. Hence $G$ is the minimal isogeny quotient of $J_C^D$ from which $[\omega]$ can be pulled-back. The toric and abelian dimension of $G$ are $(a,b) = (\dim W_{-2}K, \dim K/W_{-2}K)$.
\hfill\qed

\node{Pullback of translation invariant regular differentials.}\label{node:pullback_thirdkind} Differentials of the third kind on $C$ correspond exactly to pullbacks of translation invariant regular differentials from semi-abelian varieties. To see this, use that $J_C^D$ is universal with respect to maps from $C\setminus{}D$ to semi-abelian varieties and that translation invariant forms on $J_C^D$ pullback to differentials of the third kind with residual poles at $D$.

\node{Detecting differentials of the third kind.}\label{node:detect_third_kind} We can check whether a given differential $\omega$ is of the third kind or equivalent to one. A differential of the third kind will give an element in $F^1\H^1_{\AdR}$ and a class in $F^1\H^1_\AdR$ will contain precisely one differential of the third kind. Let $D$ be the residual polar locus of $\omega$. If $[\omega] \notin F^1\H^1_\AdR(C\setminus{}D)$ then $\omega$ is not equivalent to a differential of the third kind. If $[\omega] \in F^1\H^1_\AdR(C\setminus{}D)$ then our reduction algorithm, \S\ref{sec:reduction_algos}, 
returns a representative $\eta \in [\omega]$ of the third kind and an exact form $\dd h$ such that $\omega = \eta + \dd h$.

\node{}\label{node:detect_pullback} Consequently, we can detect whether $\omega$ is equal to, or equivalent modulo exact forms to, the pullback of a translation invariant differential form from a semi-abelian variety. Combined with the refined type of $\omega$, this is a strict refinement of the original classification~\eqref{node:classes}.

\newpage
\printbibliography
\end{document}